\documentclass{article}
\usepackage{arxiv}
\newcommand{\Ra}{\Rightarrow}
\newcommand{\p}{\mathfrak{p}}

\newcommand{\R}{\mathbb{R}}
\newcommand{\N}{\mathbb{N}}
\usepackage{latexsym, amssymb, amsfonts,amsmath}
\usepackage{amsthm}
\usepackage{graphicx}
\usepackage[all,knot,arc,import,poly]{xy}
\usepackage[utf8]{inputenc} % allow utf-8 input
\usepackage[T1]{fontenc}    % use 8-bit T1 fonts
\usepackage{hyperref}       % hyperlinks
\usepackage{url}            % simple URL typesetting
\usepackage{booktabs}       % professional-quality tables
\usepackage{amsfonts}       % blackboard math symbols
\usepackage{nicefrac}       % compact symbols for 1/2, etc.
\usepackage{microtype}      % microtypography
\usepackage{lipsum}
\usepackage{url}
\usepackage{latexsym, amssymb, amsfonts,amsmath}
\usepackage{amsthm}
\newtheorem{theorem}{Theorem}[section]
\newtheorem{lemma}[theorem]{Lemma}
\newtheorem{fact}[theorem]{Fact}
\newtheorem{definition}[theorem]{Definition}
\newtheorem{proposition}[theorem]{Proposition}

\newtheorem{remark}[theorem]{Remark}
\newtheorem{corollary}[theorem]{Corollary}
\newcommand{\ci}{\mathcal{C}^\infty}
\title{Separation Theorems in Smooth Commutative Algebra and Applications}
\author{ Jean C.~Berni\thanks{\texttt{www.ime.usp.br/$\sim$jeancb/}} \\
  Department of Mathematics\\
  Institute of Mathematics and Statistics\\
  University of S\~{a}o Paulo\\
  \texttt{jeancb@ime.usp.br} \\
   \And
 Hugo L.~Mariano\thanks{\texttt{www.ime.usp.br/$\sim$hugomar/}} \\
  Department of Mathematics\\
  Institute of Mathematics and Statistics\\
  University of S\~{a}o Paulo\\
  \texttt{hugomar@ime.usp.br} \\
}

\begin{document}

\maketitle

\begin{abstract}
In this paper we state and prove  \textit{ad hoc} ``\textbf{Separation Theorems}'' of the so-called Smooth Commutative Algebra, the Commutative Algebra of $\mathcal{C}^{\infty}-$rings. These results are formally similar to the ones we find in (ordinary) Commutative Algebra. However, their proof is not so straightforward, since it depends on the introduction of the concept of ``smooth saturation''. As an application of these theorems we present an interesting result that sheds light on the connections between the smooth Zariski spectrum  and the real smooth spectrum of a $\mathcal{C}^{\infty}-$ring.
\end{abstract}

% keywords can be removed

\section*{Introduction} \label{intro-sec}

It is a well known fact that a smooth  manifold $M$ - which is a geometrical object - can be encoded by an algebraic entity, its $\mathbb{R}$-algebra of smooth real  functions, $\mathcal{C}^\infty(M,\mathbb{R})$, since there is a canonical bijection (evaluation) $M \cong {\rm Hom}_{\rm \bf \R-Alg}(\ci(M,\R) , \R)$. This identification works even at the level of morphisms $\ci(M,M') \cong$  ${\rm Hom}_{\rm \bf \R-Alg}(\ci(M',\R),$ $\ci(M,\R))$. Moreover, geometric constructions over a manifold $M$, as its tangent bundle, $TM$, remain algebraically represented: $TM \cong {\rm Hom}_{\rm \bf \R-Alg}(\ci(M,\R),$ $\R[x]/(x^2))$. 

The set $\ci(M,\R)$  supports a far richer structure than just an $\R$-algebra: it interprets not only the real polynomial functions  but all smooth real functions $\R^n \to \R$, $n \in \mathbb{N}$. Moreover, this extended  interpretation also satisfies all the compositional identities that hold between the smooth real functions. Thus, $\ci(M,\R)$ is a natural instance of the algebraic structure called \underline{$\ci$-ring}.

A more systematic  algebraic study of these rings of smooth functions was carried out in the  mid 1960's and early 1970's. It was not until the decades of 1970's and 1980's that  a  study of the  abstract (algebraic) theory  of $\ci$-rings was made,  mainly  in order  to construct - out of the ideas of F. W. Lawvere -  topos models for Synthetic Differential Geometry (\cite{Kock},  \cite{rings1}, \cite{rings2}, \cite{MSIA}); In \cite{Dubuc}, the first steps towards an ``algebraic geometry of $\ci$-rings'' were taken, through the definition of the $\ci$-schemes.

The interest in $\ci$-rings gained strength in recent years (\cite{Bunge}, \cite{Spivak}, \cite{Borisov}, \cite{tese}, \cite{Joyce}). In \cite{Joyce},  D. Joyce presents the foundations of a version of Algebraic Geometry in which the role of rings  is replaced by $\mathcal{C}^{\infty}-$rings, with the goal to apply  these new notions and results to  Differential Geometry: this can be considered as a part of a larger (and ambitious)  program, also pursued by 
D. I. Spivak (see \cite{Spivak}), of extending/transferring Jacob Lurie's program of Derived Algebraic Geometry to  Derived Differential Geometry. 
 In \cite{Borisov}, Kremnizer and Borisov give a detailed account of six notions of radicals of an ideal of a $\mathcal{C}^{\infty}-$ring, among which we find the $\infty-$radical of an ideal of a $\mathcal{C}^{\infty}-$ring.  Here we focus on this concept, that we call `` the $\mathcal{C}^{\infty}-$radical'' of an ideal: this concept first appeared in \cite{rings1}, and it is proper to the theory of $\mathcal{C}^{\infty}-$rings, carrying some differences with respect to the usual notion of radical in ordinary Commutative Algebra - which only makes use of powers of elements.

The difference between the notions of ``radical'' and ``$\mathcal{C}^{\infty}-$radical'' ideals brings us, alone,  a whole new study of some important concepts, such as $\mathcal{C}^{\infty}-$re\-duced $\mathcal{C}^{\infty}-$rings, the ``smooth Zariski spectrum'' , this smooth version of the Zariski spectrum presents some crucial differences when compared to the ordinary Zariski spectrum, in its topological as well as in its functorial and sheaf-theoretic features (see \cite{ShvLog}, \cite{BM1}).

%In order to illustrate this difference, we refer the reader to an example of an ideal of a $\mathcal{C}^{\infty}-$ring which is radical in the ordinary sense (since it is prime), but  not $\mathcal{C}^{\infty}-$radical (see \cite{rings2}). \\

The main reference on category theory is \cite{CWM}. For Commutative algebra we refer the reader to \cite{AM}.

\textbf{Overview of the Paper:}  \\

In the {\bf Section \ref{prelim-sect}}
we provide some preliminaries on the main subject on which we build this version of Commutative Algebra, the category $\mathcal{C}^{\infty}-$rings, presenting some definitions and some of their fundamental constructions (the main definitions can be found in \cite{MSIA}, and a more detailed account of  $\mathcal{C}^{\infty}-$rings can be found, for example, in \cite{UACR}).

In the  {\bf Section \ref{smoothrings-sect}} we present o$\mathcal{C}^{\infty}$ parallels of the ordinary Commutative Algebraic construction of the ring of fractions, radicals and the concept of saturation, In order to prove the existence of the $\mathcal{C}^{\infty}-$ring of fractions (\textbf{Theorem \ref{conc}}) we make use of the notion of the $\mathcal{C}^{\infty}-$ring of $\mathcal{C}^{\infty}-$polynomials, studied in \cite{UACR}, and then we define the concept of ``smooth saturation'' (\textbf{Definition \ref{satlisa}} of \textbf{Section \ref{saty}}), pointing some of its relationships with the ordinary concept of saturation (\textbf{Theorem \ref{lara}}). We state and prove various results about this concept.

We start by recalling the ordinary Commutative Algebraic construction of the ring of fractions, radical of an ideal and the concept of saturation, motivating the introduction of their $\mathcal{C}^{\infty}$ parallels. In order to prove the existence of the $\mathcal{C}^{\infty}-$ring of fractions we make use of the notion of the $\mathcal{C}^{\infty}-$ring of $\mathcal{C}^{\infty}-$polynomials (see \cite{UACR}). We study the concept of ``$\mathcal{C}^{\infty}$-radical'' introduced in \cite{rings1}, relating it to concept of ``smooth saturation'', introduced in \cite{tese}, pointing some of its relationships with the ordinary concept of saturation in Commutative Algebra. We present many results about these three concepts. In \textbf{Section \ref{classes-sec}}  we present distinguished classes of $\mathcal{C}^{\infty}-$rings (i.e., $\mathcal{C}^{\infty}-$rings which satisfy some further axioms), as $\mathcal{C}^{\infty}-$fields, $\mathcal{C}^{\infty}-$domains and local $\mathcal{C}^{\infty}-$\-rings. We prove some results connecting the filter of closed subsets of $\mathbb{R}^n$ and the set of $\mathcal{C}^{\infty}-$radical ideals of $\mathcal{C}^{\infty}(\mathbb{R}^n)$ (in fact, a Galois connection [\textbf{Proposition \ref{pseu}}]), as well as various results about them - including the fact that the $\mathcal{C}^{\infty}-$radical of any ideal is again an ideal (\textbf{Proposition \ref{49}}). The authors know of no such proof in the current literature. From the preceding results on smooth saturation and on smooth radical ideals, we present in {\bf Section \ref{septheo-sec}} a similar version of the Separation Theorems (\textbf{Theorem \ref{TS}}) one finds in ordinary Commutative Algebra.

In \textbf{Section \ref{order-sec}} we give a survey of some order-theoretic aspects of $\mathcal{C}^{\infty}-$rings, defining fundamental concepts as the ``real $\mathcal{C}^{\infty}-$spectrum'' of a $\mathcal{C}^{\infty}-$ring and its topology (the ``Harrison smooth topology''), proving that every $\mathcal{C}^{\infty}-$ring is semi-real (\textbf{Proposition \ref{doria}}) and as an application of the Separation Theorems we prove an important result which establishes a spectral bijection from the real $\mathcal{C}^{\infty}-$spectrum of a $\mathcal{C}^{\infty}-$ring to its smooth Zariski spectrum (\textbf{Theorem \ref{exv}}).

\section{Preliminaries} \label{prelim-sect}

We provide here the main preliminary notions on   $\mathcal{C}^{\infty}-$rings, with respect to their universal algebra (cf. \cite{UACR}).

In order to formulate and study the concept of $\mathcal{C}^{\infty}-$ring, we use a first order language, $\mathcal{L}$, with a denumerable set of variables (${\rm Var}(\mathcal{L}) = \{ x_1, x_2, \cdots, x_n, \cdots\}$), whose nonlogical symbols are the symbols of $\mathcal{C}^{\infty}-$func\-tions from $\mathbb{R}^m$ to $\mathbb{R}^n$, with $m,n \in \mathbb{N}$, \textit{i.e.}, the non-logical symbols consist only of function symbols, described as follows:

For each $n \in \mathbb{N}$, the $n-$ary \textbf{function symbols} of the set $\mathcal{C}^{\infty}(\mathbb{R}^n, \mathbb{R})$, \textit{i.e.}, $\mathcal{F}_{(n)} = \{ f^{(n)} | f \in \mathcal{C}^{\infty}(\mathbb{R}^n, \mathbb{R})\}$. Thus, the set of function symbols of our language is given by:
      $$\mathcal{F} = \bigcup_{n \in \mathbb{N}} \mathcal{F}_{(n)} = \bigcup_{n \in \mathbb{N}} \mathcal{C}^{\infty}(\mathbb{R}^n).$$
Note that our set of constants is identified with the set of all $0-$ary function symbols, \textit{i.e.}, $\mathcal{C} = \mathcal{F}_{(0)} = \mathcal{C}^{\infty}(\mathbb{R}^0) \cong \mathcal{C}^{\infty}(\{ *\})$.

The terms of this language are defined, in the usual way, as the smallest set which comprises the individual variables, constant symbols and $n-$ary function symbols followed by $n$ terms ($n \in \mathbb{N}$).

Functorially, a (set-theoretic) $\mathcal{C}^{\infty}-$ring is a finite product preserving functor from the category $\mathcal{C}^{\infty}$, whose objects are of the form $\mathbb{R}^n$, $n \in \mathbb{N}$, and whose morphisms are the smooth functions between them, \textit{i.e.}, a finite product preserving functor:
$$A: \mathcal{C}^{\infty} \rightarrow {\rm \bf Set}$$

Apart from the functorial definition and the ``first-order language'' definition we just gave, there are many equivalent descriptions. We focus, first, on the universal-algebraic description of a $\mathcal{C}^{\infty}-$ring in ${\rm \bf Set}$, given in the following:

\begin{definition}\label{cabala} A \textbf{$\mathcal{C}^{\infty}-$structure} on a set $A$ is a pair $ \mathfrak{A} =(A,\Phi)$, where:

$$\begin{array}{cccc}
\Phi: & \bigcup_{n \in \mathbb{N}} \mathcal{C}^{\infty}(\mathbb{R}^n, \mathbb{R})& \rightarrow & \bigcup_{n \in \mathbb{N}} {\rm Func}\,(A^n; A)\\
      & (f: \mathbb{R}^n \stackrel{\mathcal{C}^{\infty}}{\to} \mathbb{R}) & \mapsto & \Phi(f) := (f^{A}: A^n \to A)
\end{array},$$

that is, $\Phi$ interprets the \textbf{symbols}\footnote{here considered simply as syntactic symbols rather than functions.} of all smooth real functions of $n$ variables as $n-$ary function symbols on $A$.
\end{definition}

We call a $\mathcal{C}^{\infty}-$struture $\mathfrak{A} = (A, \Phi)$ a \textbf{$\mathcal{C}^{\infty}-$ring} whenever it preserves projections and all equations between smooth functions. More precisely, we have the following:

\begin{definition}\label{CravoeCanela}Let $\mathfrak{A}=(A,\Phi)$ be a $\mathcal{C}^{\infty}-$structure. We say that $\mathfrak{A}$ (or, when there is no danger of confusion, $A$) is a \textbf{$\mathcal{C}^{\infty}-$ring} if the following is true:

$\bullet$ Given any $n,k \in \mathbb{N}$ and any projection $p_k: \mathbb{R}^n \to \mathbb{R}$, we have:

$$\mathfrak{A} \models (\forall x_1)\cdots (\forall x_n)(p_k(x_1, \cdots, x_n)=x_k).$$

$\bullet$ For every $f, g_1, \cdots g_n \in \mathcal{C}^{\infty}(\mathbb{R}^m, \mathbb{R})$ with $m,n \in \mathbb{N}$, and every $h \in \mathcal{C}^{\infty}(\mathbb{R}^n, \mathbb{R})$ such that $f = h \circ (g_1, \cdots, g_n)$, one has:
$$\mathfrak{A} \models (\forall x_1)\cdots (\forall x_m)(f(x_1, \cdots, x_m)=h(g(x_1, \cdots, x_m), \cdots, g_n(x_1, \cdots, x_m))).$$
\end{definition}

\begin{definition}Let $(A, \Phi)$ and $(B,\Psi)$ be two $\mathcal{C}^{\infty}-$rings. A function $\varphi: A \to B$ is called a \textbf{morphism of $\mathcal{C}^{\infty}-$rings} or \textbf{$\mathcal{C}^{\infty}$-homomorphism} if for any $n \in \mathbb{N}$ and any $f: \mathbb{R}^n \stackrel{\mathcal{C}^{\infty}}{\to} \mathbb{R}$, one has $\Psi(f)\circ \varphi^{(n)} = \varphi \circ \Phi(f)$, where $\varphi^{(n)} = (\varphi, \cdots, \varphi): A^n \to B^n$.
\end{definition}

\begin{remark} {\bf (on universal algebraic constructions)}
It is not difficult to see that $\mathcal{C}^{\infty}-$structures, together with their morphisms (which we call $\mathcal{C}^{\infty}-$morphisms) compose a category, that we denote by $\mathcal{C}^{\infty}{\rm \bf Str}$, and that $\mathcal{C}^{\infty}-$rings, together with all the $\mathcal{C}^{\infty}-$morphisms between $\mathcal{C}^{\infty}-$rings (which we call $\mathcal{C}^{\infty}-$homomorphisms) compose a full subcategory of $\mathcal{C}^{\infty}{\rm \bf Rng}$. In particular, since $\mathcal{C}^{\infty}{\rm \bf Rng}$ is a ``variety of algebras'', i.e. it is a class of $\mathcal{C}^{\infty}-$structures which satisfies a given set of equations, (or equivalently,  by \textbf{Birkhoff's HSP Theorem}) it is closed under substructures, homomorphic images and products. Moreover:\\
$\bullet$ $\mathcal{C}^{\infty}{\rm \bf Rng}$ is a concrete category and the forgetful functor, $ U :  \mathcal{C}^{\infty}{\rm \bf Rng}  \to {\rm \bf Set}$ creates directed inductive colimits;\\
$\bullet$ Each set $X$ freely generates a $\mathcal{C}^{\infty}$-ring. In particular, the free $\mathcal{C}^{\infty}$-ring on $n$ generators is $\mathcal{C}^{\infty}(\mathbb{R}^n)$, $n \in \mathbb{N}$;\\
$\bullet$ Every $\mathcal{C}^{\infty}-$ring is the homomorphic image of some free $\mathcal{C}^{\infty}-$ring determined by some set, being isomorphic to the quotient of a free $\mathcal{C}^{\infty}-$ring by some congruence;\\
$\bullet$ The congruences of $\mathcal{C}^{\infty}-$rings are classified by their ``ring-theoretical'' ideals;\\
$\bullet$ In $\mathcal{C}^{\infty}{\rm \bf Rng}$ one defines ``the $\mathcal{C}^{\infty}-$coproduct'' between two $\mathcal{C}^{\infty}-$rings $\mathfrak{A} = (A,\Phi)$ and $\mathfrak{B}=(B,\Psi)$, denoted by $A \otimes_{\infty} B$;\\
$\bullet$ Using free $\mathcal{C}^{\infty}-$rings and the $\mathcal{C}^{\infty}-$coproduct, one gets the ``$\mathcal{C}^{\infty}-$ring of polynomials'' on any set $S$ of variables with coefficients in $A$, given by $A\{ x_s \mid s \in S\} = A\otimes_{\infty}\mathcal{C}^{\infty}(\mathbb{R}^S)$.
\end{remark}

\section{Smooth Rings of Fractions}\label{smoothrings-sect}

\hspace{0.5cm} We begin by giving a description of the fundamental concept of ``smooth ring of fractions'', presenting a slight modification of the axioms given in \cite{rings1}. In order to show that the $\mathcal{C}^{\infty}-$ring of fractions exists in the category of $\mathcal{C}^{\infty}-$rings, we use the $\mathcal{C}^{\infty}-$ring of $\mathcal{C}^{\infty}-$polynomials. The definitions and results we state here may be found with more detail in \cite{BMTSCA}.\\

We turn to the discussion of how to obtain the ring of fractions of a $\mathcal{C}^{\infty}-$ring $A$ with respect to some of its subsets, $S \subseteq U(A, \Phi)$.\\

For any commutative unital ring $R$, one proves that $S^{-1}R$ has the following universal property (cf. \textbf{Proposition 3.1} of \cite{AM}):

\begin{proposition}\label{Nab}Given a ring homomorphism $g: R \to B$  such that $(\forall s \in S)(g(s) \in B^{\times})$, there is a unique ring homomorphism $\widetilde{g}: S^{-1}R \to B$ such that the following triangle commutes:

$$\xymatrixcolsep{5pc}\xymatrix{
A \ar[r]^{\eta_S} \ar[dr]_{g} & S^{-1}R \ar[d]^{\exists ! \widetilde{g}}\\
   & B}$$
\end{proposition}

In order to extend the notion of the ring of fractions to the category $\mathcal{C}^{\infty}{\rm \bf Rng}$ we make use of the universal property described in \textbf{Proposition \ref{Nab}}, as we can see in the following:\\

\begin{definition}\label{Alem}Let $A$ be a $\mathcal{C}^{\infty}-$ring and $S \subseteq A$ one of its subsets. The \index{$\mathcal{C}^{\infty}-$ring of fractions}$\mathcal{C}^{\infty}-$\textbf{ring of fractions} of $A$ with respect to $S$ is a $\mathcal{C}^{\infty}-$ring together with a $\mathcal{C}^{\infty}-$homo\-morphism $\eta_S: A \to A\{ S^{-1}\}$ with the following properties:
\begin{itemize}
  \item[(1)]{$(\forall s \in S)(\eta_S(s) \in (A\{ S^{-1}\})^{\times})$}
  \item[(2)]{If $\varphi: A \to B$ is any $\mathcal{C}^{\infty}-$homomorphism such that for every $s \in S$ we have $\varphi(s) \in B^{\times}$, then there is a unique $\mathcal{C}^{\infty}-$homomorphism $\widetilde{\varphi}: A\{ S^{-1}\} \to B$ such that the following triangle commutes:
      $$\xymatrixcolsep{5pc}\xymatrix{
      A \ar[r]^{\eta_S} \ar[rd]^{\varphi} & A\{ S^{-1}\} \ar[d]^{\widetilde{\varphi}}\\
        & B}$$}
\end{itemize}

By this universal property, the $\mathcal{C}^{\infty}-$ring of fractions is unique, up to (unique) isomorphisms.
\end{definition}

Now we prove the existence of such a $\mathcal{C}^{\infty}-$ring of fractions by constructing it. Recall that we can only make use of the constructions available within the category $\mathcal{C}^{\infty}{\rm \bf Rng}$, such as the free $\mathcal{C}^{\infty}-$ring on a set of generators, their coproduct, their quotients and others described in \cite{UACR}.\\

Recall that given any set $S$, the  $\mathcal{C}^{\infty}-$ring of ``smooth polynomials'' in the set $S$ of variables over a $\mathcal{C}^{\infty}-$ring $A$ is obtained as follows:\\

We consider $\mathcal{C}^{\infty}(\mathbb{R}^S)$,  the free $\mathcal{C}^{\infty}-$ring on the set $S$ of generators, together with its canonical map, $\jmath_S: S \to \mathcal{C}^{\infty}(\mathbb{R}^S)$. If we denote by:

$$\xymatrixcolsep{5pc}\xymatrix{
A \ar[dr]^{\iota_A} & \\
   & A \otimes_{\infty} \mathcal{C}^{\infty}(\mathbb{R}^S)\\
\mathcal{C}^{\infty}(\mathbb{R}^S) \ar[ur]_{\iota_{\mathcal{C}^{\infty}(\mathbb{R}^S)}}
}$$

the coproduct of $A$ and $\mathcal{C}^{\infty}(\mathbb{R}^S)$, we define:
$$x_s := \iota_{\mathcal{C}^{\infty}(\mathbb{R}^S)}(\jmath_S(s)).$$

By definition, we have:

$$A\{ x_s | s \in S \} := A \otimes_{\infty} \mathcal{C}^{\infty}(\mathbb{R}^S).$$

so we can consider the quotient:

$$\dfrac{A\{x_s | s \in S \}}{\langle \{ x_s \cdot \iota_A(s) - 1 | s \in S \}\rangle},$$

where $\langle \{ x_s \cdot \iota_A(s) - 1 | s \in S\}\rangle$ is the ideal of $A$ generated by $\{ x_s \cdot \iota_A(s) - 1 | s \in S\}$. \\

In these settings we can formulate the following:\\

\begin{theorem}\label{conc}Let $A$ be a $\mathcal{C}^{\infty}-$ring and $S \subseteq A$ be any of its subsets. The \textbf{$\mathcal{C}^{\infty}-$ring of fractions of $A$ with respect to $S$} is given concretely by the $\mathcal{C}^{\infty}-$ring:
$$A\{ S^{-1}\}:= \dfrac{A\{ x_s | s \in S \}}{\langle \{x_s \cdot \iota_A(s) - 1 | s \in S \}\rangle}$$
together with the $\mathcal{C}^{\infty}-$homomorphism:
$$\eta_S:= q \circ \iota_A : A \to A\{ S^{-1}\},$$
where $q: A\{ x_s | s \in S\} \to \dfrac{A\{x_s | s \in S \}}{\langle \{ x_s \cdot \iota_A(s) - 1 | s \in S\}\rangle}$ is the canonical quotient map and $\iota_A: A \to A\{ x_s | s \in S\} = A \otimes_{\infty} \mathcal{C}^{\infty}(\mathbb{R}^S)$ is the canonical coproduct homomorphism corresponding to $A$.
\end{theorem}
\begin{proof}
It suffices to show that $q \circ \iota_A: A \rightarrow \dfrac{A\{ x_s | s \in S \}}{\langle \{x_s \cdot s - 1 \}\rangle}$ satisfies the universal property described in the \textbf{Definition \ref{Alem}}.\\

For a detailed proof, see p. 4 of \cite{BMTSCA}.

\end{proof}

\begin{remark}As we shall see later on, in \textbf{Proposition \ref{fact}}, the universal map $\eta_S: A \to A\{ S^{-1}\}$ is an epimorphism in $\mathcal{C}^{\infty}{\rm \bf Rng}$, even though it is not always a surjective map.
\end{remark}

\subsection{Smooth Saturation}\label{saty}

A very useful concept for Commutative Algebra was given by A. Grothendieck and J. Dieudonn\'{e} in \cite{EGA}, namely the concept of ``a saturated (multiplicative) set''. For the reader's benefit, we give the following:\\

\begin{definition}Let $A$ be any commutative ring with unity and $S \subseteq A$ any of its subsets. Denoting by $\langle S \rangle$ the multiplicative submonoid of $A$ generated by $S$, we define \textbf{the saturation of} $S$ as follows:
$$S^{{\rm sat}} = \{ a \in A | (\exists d \in A)( a \cdot d \in \langle S \rangle) \}.$$

In other words, the saturation of a set $S$ is the set of all divisors of  elements in $\langle S \rangle$.
\end{definition}

One easily sees that the saturation of a subset $S$ of a commutative ring $A$ is equal to the pre-image of the invertible elements of $A[S^{-1}]$ by the canonical map $\eta_S: A \to A\{ S^{-1}\}$, \textit{i.e.}, $S^{{\rm sat}} = \eta_S^{\dashv}[(A[S^{-1}])^{\times}]$. \\

We are going to use this characterization in order to introduce the concept of ``the \textit{smooth} saturation of a subset $S$ of a $\mathcal{C}^{\infty}-$ring $A$'', that we are going to denote by $S^{\infty-{\rm sat}}$.\\

First we need the following:

\begin{proposition}\label{fati}Let $(A,\Phi)$ be a $\mathcal{C}^{\infty}-$ring and $S \subseteq A$ be any subset. If both $(F, \sigma)$ and $(F', \sigma')$ satisfy the \textbf{Definition \ref{Alem}}, then $\sigma^{\dashv}[F^{\times}] = {\sigma'}^{\dashv}[{F'}^{\times}].$
\end{proposition}
\begin{proof}
See the proof of \textbf{Proposition 2}, p. 19 of \cite{BMTSCA}.
\end{proof}

Now we give the following:

\begin{definition}\label{satlisa}Let $A$ be a $\mathcal{C}^{\infty}-$ring, $S \subseteq A$ and $(F,\sigma)$ be a ring of fractions of $A$ with respect to $S$. The \index{smooth saturation}\textbf{smooth saturation} of $S$ in $A$ is:
$$S^{\infty-{\rm sat}}:= \{ a \in A | \sigma(a) \in F^{\times}\}.$$
\end{definition}

In virtue of the \textbf{Proposition \ref{fati}}, the set $S^{\infty-{\rm sat}}$ does not depend on any particular choice of the representation of the ring of fractions, rather it depends only on $A$ and $S$.\\

\begin{remark}Since for every $s \in S$, $\eta_S(s) \in (A\{ S^{-1}\})^{\times}$, from now on we are going to use the more suggestive \underline{notation}:

$$\left( \forall s \in S\right)\left(\dfrac{1}{\eta_S(s)} \stackrel{\cdot}{=} \eta_S(s)^{-1}\right),$$

and for any $a \in A$ and $s \in S$ we are going to denote:

$$\dfrac{\eta_S(a)}{\eta_S(s)} \stackrel{\cdot}{=} \eta_S(a)\cdot{\eta_S(s)^{-1}}.$$
\end{remark}

\begin{definition}Let $A$ be a $\mathcal{C}^{\infty}-$ring and let $S \subseteq A^{\times}$ be any subset. The \textbf{smooth saturation of $S$} is $S^{\infty-{\rm sat}} = {\eta_S}^{\dashv}[A\{ S^{-1}\}^{\times}]$, where $\eta_S\,: A \to A\{ S^{-1}\}$ is the canonical map of the ring of fractions of $A$ with respect to $S$.
\end{definition}

\begin{remark}\label{bols}Let $A$ be a $\mathcal{C}^{\infty}-$ring, $S \subseteq A$ and consider the forgetful functor:
$$\begin{array}{cccc}
    \mathcal{U}: & \mathcal{C}^{\infty}{\rm \bf Rng} & \to & \mathbb{R}-{\rm \bf Alg} \\
     & A & \mapsto & \mathcal{U}(A)\\
     & A \stackrel{f}{\to} B & \mapsto & \mathcal{U}(A) \stackrel{\mathcal{U}(f)}{\to} \mathcal{U}(B)\\
  \end{array}$$
We have always:
$$S^{{\rm sat}} \subseteq S^{\infty-{\rm sat}}$$
\end{remark}

%%%%%%%%%%%%%%%%%%%%%%%%%%%%%%%%%%%%%%%%%%%%%%%%%%%%%%%

\begin{theorem}\label{lara}Let $A$ be a $\mathcal{C}^{\infty}-$ring, $S \subseteq A$ and consider the forgetful functor:
$$\begin{array}{cccc}
    \mathcal{U}: & \mathcal{C}^{\infty}{\rm \bf Rng} & \to & \mathbb{R}-{\rm \bf Alg} \\
     & A & \mapsto & \mathcal{U}(A)\\
     & A \stackrel{f}{\to} B & \mapsto & \mathcal{U}(A) \stackrel{\mathcal{U}(f)}{\to} \mathcal{U}(B)\\
  \end{array}$$
Since $\eta_S^{\infty}: A \to A\{ S^{-1}\}$ is such that $\eta_S^{\infty}[S] \subseteq (A\{ S^{-1}\})^{\times}$, then $(\eta_S^{\infty})[S] \subseteq (\mathcal{U}(A\{ S^{\-1}\}))^{\times}$, so by the universal property of the ring of fractions:
$$\eta_S : \mathcal{U}(A) \to \mathcal{U}(A)[S^{-1}],$$
there is a unique $\R-$algebras homomorphism, ${\rm Can}: \mathcal{U}(A)[S^{-1}] \to \mathcal{U}(A\{ S^{-1}\})$ such that the following diagram commutes:
$$\xymatrix{
\mathcal{U}(A) \ar[r]^{\eta_S} \ar[dr]_{\mathcal{U}(\eta_S^{\infty})} & \mathcal{U}(A)[S^{-1}] \ar@{.>}[d]^{{\rm Can}}\\
  & \mathcal{U}(A\{ S^{-1}\})
}$$
The following assertions are equivalent:
\begin{itemize}
  \item[(1)]{$S^{{\rm sat}} = S^{\infty-{\rm sat}}$;}
  \item[(2)]{${\rm Can}$ is an isomorphism of $\mathbb{R}-$algebra.}
\end{itemize}
\end{theorem}
\begin{proof}See p. 20 of \cite{BMTSCA}.
%Ad (2) $\to$ (1). Since the above diagram commutes, we have:
%$$S^{\infty-{\rm sat}} = {\eta_S^{\infty}}^{\dashv}[(A\{S^{-1}\})^{\times}] = \eta_S^{\dashv}[{\rm Can}^{\dashv}[(A\{S^{-1}\})^{\times}]] \stackrel{(2)}{=} \eta_S^{\dashv}[(A[S^{-1}])^{\times}] = S^{{\rm sat}}$$

%Ad (1) $\to$ (2): We already know that ${\rm Can}: \mathcal{U}(A)[S^{-1}] \to \mathcal{U}(A\{S^{-1}\})$ is an $\R-$algebras morphism. Note that since $S^{{\rm sat}} = S^{\infty-{\rm sat}}$, the morphism:
%$$\mathcal{U}(\eta_S^{\infty}): \mathcal{U}(A) \to \mathcal{U}(A\{ S^{-1}\})$$
%is such that:
%\begin{itemize}
 % \item[(i)]{$(\forall \varphi' \in \mathcal{U}(A)[S^{-1}])(\exists a \in U(A))(\exists b \in S^{{\rm sat}})(\mathcal{U}(\eta_S^{\infty})(b) \cdot \varphi' = \mathcal{U}(\eta_S^{\infty})(a))$;}
%  \item[(ii)]{$(\forall a \in U(A))(\mathcal{U}(\eta_S^{\infty})(a)=0 \to (\exists \lambda' \in S^{{\rm sat}})(\lambda' \cdot a = 0))$;}
%  \item[(iii)]{$S \subseteq S^{{\rm sat}}$.}
%\end{itemize}
%which are precisely the hypotheses of \textbf{Theorem \ref{37}}, so $\mathcal{U}(\eta_S^{\infty}): \mathcal{U}(A) \to \mathcal{U}(A\{ S^{-1}\})$ is isomorphic to the localization. This fact implies that since ${\rm Can}$ is the only ring homomorphism which makes the diagram commute, ${\rm Can}$ must be the unique isomorphism between $\mathcal{U}(A)[S^{-1}]$ and $\mathcal{U}(A\{ S^{-1}\})$
\end{proof}

In what follows we give some properties relating the inclusion relation among the subsets of a $\mathcal{C}^{\infty}-$ring and their smooth saturations.\\

\begin{proposition}Let $A$ be a $\mathcal{C}^{\infty}-$ring and $T, S \subseteq A$ be any of its subsets. Then:
\begin{itemize}
  \item[(i)]{$A^{\times} \subseteq S^{\infty-{\rm sat}}$}
  \item[(ii)]{$S \subseteq S^{\infty-{\rm sat}}$}
  \item[(iii)]{$S \subseteq T$ implies $S^{\infty-{\rm sat}} \subseteq T^{\infty-{\rm sat}}$}
  \item[(iv)]{$S^{\infty-{\rm sat}} = \langle S \rangle^{\infty-{\rm sat}}$, where $\langle S \rangle$ is the submonoid generated by $S$.}
\end{itemize}
\end{proposition}
\begin{proof}See p. 21 of \cite{BMTSCA}.
\end{proof}

Some necessary and sufficient conditions for the $\mathcal{C}^{\infty}-$homomorphism $\eta_S : A \to A\{ S^{-1}\}$ be a $\mathcal{C}^{\infty}-$iso\-morphism is given below:\\

\begin{proposition}Let $A$ be a $\mathcal{C}^{\infty}-$ring and $S \subseteq A$ any of its subsets. The following assertions are equivalent:
\begin{itemize}
  \item[(i)]{$\eta_S : A \to A\{ S^{-1}\}$ is an isomorphism;}
  \item[(ii)]{$S^{\infty-{\rm sat}} \subseteq A^{\times}$;}
  \item[(iii)]{$S^{\infty-{\rm sat}} = A^{\times}$}
\end{itemize}
\end{proposition}
\begin{proof}

Ad (i) $\to$ (iii): Since $\eta_S$ is an isomorphism, both $\eta_S^{-1}$ and $\eta_S$ preserve the invertible elements, so $A^{\times} = \eta_S^{-1}[(A\{S^{-1}\})^{\times}]=\eta_S^{\dashv}[(A\{ S^{-1}\})^{\times}] = S^{{\rm sat}}$.\\

Ad (ii) $\leftrightarrow$ (iii): Since we always have $A^{\times} \subseteq S^{\infty-{\rm sat}}$, by (ii) we conclude that $A^{\times} = S^{\infty-{\rm sat}}$.

Ad (iii) $\to$ (i): Suppose that $S^{\infty-{\rm sat}} = A^{\times}$. We need to show $\eta_{A^{\times}} : A \to A\{ (A^{\times})^{-1}) \}$ is an isomorphism.\\

First we note that ${\rm id}_A : A \to A$ has the universal property of the ring of fractions of $A$ with respect to $A^{\times}$. Indeed, given any $\psi: A \to B$ such that $\psi[A^{\times}] \subseteq B^{\times}$ (i.e., any $\psi$ which is a $\mathcal{C}^{\infty}-$homomorphism), there exists a unique $\mathcal{C}^{\infty}-$homomorphism from $A$ to $B$, namely $\psi: A \to B$, such that the following diagram commutes:
$$\xymatrix{
A \ar[r]^{{\rm id}_A} \ar[rd]_{\psi} & A \ar[d]_{\psi}\\
           & B
}$$
It follows that $({\rm id}_A\, : A \to A) \cong (\eta_{A^{\times}} : A \to  A\{ (A^{\times})^{-1} \})$, since both satisfy the same universal property. Thus, $\eta_{A^{\times}}$ is the composition of a $\mathcal{C}^{\infty}-$isomorphism with ${\rm id}_A$, hence it is a $\mathcal{C}^{\infty}-$isomorphism.
%It follows that:

%$$\varphi \circ \eta_S = \varphi \circ \eta_{A\{ (A^{\times})^{-1}\}} = {\rm id}_A$$

Moreover, by \textbf{Proposition \ref{Alem}}, we conclude that $A^{\times} = \eta_{A^{\times}}[(A\{ (A^{\times})^{-1}\})^{\times}] = (A^{\times})^{\infty-{\rm sat}}$.
\end{proof}

Next we prove that the smooth saturation of the smooth saturation of a set is again the smooth saturation of this set.\\

\begin{proposition}\label{satsatsat}Let $A$ be a $\mathcal{C}^{\infty}-$ring and $S \subseteq A$ be any of its subsets. Then
$$(S^{\infty-{\rm sat}})^{\infty-{\rm sat}} = S^{\infty-{\rm sat}}$$
\end{proposition}
\begin{proof}
Since $S \subseteq S^{\infty-{\rm sat}}$, there exists a unique morphism $\mu_{SS^{\infty-{\rm sat}}}: A\{ S^{-1}\} \to A\{ (S^{\infty-{\rm sat}})^{-1} \}$ such that the following diagram commutes:
$$\xymatrix{
A \ar[r]^{\eta_S} \ar[dr]_{\eta_{S^{\infty-{\rm sat}}}} & A\{ S^{-1}\} \ar[d]^{\mu_{SS^{\infty-{\rm sat}}}}\\
            & A\{(S^{\infty-{\rm sat}})^{-1} \}
}$$

Now, $\eta_S[S^{\infty-{\rm sat}}] \subseteq (A\{ S^{-1}\})^{\times}$ by the very definition of $S^{\infty-{\rm sat}}$, so, by the universal property of $\eta_{S^{\infty-{\rm sat}}}: A \to A\{ (S^{\infty-{\rm sat}})^{-1} \}$, there exist a unique $\nu : A\{ (S^{\infty-{\rm sat}})^{-1}\} \to A\{ S^{-1}\}$ such that the following diagram commutes:

$$\xymatrix{
A \ar[dr]^{\eta_S} \ar[r]_{\eta_{S^{\infty-{\rm sat}}}} & A\{(S^{\infty-{\rm sat}})^{-1}\}  \ar[d]^{\nu}\\
            & A\{ S^{-1}\}
}$$

We have, then, the following commuting diagrams:

$$\begin{array}{cc}
\xymatrix{
          & A\{(S^{\infty-{\rm sat}})^{-1}\} \ar[d]_{\nu} \ar@/^3pc/[dd]^{{\rm id}_{A\{ (S^{\infty-{\rm sat}})^{-1} \}}} \\
A \ar[ur]^{\eta_{S^{\infty-{\rm sat}}}} \ar[r]^{\eta_S} \ar[dr]_{\eta_{S^{\infty-{\rm sat}}}} & A\{ S^{-1}\} \ar[d]_{\mu_{SS^{\infty-{\rm sat}}}}\\
          & A\{ (S^{\infty-{\rm sat}})^{-1} \}} & \xymatrix{
          & A\{S^{-1}\} \ar[d]_{\mu_{SS^{\infty-{\rm sat}}}} \ar@/^3pc/[dd]^{{\rm id}_{A\{ S^{-1} \}}} \\
A \ar[ur]^{\eta_{S}} \ar[r]^{\eta_{S^{\infty-{\rm sat}}}} \ar[dr]_{\eta_{S}} & A\{ (S^{\infty-{\rm sat}})^{-1}\} \ar[d]_{\nu}\\
          & A\{ S^{-1} \}}
          \end{array}$$

So $(\mu_{SS^{\infty-{\rm sat}}})^{-1} = \nu$, and $\mu_{SS^{\infty-{\rm sat}}}$ is an isomorphism. Hence $A\{ S^{-1}\} \cong A\{ (S^{\infty-{\rm sat}})^{-1}\}$, so by \textbf{Proposition \ref{Alem}},

$$(S^{\infty-{\rm sat}})^{\infty-{\rm sat}} = S^{\infty-{\rm sat}}$$
\end{proof}

\begin{proposition}\label{fact} Let $A$ be a $\mathcal{C}^{\infty}-$ring, $S \subseteq A$ be any of its subsets, and $\eta_S : A \to A\{ S^{-1}\}$ be the canonical morphism of the ring of fractions. If $g,h: A\{ S^{-1}\} \to B$ are two morphisms such that $g \circ \eta_S = h \circ \eta_S$ then $g=h$. In other words, $\eta_S : A \to A\{ S^{-1}\}$ is an epimorphism.
\end{proposition}
\begin{proof}
Note that since $g \circ \eta_S$ is such that $(g \circ \eta_S)[S] \subseteq B^{\times}$, there exists a unique morphism $\widetilde{t}: A\{ S^{-1}\} \to B$ such that $\widetilde{t} \circ \eta_S = g \circ \eta_S$. By hypothesis we have $h \circ \eta_S = g \circ \eta_S$, so  $g$ has the property which determines $\widetilde{t}$. Hence $g = \widetilde{t} = h$\\
\end{proof}

\begin{proposition}\label{zero}Let $A$ be a $\mathcal{C}^{\infty}-$ring and $S,T \subseteq A$ two of its subsets. The following assertions are equivalent:
\begin{itemize}
  \item[(i)]{$S^{\infty-{\rm sat}} \subseteq T^{\infty-{\rm sat}}$}
  \item[(ii)]{There is a unique morphism $\mu: A\{ S^{-1}\} \to A\{ T^{-1}\}$ such that the following diagram commutes:
      $$\xymatrix{
      A \ar[dr]^{\eta_T}\ar[r]^{\eta_S} & A\{ S^{-1}\} \ar[d]^{\mu}\\
          & A\{ T^{-1}\}
      }$$}
\end{itemize}
\end{proposition}
\begin{proof}See p. 24 of \cite{BMTSCA}.

\end{proof}

\begin{corollary}\label{cem}The following assertions are equivalent:
\begin{itemize}
  \item[(i)]{$S^{\infty-{\rm sat}} = T^{\infty-{\rm sat}}$}
  \item[(ii)]{There is an isomorphism $\mu: A\{ S^{-1}\} \to A\{ T^{-1}\}$ such that the following diagram commutes:
      $$\xymatrix{
      A \ar[r]^{\eta_S} \ar[dr]^{\eta_T} & A\{ S^{-1}\} \ar[d]^{\mu}\\
         & A\{ T^{-1}\}
      }$$}
  \item[(iii)]{There is a unique isomorphism $\mu: A\{ S^{-1}\} \to A\{ T^{-1}\}$ such that the following diagram commutes:
      $$\xymatrix{
      A \ar[r]^{\eta_S} \ar[dr]^{\eta_T} & A\{ S^{-1}\} \ar[d]^{\mu}\\
         & A\{ T^{-1}\}
      }$$}
\end{itemize}
\end{corollary}
\begin{proof}See p. 27 of \cite{BMTSCA}.
\end{proof}

\begin{proposition}Let $A$ be a $\mathcal{C}^{\infty}-$ring and $S,T$ two of its subsets such that $S \subseteq T$. The following assertions are equivalent:
\begin{itemize}
  \item[(i)]{$\mu_{ST}: A\{ S^{-1}\} \to A\{ T^{-1}\}$ is an isomorphism;}
  \item[(ii)]{$S \subseteq T \subseteq S^{\infty-{\rm sat}}$;}
  \item[(iii)]{$T^{\infty-{\rm sat}} = S^{\infty-{\rm sat}}$.}
\end{itemize}
\end{proposition}
\begin{proof}
Ad (ii) $\to$ (iii): Since $S \subseteq T$ we have $S^{\infty-{\rm sat}} \subseteq T^{\infty-{\rm sat}}$, and since $T \subseteq S^{\infty-{\rm sat}}$ we have $T^{\infty-{\rm sat}} \subseteq (S^{\infty-{\rm sat}})^{\infty-{\rm sat}} = S^{\infty-{\rm sat}}$. Hence:
$$S^{\infty-{\rm sat}} \subseteq T^{\infty-{\rm sat}} \subseteq S^{\infty-{\rm sat}}$$
and
$$S^{\infty-{\rm sat}} = T^{\infty-{\rm sat}}.$$

Ad (iii) $\to$ (ii): We always have $T \subseteq T^{{\rm sat}}$, and since $T^{{\rm sat}} = S^{{\rm sat}}$, it follows that $T \subseteq S^{{\rm sat}}$.\\

(i) $\leftrightarrow$ (iii) was established in \textbf{Corollary \ref{cem}}.

\end{proof}

\begin{proposition}\label{marina}Let $A$ be a $\mathcal{C}^{\infty}-$ring and $S \subseteq A$. Whenever $\{ S_i\}_{i \in I}$ is a directed system such that:
$$S = \bigcup_{i \in I} S_i$$
we have:
$$S^{\infty-{\rm sat}} = \left( \bigcup_{i \in I} S_i\right)^{\infty-{\rm sat}} = \bigcup_{i \in I} {S_i}^{\infty-{\rm sat}}$$
\end{proposition}
\begin{proof}It is clear that:
$$\bigcup_{i \in I} {S_i}^{\infty-{\rm sat}} \subseteq S^{\infty-{\rm sat}}.$$

In order to prove the other inclusion, we shall use the fact that $A\{S^{-1}\}$ is isomorphic to the vertex of the following directed colimit:

$$\xymatrix{
  & \varinjlim_{i \in I} A\{ {S_i}^{-1}\} &  \\
A\{{S_i}^{-1} \} \ar[rr]^{\alpha_{ij}} \ar@/^/[ur]^{\alpha_i} & & A\{{S_j}^{-1} \} \ar@/_/[ul]_{\alpha_j}
}$$

Note that $\eta_S : A \to A\{ S^{-1}\}$ is such that for any $i \in I$, $\eta_S[S_i] \subseteq \eta_S[S] \subseteq (A\{ S^{-1}\})^{\times}$, so by the universal property of $\eta_{S_i} : A \to A\{ (S_i)^{-1}\}$, there is a unique $\mathcal{C}^{\infty}-$rings homomorphism $\varphi_i: A\{ (S_i)^{-1}\} \to A\{S^{-1}\}$ such that the following triangle commutes:
$$\xymatrix{
A \ar[r]^{\eta_{S_i}} \ar[dr]^{\eta_S} & A\{(S_i)^{-1} \} \ar[d]^{\varphi_i}\\
      & A\{ S^{-1}\}}$$
so
$$\xymatrix{
  & A\{ S^{-1}\} &  \\
A\{{S_i}^{-1} \} \ar[rr]^{\varphi_{ij}} \ar@/^/[ur]^{\varphi_i} & & A\{{S_j}^{-1} \} \ar@/_/[ul]_{\varphi_j}
}$$
commutes for every $i,j \in I$ such that $i \leq j$, since $\eta_{S_i}: A \to A\{ {S_i}^{-1}\}$ is an epimorphism.

Thus, by the universal property of the colimit, there is a unique $\mathcal{C}^{\infty}-$homomorphism:

$$\varphi: \varinjlim_{i \in I} A\{ {S_i}^{-1}\} \rightarrow A\{ S^{-1}\}$$

such that:
 $$(\forall i \in I)(\varphi \circ \alpha_i = \varphi_i).$$

On the other hand, given $s \in S = \bigcup_{i \in I}S_i$ there is $i \in I$ such that $s \in S_i$, so $\eta_{S_i}(s) \in A\{ {S_i}^{-1}\}^{\times}$. Taking $\widetilde{\eta}:= \alpha_i \circ \eta_{S_i}: A \rightarrow \varinjlim_{i \in I} A\{ {S_i}^{-1}\}$, we have:

$$\widetilde{\eta}(s) \in \left( \varinjlim_{i \in I} A\{ {S_i}^{-1}\} \right)^{\times}.$$

By the universal property of $\eta_S: A \rightarrow A\{ S^{-1}\}$, there is a unique $\mathcal{C}^{\infty}-$homo\-morphism $\psi: A\{ S^{-1}\} \rightarrow \varinjlim_{i \in I} A\{ {S_i}^{-1}\}$ such that:

$$\xymatrix{
A \ar[r]^{\eta_{S}} \ar[dr]^{\widetilde{\eta}} & A\{S^{-1} \} \ar[d]^{\psi}\\
      & \varinjlim_{i \in I} A\{ {S_i}^{-1}\}}$$

It is easy to see, by the universal properties involved, that $\varphi$ and $\psi$ are inverse $\mathcal{C}^{\infty}-$isomorphism and that:

$$\widetilde{\eta}^{\dashv}\left[ \left( \varinjlim_{i \in I} A\{ {S_i}^{-1}\}\right)^{\times}\right] = \bigcup_{i \in I} {S_i}^{\infty-{\rm sat}}.$$

Thus $S=\bigcup_{i \in I} {S_i}^{\infty-{\rm sat}}$, as we claimed.
\end{proof}

%%%%%%%%%%%%%%%%%%%%%%%%%%%%%%%%%

%%%%%%%%%%%%%%%%%%%%%%%%%%%%%%%%

%%%%%%%%%%%%%%%%%%%%%%%%%%%%%%%%

%Na versão final colocar isto aqui no capítulo 1:

%There is also the construction of the universal solution to freely adjoining an element $t$ to a $\mathcal{C}^{\infty}-$ring  $A$ now \textit{in  the category $\mathcal{C}^{\infty}{\rm {\bf Rng}}$}, that we shall denote by $A\{ t\}$. The following result gives us an explicit construction.\\

Our next goal is to give a characterization of ring of fractions in $\mathcal{C}^{\infty}{\rm \bf Rng}$ using a similar axiomatization one has in Commutative Algebra. In order to motivate it, we first present some characterizations of rings of fractions in ${\rm \bf CRing}$.

\begin{proposition}
Let $A$ be a commutative ring with unity, $S \subseteq A$ a multiplicative subset, and $\eta: A \to A[S^{-1}]$ its localization. Then:
\begin{itemize}
\item[(i)]{$(\forall \beta \in A[S^{-1}])(\exists c \in S)(\exists d \in A)(\beta \cdot \eta(c) = \eta(d))$;}
\item[(ii)]{$(\forall \beta \in A)(\eta(\beta)=0 \to (\exists c \in S)(c \cdot \beta = 0)$}
\end{itemize}
\end{proposition}
\begin{proof}See p. 31 of \cite{BMTSCA}.
%Ad $(i)$: Given any $\beta \in A[S^{-1}]$ we can write $\beta =\dfrac{a}{s}$ for some $a \in A$ and $s \in S$. Take $c=s$ and $d=a$, so we have:
%$$\beta \cdot \eta(c) = \dfrac{a}{s} \cdot \eta(c) = \dfrac{a}{s}\cdot \dfrac{s}{1} = \dfrac{a}{1} = \eta(a),$$

%so it suffices to take $d = a$.\\

%Ad $(ii)$: Given any $\beta \in A$ such that $\eta(b)=0$ we have, by the very construction of $A[S^{-1}]$, that $(\exists u \in S)(u(\beta\cdot 1 - 0 \cdot 1) = 0)$, so it is enough to take $c=u$ in order to get $c \cdot \beta = u \cdot \beta = 0$.
\end{proof}

Conversely we have:

\begin{proposition}Let $A$ be a commutative ring with unity and $S \subseteq A$ a multiplicative submonoid. If $\varphi : A \to B$ is a ring homomorphism such that $\varphi[A] \subseteq B^{\times}$ and:
\begin{itemize}
\item[(i)]{$(\forall \beta \in B)(\exists c \in S)(\exists d \in A)(\beta \cdot \varphi(c) = \varphi(d))$}
\item[(ii)]{$(\forall \beta \in A)(\varphi(\beta)=0 \to (\exists c \in S)(c \cdot \beta = 0))$}
\end{itemize}
then $B \cong A[S^{-1}]$.
\end{proposition}
\begin{proof}
See p. 32 of \cite{BMTSCA}.
%Since $\varphi: A \to B$ is such that $\varphi[A] \subseteq B$, by the universal property of the localization there is a unique $\psi: A[S^{-1}] \to B$ such that the diagram below commutes:
%$$\xymatrix{
%A \ar[r]^{\eta} \ar[rd]^{\varphi} & A[S^{-1}] \ar@{-->}[d]^{\exists ! \psi}\\
 % & B
%  }$$

%Now we claim that $\psi$ is a surjection. Given any $z \in B$, since $\varphi: A \to B$ satisfies $(i)$, there are elements $c \in S$ and $d \in A$ such that $z \cdot \varphi(c) = \varphi(d)$, that is to say $z = \varphi(d)\cdot (\varphi(c))^{-1}$. Taking $w = \eta(d)\cdot (\eta(c))^{-1} \in A[S^{-1}]$ we get $\psi(w) = \psi(\eta(d)\cdot (\eta(c))^{-1}) = \psi(\eta(d))\cdot \psi(\eta(c))^{-1} = \varphi(d)\cdot (\varphi(c))^{-1} = z$.\\

%We claim, also, that $\psi$ is an injective ring homomorphism. Let $w \in A[S^{-1}]$ be an element such that $\psi(w)=0$. Since $w \in A[S^{-1}]$ there are $c \in S$ and $d \in A$ such that $w = \eta(d)\cdot (\eta(c))^{-1}$. The condition $\psi(w)=0$ means that $\psi(\eta(d)\cdot (\eta(c))^{-1})= \psi(\eta(d))\cdot \psi((\eta(c)))^{-1} = 0$, so $\psi(\eta(d))=0$, and since $\psi \circ \eta = \varphi$, it follows that $\varphi(d)=0$. By the property $(ii)$, there exists $c' \in S$ such that $c' \cdot d = 0$, so $\eta(c')\cdot \eta(d) = \eta(c' \cdot d) =  0$, and since $\eta(c') \in A[S^{-1}]^{\times}$, it follows that $\eta(d)=0$. Thus $w = \eta(d)\cdot (\eta(c'))^{-1} = 0\cdot (\eta(c'))^{-1} = 0$.\\

%It follows that $\psi$ is an isomorphism between $B$ and $A[S^{-1}]$.
\end{proof}

The two preceding results gives us the following:

\begin{theorem}\label{37}Let $A$ be a commutative ring with unity and $S \subseteq A$. Then $\varphi: A \to B$ is isomorphic to the localization $\eta: A \to A[S^{-1}]$ if and only if:
\begin{itemize}
\item[(i)]{$(\forall b \in B)(\exists c \in S)(\exists d \in A)(b \cdot \varphi(c) = \varphi(d))$}
\item[(ii)]{$(\forall b \in A)(\varphi(b)=0 \to (\exists c \in S)(c \cdot b = 0))$}
\end{itemize}
hold.
\end{theorem}

Thus we have obtained a characterization of the localization of a commutative ring. For $\mathcal{C}^{\infty}-$rings we have the analogous result, that generalizes \textbf{Theorem 1.4} of \cite{rings1}, in the sense that it is an equivalence (an ``if and only if'' statement) and that  $S$ needs not to be a singleton:

\begin{theorem}\label{38}Let $A$ be a $\mathcal{C}^{\infty}-$ring $\Sigma \subset A$ a set. Then $\varphi: A \to B$ is isomorphic to the smooth localization $\eta: A \to A\{\Sigma^{-1}\}$ if and only if:
\begin{itemize}
\item[(i)]{$(\forall b \in B)(\exists c \in \Sigma^{\infty-{\rm sat}})(\exists d \in A)(b \cdot \varphi(c) = \varphi(d))$}
\item[(ii)]{$(\forall b \in A)(\varphi(b)=0 \to (\exists c \in \Sigma^{\infty-{\rm sat}})(c \cdot b = 0))$}
\end{itemize}
hold.
\end{theorem}

We postpone the proof of this theorem, giving it right after \textbf{Remark \ref{raca}}.

\begin{theorem}\label{340}Let $A$, $\widetilde{A}$ be $\mathcal{C}^{\infty}-$rings and let $\eta : A \to \widetilde{A}$ be a $\mathcal{C}^{\infty}-$rings homomorphism such that:
\begin{itemize}
  \item[(i)]{$(\forall d \in \widetilde{A})(\exists b \in A)(\exists c \in A)(\eta(c) \in \widetilde{A}^{\times} \& (d \cdot \eta(c) = \eta(b)))$;}
  \item[(ii)]{$(\forall b \in A)((\eta(b) = 0_{\widetilde{A}}) \to (\exists c \in A)((\eta(c) \in \widetilde{A}^{\times}) \wedge (b \cdot c = 0_A)))$}
\end{itemize}
Then $\eta : A \to \widetilde{A}$ is isomorphic to ${\rm Can}_{S_{\eta}} : A \to A\{ {S_{\eta}}^{-1}\}$, where $S_{\eta} = \eta^{\dashv}[\widetilde{A}^{\times}]$.
\end{theorem}
\begin{proof}
First we show that $\eta: A \to \widetilde{A}$ has the universal property which characterize ${\rm Can}_{S_{\eta}}$.

Let $f: A \to B$ be a $\mathcal{C}^{\infty}-$rings homomorphism such that $f[S_{\eta}] \subseteq B^{\times}$. We are going to show there is a unique $\mathcal{C}^{\infty}-$rings homomorphism $\widetilde{f} : \widetilde{A} \to B$ such that the following diagram commutes:
$$\xymatrix{
  A \ar[r]^{\eta} \ar[rd]_{f} & \widetilde{A} \ar@{.>}[d]^{\widetilde{f}}\\
      & B
}$$
i.e., such that $\widetilde{f}\circ \eta = f$.

First we notice that:
$$\eta[S_{\eta}] = \eta[\eta^{\dashv}[\widetilde{A}^{\times}]] \subseteq \widetilde{A}^{\times}.$$

\textbf{Candidate and Uniqueness of $\widetilde{f}$:} Let $\widetilde{f}_1, \widetilde{f}_2: \widetilde{A} \to B$ be two $\mathcal{C}^{\infty}-$rings homomorphisms such that the following diagram commutes:

$$\xymatrix{
  A \ar[r]^{\eta} \ar@/_/[dr]_f & \widetilde{A} \ar[d]^{\widetilde{f}_1} \ar@<-1ex>[d]_{\widetilde{f}_2} & \ar[l]_{\eta} A \ar@/^/[dl]^{f}\\
      & B &
}$$
that is to say, such that $\widetilde{f}_1 \circ \eta = f = \widetilde{f}_2 \circ \eta$.

Given any $d \in \widetilde{A}$, by the hypothesis (i) there exist $b,c \in A$, $\eta(c) \in \widetilde{A}^{\times}$ such that $d = \dfrac{\eta(b)}{\eta(c)}$, so:
\begin{multline*}\widetilde{f}_1(d) = \widetilde{f}_1(\eta(b)\cdot \eta(c)^{-1}) = \widetilde{f}_1(\eta(b))\cdot \widetilde{f}_1(\eta(c))^{-1}= \\
 = (\widetilde{f}_1 \circ \eta)(b) \cdot ((\widetilde{f}_1 \circ \eta)(c))^{-1} = f(b)\cdot f(c)^{-1} = (\widetilde{f}_2 \circ \eta)(b) \cdot ((\widetilde{f_2} \circ \eta)(c))^{-1} =\\ =\widetilde{f}_2(\eta(b))\cdot \widetilde{f}_2(\eta(c)^{-1}) = \widetilde{f}_2(\eta(b)\cdot\eta(c)^{-1}) = \widetilde{f}_2(d)
\end{multline*}
so we conclude that $(\forall d \in \widetilde{A})(\widetilde{f}_1(d) = \widetilde{f}_2(d))$ and $\widetilde{f}(\frac{\eta(b)}{\eta(c)}) = f(b)\cdot f(c)^{-1}$. Thus $\widetilde{f}_1 = \widetilde{f}_2$.\\

\textbf{Existence of $\widetilde{f}$:} We know that for every $d \in \widetilde{A}$ there exist $b,c \in A$, $\eta(c) \in \widetilde{A}^{\times}$, such that $ d = \eta(b)\cdot \eta(c)^{-1}$. Define the following relation: $\widetilde{f} = \{ (d, f(b)\cdot f(c)^{-1}) | d \in \widetilde{A} \} \subseteq \widetilde{A} \times B$. We claim that $\widetilde{f} $ is a function.\\

\textbf{Claim:} $\widetilde{f}$ is a univocal relation.\\

Given $d \in \widetilde{A}$, let $b,c,b',c' \in A$, $\eta(c), \eta(c') \in \widetilde{A}^{\times}$ be such that $\eta(b)\cdot \eta(c)^{-1} = d = \eta(b')\cdot \eta(c')^{-1}$, so  $\eta(b)\cdot \eta(c') = \eta(b')\cdot \eta(c)$. Then we have $\eta(b \cdot c') = \eta(b) \cdot \eta(c') = \eta(b')\cdot \eta(c)$, so $\eta(b \cdot c' - b' \cdot c)= 0$. By (ii), since $\eta(b \cdot c' - b' \cdot c)= 0$, there is some $a \in A$ such that $\eta(c) \in \widetilde{A}^{\times}$, i.e., $a \in \eta^{\dashv}[\widetilde{A}^{\times}] = S_{\eta}$, and $a \cdot (b\cdot c' - b' \cdot c) = 0$, i.e., $b\cdot c' \cdot a = b' \cdot c \cdot a$. We now have:

$$ f(b) \cdot f(c') \cdot f(a) = f(b \cdot c' \cdot a) = f(b' \cdot c \cdot a) = f(b') \cdot f(c) \cdot f(a).$$

Now since $f[S_{\eta}] \subseteq B^{\times}$ and $a \in S_{\eta}$, from the preceding equations we obtain, by cancelling $f(a)$:
$$f(b)\cdot f(c') = f(b')\cdot f(c)$$
and since $\eta(c), \eta(c') \in \widetilde{A}^{\times}$, $c, c' \in \eta^{\dashv}[\widetilde{A}^{\times}] = S_{\eta}$ so $f(c), f(c') \in B^{\times}$, so:
$$f(b)\cdot f(c)^{-1} = f(b')\cdot f(c')^{-1}.$$

We conclude that $(\forall d \in \widetilde{A})(((d, f(b)\cdot f(c)^{-1}) \in \widetilde{f} ) \wedge ((d, f(b')\cdot f(c')^{-1}) \in \widetilde{f} ) \to f(b) \cdot f(c)^{-1} = f(b') \cdot f(c')^{-1})$.\\

The fact that $\widetilde{f} $ is a total relation follows immediately from item (i).Thus, we denote $ (d, f(b)\cdot f(c)^{-1}) \in \widetilde{f} $ simply by $\widetilde{f} (d) = f(b)\cdot f(c)^{-1}$, as usually we do for functions.\\

Therefore, there exists exactly one function $\widetilde{f}: \widetilde{A} \to B$ such that $\widetilde{f} \circ \eta = f$. \\

Now we show that $\widetilde{A} \cong A\{ S_{\eta}^{-1}\}$. We begin by noticing that ${\rm Can}_{S_{\eta}}[S_{\eta}] \subseteq A\{ S_{\eta}^{-1}\}^{\times}$ - by the very definition of ${\rm Can}_{S_{\eta}}$. For what we have seen above, there exists a unique function $\widetilde{{\rm Can}_{S_{\eta}}}: \widetilde{A} \to A\{ S_{\eta}^{-1} \}$ such that $\widetilde{{\rm Can}_{S_{\eta}}}\circ \eta = {\rm Can}_{S_{\eta}}$. Now, from the universal property of ${\rm Can}_{S_{\eta}}$ there exists a unique $\mathcal{C}^{\infty}-$ring homomorphism, $\widehat{\eta}: A\{ S_{\eta}^{-1}\} \to \widetilde{A}$, such that $\widehat{\eta}\circ {\rm Can}_{S_{\eta}} = \eta$.\\

\textbf{Claim:} $\widehat{\eta}$ is a bijection whose inverse is $\widetilde{{\rm Can}_{S_{\eta}}}$, and that will prove that $\widetilde{{\rm Can}_{S_{\eta}}}$ is a $\mathcal{C}^{\infty}-$rings isomorphism.\\

Now, $(\widehat{\eta}\circ \widetilde{{\rm Can}_{S_{\eta}}}) \circ \eta = \widehat{\eta} \circ (\widetilde{{\rm Can}_{S_{\eta}}} \circ \eta) = \widehat{\eta} \circ {\rm Can}_{S_{\eta}} = \eta = {\rm id}_{\widetilde{A}} \circ \eta$, so:
$$(\widehat{\eta}\circ \widetilde{{\rm Can}_{S_{\eta}}}) \circ \eta = {\rm id}_{\widetilde{A}} \circ \eta.$$

We have seen, however, that there is exactly one function $\widetilde{\varphi}$ such that $\widetilde{\varphi} \circ \eta = \eta$, so it follows that ${\rm id}_{\widetilde{A}} = \widehat{\eta} \circ \widetilde{{\rm Can}_{S_{\eta}}}$.\\

On the other hand,

$$\widetilde{{\rm Can}_{S_{\eta}}} \circ \widehat{\eta} \circ {\rm Can}_{S_{\eta}} = \widetilde{{\rm Can}_{S_{\eta}}} \circ \eta = {\rm Can}_{S_{\eta}} = {\rm id}_{A\{ S_{\eta}^{-1} \}} \circ {\rm Can}_{S_{\eta}}.$$

Once again, by the universal property of ${\rm Can}_{S_{\eta}}$ we have:

$${\rm id}_{A\{S_{\eta}^{-1}\}} = \widetilde{{\rm Can}_{S_{\eta}}} \circ \widehat{\eta}.$$

Hence $\widetilde{{\rm Can}_{S_{\eta}}}$ is the $\mathcal{C}^{\infty}-$rings isomorphism between $\eta$ and ${\rm Can}_{S_{\eta}}$, that is, it is a $\mathcal{C}^{\infty}-$rings isomorphism such that the following diagram commutes:
  $$\xymatrix{A \ar[r]^{\eta} \ar[dr]^{{\rm Can}_{S_{\eta}}} & \widetilde{A} \ar@{>->>}[d]^{\widetilde{{\rm Can}_{S_{\eta}}}}\\
        &  A\{ S_{\eta}^{-1}\}
  }.$$
\end{proof}

In order to smoothly localize larger subsets of $\mathcal{C}^{\infty}(\R^n)$ for some $n \in \mathbb{N}$, say $\Sigma$, which is a set that contains possibly a non-countable amount of elements,  we can proceed as follows. First notice that we can obtain the $\mathcal{C}^{\infty}-$ring of fractions of $\mathcal{C}^{\infty}(\mathbb{R}^n)$ with respect to the singleton $\Sigma = \{f: \R^n \to \R\}$, provided that $f \not\equiv 0$. Whenever $\Sigma = \{ f_1, \cdots, f_k \}$ for some $k \in \N$, inverting $\Sigma$ is equivalent to inverting $\prod \Sigma = f_1 \cdot f_2 \cdots f_{k-1}\cdot f_k$. In the case that $\Sigma$ is infinite, first we decompose it as the union of its finite subsets:

$$\Sigma = \bigcup_{\Sigma' \subseteq_{\rm fin} \Sigma} \Sigma'$$

Note that $\mathcal{S} =\{ \Sigma' \subseteq \Sigma | \Sigma' \,\,\mbox{is}\,\, \mbox{finite}\}$ is partially ordered by the  inclusion relation. Also,  whenever $\Sigma' \subseteq \Sigma''$, since $\eta_{\Sigma''}[\Sigma']\subseteq \eta_{\Sigma''}[\Sigma'']\subseteq (A\{ {\Sigma''}^{-1}\})^{\times}$, by the universal property of $\eta_{\Sigma''}: A \to A\{ {\Sigma''}^{-1}\}$, there is a unique $\mathcal{C}^{\infty}-$homomorphism $\alpha_{\Sigma' \Sigma''}: A\{ {\Sigma'}^{-1}\} \rightarrow A\{ {\Sigma''}^{-1}\}$ such that the following diagram commutes:

$$\xymatrixcolsep{5pc}\xymatrix{
A \ar[r]^{\eta_{\Sigma'}} \ar[dr]_{\eta_{\Sigma''}} & A\{ {\Sigma'}^{-1}\} \ar@{-->}[d]^{\exists ! \alpha_{\Sigma' \Sigma''}}\\
 & A\{ {\Sigma''}^{-1}\}}$$

It is simple to prove, using the ``uniqueness part'' of the $\mathcal{C}^{\infty}-$homomorphism obtained via universal property, that for any finite $\Sigma'$ we have $\alpha_{\Sigma' \Sigma'} = {\rm id}_{A\{ {\Sigma'}^{-1}\}}$, and given any finite $\Sigma', \Sigma''$ and $\Sigma'''$ such that $\Sigma' \subseteq \Sigma'' \subseteq \Sigma'''$,  $\alpha_{\Sigma'' \Sigma'''} \circ \alpha_{\Sigma' \Sigma''} = \alpha_{\Sigma'\Sigma'''}$, so we have an inductive system:
$$\{ \alpha_{\Sigma' \Sigma''}: A\{ {\Sigma'}^{-1}\} \to A\{ {\Sigma''}^{-1}\} | (\Sigma', \Sigma'' \in \mathcal{S}) \& (\Sigma' \subseteq \Sigma'')\}$$

We take, thus:

$$A\{ \Sigma^{-1}\} = \varinjlim_{\Sigma' \subseteq_{\rm fin} \Sigma} A\{ {\Sigma'}^{-'}\} = \varinjlim_{\Sigma' \subseteq_{\rm fin} \Sigma} A \left\{ {\prod \Sigma'}^{-1} \right\}$$

$$\xymatrixcolsep{5pc}\xymatrix{
 & A\{ \Sigma^{-1}\} & \\
A\{ {\Sigma'}^{-'}\} \ar[ur]^{\alpha_{\Sigma'}} \ar[rr]_{\alpha_{\Sigma' \Sigma''}} & & A\{ {\Sigma''}^{-1}\} \ar[ul]_{\alpha_{\Sigma''}}}$$

Hence, given any $\mathcal{C}^{\infty}-$ring $A$ and any $S \subseteq A$, we can construct:

$$\eta_{\Sigma} = \alpha_{\Sigma'}\circ \eta_{\Sigma'}: A \rightarrow A\{ {\Sigma}^{-1}\} $$

It is easy to prove that $\eta_{\Sigma}$ has the universal property which characterizes the $\mathcal{C}^{\infty}-$ring of fractions of $A$ with respect to $\Sigma$.\\

\begin{remark}In the case that $A = \mathcal{C}^{\infty}(\mathbb{R}^n)$ and  $\Sigma = \{ f: \mathbb{R}^n \to \mathbb{R}\}$, we have $\Sigma^{\infty-{\rm sat}} = \{ g \in \mathcal{C}^{\infty}\,(\mathbb{R}^n) | U_f \subseteq U_g\} = \{ g \in \mathcal{C}^{\infty}\,(\mathbb{R}^n) | Z(g) \subseteq Z(f) \}$.
\end{remark}

The following theorem gives us concretely the $\mathcal{C}^{\infty}-$ring of fractions of a finitely generated free $\mathcal{C}^{\infty}-$ring with respect to one element.\\

\begin{theorem}Let $\varphi \in \mathcal{C}^{\infty}(\mathbb{R}^n)$ and $U_{\varphi} = \{ x \in \mathbb{R}^n | \varphi(x)\neq 0\} = {\rm Coz}(\varphi)$. Then:
$$\mathcal{C}^{\infty}(U_{\varphi}) \cong \dfrac{\mathcal{C}^{\infty}(\mathbb{R}^{n+1})}{\langle \{ y \cdot \varphi(x) - 1 \}\rangle}$$
\end{theorem}

Let us now consider the following result, credited to Ortega and Mu\~{n}oz by I. Moerdijk and G. Reyes in \cite{rings1}:

\begin{theorem}Let $U \subseteq \mathbb{R}^n$ be open, and $g \in \mathcal{C}^{\infty}(U)$. Then there are $h,k \in \mathcal{C}^{\infty}(\mathbb{R}^n)$ with $U_k = U$ and $g \cdot k \displaystyle\upharpoonright_{U} \equiv h \displaystyle\upharpoonright_{U}$. where $U_k = \mathbb{R}^n \setminus Z(k)$ and $Z(k) = \{ x \in \mathbb{R}^n | k(x)=0\}$.
\end{theorem}

%\begin{theorem}\label{creio}(\textbf{Lemma 2} of \cite{JR} )Let $\varphi \in \mathcal{C}^{\infty}(\mathbb{R}^n)$ and $U_{\varphi} = \{ x \in \mathbb{R}^n | f(x) \neq 0\}$. We have:
%$$\mathcal{C}^{\infty}(U_{\varphi}) \cong \dfrac{\mathcal{C}^{\infty}(\mathbb{R}^{n+1})}{\langle y \cdot \varphi(x) - 1\rangle}$$
%and

%$$\begin{array}{cccc}
%    \eta_{\varphi}: & \mathcal{C}^{\infty}(\mathbb{R}^n) & \rightarrow & \mathcal{C}^{\infty}(U_{\varphi}) \\
%     & g & \mapsto & g\upharpoonright_{U_{\varphi}}
%  \end{array}$$
%
%\end{theorem}

\begin{theorem}\label{cara} Let $A$ be a $\mathcal{C}^{\infty}-$ring and $S \subseteq A$. An element $\lambda = \dfrac{\eta_S(c)}{\eta_S(b)}$ (with $c \in A$ and $b \in S^{\infty-{\rm sat}}$) is invertible in $A\{ S^{-1}\}$ if, and only if, there are elements $d \in S^{\infty-{\rm sat}}$ and $c' \in A$ such that $dc'c \in S^{\infty-{\rm sat}}$, that is,
$$\dfrac{\eta_S(c)}{\eta_S(b)} \in (A\{ S^{-1}\})^{\times} \iff (\exists d \in S^{\infty-{\rm sat}})(\exists c' \in A)(d\cdot c' \cdot c \in S^{\infty-{\rm sat}}).$$
\end{theorem}
\begin{proof}
Suppose $\dfrac{\eta_S(c)}{\eta_S(b)} \in (A\{ S^{-1}\})^{\times}$, so there are $c' \in A$ and $b' \in S^{\infty-{\rm sat}}$ such that:
$$\dfrac{\eta_S(c)}{\eta_S(b)} \cdot \dfrac{\eta_S(c')}{\eta_S(b')} = 1_{A\{ S^{-1}\}} = \eta_S(1_A).$$
$$\eta_S(c \cdot c') = \eta_S(b \cdot b')$$
$$\eta_S(c \cdot c' - b \cdot b') = 0$$
By \textbf{Theorem \ref{38}}, there is some $d \in S^{\infty-{\rm sat}}$ such that:
$$d \cdot (c \cdot c' - b \cdot b') = 0$$
$$d \cdot c \cdot c' = d \cdot b \cdot b' \in S^{\infty-{\rm sat}}$$
where $d \cdot b \cdot b' \in S^{\infty-{\rm sat}}$ because it is a product of elements of $S^{\infty-{\rm sat}}$, which is a submonoid of $A$.\\

Conversely, suppose that $\dfrac{\eta_S(c)}{\eta_S(b)} \in A\{ S^{-1}\}$ with $b \in S^{\infty-{\rm sat}}$ is an element for which there are elements $d \in S^{\infty-{\rm sat}}$ and $c' \in A$ such that $d \cdot c \cdot c' \in S^{\infty-{\rm sat}}$. We have $\eta_S(d \cdot c' \cdot c) \in (A\{ S^{-1}\})^{\times}$ and $b \in S^{\infty-{\rm sat}}$, so $\eta_S(b) \in (A\{ S^{-1}\})^{\times}$, hence
$$\dfrac{\eta_S(d \cdot c' \cdot c)}{\eta_S(b)} \in (A\{ S^{-1}\})^{\times}.$$

Since
$$\dfrac{\eta_S(c)}{\eta_S(b)} \cdot \eta_S(d \cdot c') = \dfrac{\eta_S(d \cdot c' \cdot c)}{\eta_S(b)} \in (A\{ S^{-1}\})^{\times}$$
it follows that $\dfrac{\eta_S(c)}{\eta_S(b)} \in (A\{ S^{-1}\})^{\times}$, for if $\alpha \cdot \beta$ is invertible, then both $\alpha$ and $\beta$ are invertible.
\end{proof}

Now we have the following:

\begin{proposition}\label{pani}Let $U \subseteq \mathbb{R}^n$ be any open subset and define $S_U = \{ g \in \mathcal{C}^{\infty}(\mathbb{R}^n) | U \subseteq U_g\} \subseteq \mathcal{C}^{\infty}(\mathbb{R}^n)$. The $\mathcal{C}^{\infty}-$ring of fractions of $\mathcal{C}^{\infty}(\mathbb{R}^n)$ with respect to the set $S_U$:
$$\eta_{S_U}: \mathcal{C}^{\infty}(\mathbb{R}^n) \to \mathcal{C}^{\infty}(\mathbb{R}^n)\{ {S_U}^{-1}\}$$
is isomorphic to the restriction map:
$$\begin{array}{cccc}
    \rho : & \mathcal{C}^{\infty}(\mathbb{R}^n) & \rightarrow & \mathcal{C}^{\infty}(U) \\
     & h & \mapsto & h \displaystyle\upharpoonright_{U}
  \end{array}$$
\end{proposition}
\begin{proof}
See p. 39 of \cite{BMTSCA}.
%We are going to show that $\rho \cong \eta_{S_U}$ using \textbf{Theorem \ref{340}} with $A = \mathcal{C}^{\infty}(\mathbb{R}^n)$, $\widetilde{A} = \mathcal{C}^{\infty}\,(U)$ and $\eta = \rho$.\\

%Note that $\rho[S_U] \subseteq \mathcal{C}^{\infty}(U)^{\times}$.\\

%Let's verify the first item, (i):\\

%By the \textbf{Theorem 1.3} of \cite{rings1}, given any $g \in \mathcal{C}^{\infty}(U)$, there are $h,k \in \mathcal{C}^{\infty}(\mathbb{R}^n)$ with $U_k = U$, such that $g \cdot k\upharpoonright_U = h\upharpoonright_U$. Since $U=U_k$, $\rho(k) \in \mathcal{C}^{\infty}(U)^{\times}$, and we have, thus:
%$$(\rho(k) \in \mathcal{C}^{\infty}(U)^{\times})\&(g \cdot \rho(k) = \rho(h)).$$

%For item (ii), suppose that $g \in \mathcal{C}^{\infty}(\mathbb{R}^n)$ is such that $g\upharpoonright_U = 0$, so $U \subseteq Z(g)$. We know, by  a well-known theorem proved by Whitney (see \textbf{Theorem 5.0.5} of \cite{tese}), that given the open subset $U$ of $\mathbb{R}^n$, there is some $h \in \mathcal{C}^{\infty}(\mathbb{R}^n)$ such that $U=U_h$, so $Z(h \cdot g) = Z(h) \cup Z(g) = U \cup (\mathbb{R}^n \setminus U) = \mathbb{R}^n$. We have both:

%$$\rho(g) = 0 \in \mathcal{C}^{\infty}(U)$$

%and

%$$g \cdot h = 0 \in \mathcal{C}^{\infty}(\mathbb{R}^n),$$

%so item (ii) is also fulfilled. By \textbf{Theorem \ref{340}}, the result follows.
\end{proof}

\begin{remark}Let $A$ be a $\mathcal{C}^{\infty}-$ring and $a \in A$. In general, the $\mathcal{C}^{\infty}-$ring of fractions of $A$ with respect to $a$ \textbf{is not a local $\mathcal{C}^{\infty}-$ring}. Let us consider the case on which $A = \mathcal{C}^{\infty}(\mathbb{R}^n)$ and $a = f: \mathbb{R}^n \to \mathbb{R}$ is such that $\neg (f \equiv 0)$. By the \textbf{Theorem 1.3} of \cite{rings1}, $A\{ a^{-1}\} \cong \mathcal{C}^{\infty}(\mathbb{R}^n)\{ f^{-1}\}$, and $\mathcal{C}^{\infty}(\mathbb{R}^n)\{ f^{-1}\} \cong \mathcal{C}^{\infty}(U_f)$, where $U_f = {\rm Coz}\,(f) = \mathbb{R}^n \setminus Z(f)$. For every $x \in U_f$, we have a maximal ideal:
$$\mathfrak{m}_x = \{ g \in \mathcal{C}^{\infty}(U_f) | g(x) = 0\},$$
hence a \textit{continuum} of maximal ideals.
\end{remark}

\begin{remark}\label{raca}In the context of \textbf{Proposition \ref{pani}}, note that if $f \in \mathcal{C}^{\infty}(\mathbb{R}^n)$ is such that $U=  U_f$, we have $S_{U_f} = \{ f\}^{\infty-{\rm sat}}$.
\end{remark}

Now we are ready to prove \textbf{Theorem \ref{38}}.\\

\textbf{Proof of Theorem \ref{38}:} Let $A$ be a $\mathcal{C}^{\infty}-$ring. We are going to prove the result first for free $\mathcal{C}^{\infty}-$rings, then for quotients and finally for colimits.\\

The case where $A = \mathcal{C}^{\infty}(\mathbb{R}^n)$ was proved in \textbf{Proposition \ref{pani}};\\

Suppose, now, that $A = \dfrac{\mathcal{C}^{\infty}(\mathbb{R}^n)}{I}$ for some ideal $I$. By \textbf{Corollary \ref{Jeq}}, it follows that:

$$\dfrac{\mathcal{C}^{\infty}(\mathbb{R}^n)}{I}\lbrace {f + I}^{-1} \rbrace \cong \dfrac{\mathcal{C}^{\infty}(\mathbb{R}^n)\{ f^{-1}\}}{\langle \eta_f[I]\rangle}$$

so items (i) and (ii) of \textbf{Theorem \ref{340}} hold for the quotient.\\

Up to this point, we have proved the result for a finitely generated $\mathcal{C}^{\infty}-$ring,  $A$, and $\Sigma = \{ f\}$. Given any finite $S$, say $S = \{ f_1, \cdots, f_{\ell}\}$, since $A\{ S^{-1}\} = A\{ {\prod_{i = 1}^{\ell}f_i}^{-1}\}$, the result follows too.\\

Let $A$ be a finitely generated $\mathcal{C}^{\infty}-$ring and  $S \subseteq A$ be any set. Write:

$$S = \bigcup_{S' \subseteq_{\rm fin} S}S'$$

Since $A\{ S^{-1}\} \cong \varinjlim_{S' \subseteq S} A\{ {S'}^{-1}\}$, the items (i) and (ii) hold for $A\{ S^{-1}\}$, for a finitely generated $\mathcal{C}^{\infty}-$ring $A$ and any set $S$.\\

Finally, given any $\mathcal{C}^{\infty}-$ring $B$ and any set $S \subseteq B$, write $B$ as a directed colimit of its finitely generated $\mathcal{C}^{\infty}-$subrings

$$\xymatrixcolsep{5pc}\xymatrix{
 & B & \\
B_i \ar@/^/[ur]^{\jmath_i} \ar[rr]_{\beta_{ij}} & & B_j \ar@/_/[ul]_{\jmath_j}}$$

so $B \cong \varinjlim_{B_i \subseteq_{\rm f.g.} B}B_i$ and define $S_i = \jmath_i^{\dashv}[S]$. Since items (i) and (ii) hold for every $B_i\{ {S_i}^{-1}\}$, the same is true for

$$B\{ S^{-1}\} \cong \varinjlim_{B_i \subseteq_{\rm fin} B} B_i\{ {S_i}^{-1}\},$$

and the result follows. \hfill $\square$ \\

Now we analyze the concept of a ``$\mathcal{C}^{\infty}-$radical ideal'' in the theory of $\mathcal{C}^{\infty}-$rings, which plays a similar role to the one played by radical ideals in Commutative Algebra. This concept is presented by I. Moerdijk and G. Reyes in \cite{rings1} in 1986 and explored in more details in \cite{rings2}. \\

Contrary to the concepts of $\mathcal{C}^{\infty}-$fields, $\mathcal{C}^{\infty}-$domains and local $\mathcal{C}^{\infty}-$rings, the concept of a $\mathcal{C}^{\infty}-$radical of an ideal can not be brought from Commutative Algebra via the forgetful functor. Recall that the radical of an ideal $I$ of a commutative unital ring $R$ is given by:

$$\sqrt{I} = \{ x \in R | (\exists n \in \mathbb{N})(x^n \in I)\}$$

and this concept can be characterized by:

$$\sqrt{I} = \bigcap \{ \mathfrak{p} \in {\rm Spec}\,(R) | I \subseteq \mathfrak{p}  \} = \{ x \in R | \left( \dfrac{R}{I}\right)[(x+I)^{-1}] \cong 0\}.$$

We use the latter equality in order to motivate our definition.\\

\begin{definition}\label{defrad}Let $A$ be a $\mathcal{C}^{\infty}-$ring and let $I \subseteq A$ be a proper ideal. The \index{$\mathcal{C}^{\infty}-$radical}\textbf{$\mathcal{C}^{\infty}-$radical of $I$} is given by:

$$\sqrt[\infty]{I}:= \{ a \in A | \left( \dfrac{A}{I}\right)\{ (a+I)^{-1}\} \cong 0\}$$
\end{definition}

Contrary to what happens in ordinary Commutative Algebra, it is not evident, up to this point, that whenever $I$ is an ideal of a $\mathcal{C}^{\infty}-$ring, $\sqrt[\infty]{I}$ is also an ideal. Later on we prove this fact using ``steps'': first we prove the result for free $\mathcal{C}^{\infty}-$rings, then for colimits and finally for quotients.\\ 

\begin{proposition}\label{alba}Let $A$ be a $\mathcal{C}^{\infty}-$ring and let $I \subseteq A$ be any ideal. We have the following equalities:
  $$\sqrt[\infty]{I} = \{ a \in A | (\exists b \in I)\wedge(\eta_a(b) \in (A\{ a^{-1}\})^{\times}) \} = \{ a \in A | I \cap \{ a\}^{\infty-{\rm sat}} \neq \varnothing\}$$
  where $\eta_a : A \to A\{ a^{-1}\}$ is the morphism of fractions with respect to $\{ a\}$.\\
\end{proposition}
\begin{proof}See \textbf{Proposition 19}, p. 63 of \cite{BMTSCA}.
\end{proof}

\begin{theorem}\label{Hilfssatz} Let $A$ and $B$ be two $\mathcal{C}^{\infty}-$rings, $S \subseteq A$, $p \twoheadrightarrow B$ a surjective map and $\eta_S : A \to A\{ S^{-1}\}$ and $\eta_{p[S]} : B \to B\{ p[S]^{-1}\}$ be the canonical $\mathcal{C}^{\infty}-$ring of fractions. There is a unique surjective map $q: A\{ S^{-1}\} \twoheadrightarrow B\{ p[S]^{-1}\}$ such that the following square commutes:
$$\xymatrixcolsep{5pc}\xymatrix{
A \ar[r]^{\eta_S} \ar@{->>}[d]^{p} & A\{ S^{-1}\} \ar@{->>}[d]^{q}\\
B \ar[r]^{\eta_{p[S]}} & B\{ p[S]^{-1}\}}$$
\end{theorem}
\begin{proof}
See p. 65 of \cite{BMTSCA}.
\end{proof}

\subsection{A Category of Pairs}

Let us now consider the category whose objects are pairs $(A,S)$, where $A$ is any $\mathcal{C}^{\infty}-$ring and $S \subseteq A$ is any of its  subsets, and whose morphisms between $(A,S)$ and $(B,T)$ are precisely the $\mathcal{C}^{\infty}-$morphisms $f: A \to B$ such that $f[S] \subseteq T$. In order to distinguish the latter morphism of the category of pairs from the morphisms of $\mathcal{C}^{\infty}{\rm \bf Rng}$, we shall denote it by $f_{ST}: (A,S) \to (B,T)$. We denote this category by $\mathcal{C}^{\infty}_2$. We can define, naturally, a forgetful functor $U: \mathcal{C}^{\infty}_2 \to \mathcal{C}^{\infty}{\rm \bf Rng}$, which maps any pair $(A,S)$ to $A$ and any $f_{ST}: (A,S) \to (B,T)$ to $f_{ST}: A \to B$.\\

There are some natural ways to define full and faithful functors fom $\mathcal{C}^{\infty}{\rm \bf Rng}$ to $\mathcal{C}^{\infty}_{2}$, such as:

$$\begin{array}{cccc}
 \imath_0 : & \mathcal{C}^{\infty}{\rm \bf Rng} & \rightarrow & \mathcal{C}^{\infty}_2  \\
     & (A \stackrel{f}{\rightarrow} B) & \mapsto & ((A,\varnothing) \stackrel{f_{\varnothing \varnothing}}{\rightarrow} (B,\varnothing))
\end{array}$$

and

$$\begin{array}{cccc}
 \imath_1 : & \mathcal{C}^{\infty}{\rm \bf Rng} & \rightarrow & \mathcal{C}^{\infty}_2  \\
     & (A \stackrel{f}{\rightarrow} B) & \mapsto & ((A,A) \stackrel{f_{A B}}{\rightarrow} (B,B))
\end{array}$$

It is not difficult to see that the forgetful functor, $U: \mathcal{C}^{\infty}_2 \to \mathcal{C}^{\infty}{\rm \bf Rng}$, has $\imath_0$ as a left adjoint and $\imath_1$ as a right adjoint, so $\imath_0 \dashv U \dashv \imath_1$.

Another way to define a full and faithful functor is given by:

$$\begin{array}{cccc}
 \jmath : & \mathcal{C}^{\infty}{\rm \bf Rng} & \rightarrow & \mathcal{C}^{\infty}_2  \\
     & (A \stackrel{f}{\rightarrow} B) & \mapsto & ((A,A^{\times}) \stackrel{f_{A^{\times} B^{\times}}}{\rightarrow} (B,B^{\times}))
\end{array}$$

In this section we are going to define a functor $F: \mathcal{C}^{\infty}_{2} \to \mathcal{C}^{\infty}{\rm \bf Rng}$ which has, as a right adjoint, the functor $\jmath$, that is, given any $(A,S)\in {\rm Obj}\,(\mathcal{C}^{\infty}_2)$ and any $B \in {\rm Obj}\,(\mathcal{C}^{\infty}{\rm \bf Rng})$, there is a natural bijection, ${\rm Hom}_{\mathcal{C}^{\infty}{\rm \bf Rng}}(F(A,S),B) \cong {\rm Hom}_{\mathcal{C}^{\infty}_2}((A,S),\jmath(B))$. As a consequence of $F: \mathcal{C}^{\infty}_2 \to \mathcal{C}^{\infty}{\rm \bf Rng}$ being a left adjoint, it must preserve all colimits - a fact that will become clear along this section.\\

The proofs of the results of this section can be found in \textbf{Section 3.2} of \cite{BMTSCA}.\\

Given $(A,S) \in {\rm Obj}\,(\mathcal{C}^{\infty}_{2})$, let:
$$F(A,S) = A\{ S^{-1}\}.$$

Also, given $(A,S), (B,T) \in {\rm Obj}\,(\mathcal{C}^{\infty}_{2})$ and $(A,S) \stackrel{f_{ST}}{\rightarrow} (B,T) \in {\rm Mor}\,(\mathcal{C}^{\infty}_{2})$, since $f_{ST}[S] \subseteq T$ and $\eta_T[T] \subseteq {B\{ T^{-1}\}}^{\times}$, we have $(\eta_T \circ f_{ST})[S] \subseteq {B\{ T^{-1}\}}^{\times}$, so by the universal property of $\eta_S: A \to A\{ S^{-1}\}$ there is a unique $\mathcal{C}^{\infty}-$homomorphism:
$$\widetilde{f_{ST}}: A\{ S^{-1}\} \to B\{ T^{-1}\}$$
such that the following diagram commutes:
$$\xymatrixcolsep{5pc}\xymatrix{
A \ar[d]_{f_{ST}} \ar[r]^{\eta_S} & A\{ S^{-1}\} \ar[d]^{\widetilde{f_{ST}}}\\
B \ar[r]^{\eta_T} & B\{ T^{-1}\}}$$

Let:

$$F((A,S) \stackrel{f_{ST}}{\rightarrow} (B,T)) := A\{ S^{-1}\} \stackrel{\widetilde{f_{ST}}}{\rightarrow} B\{ T^{-1}\}.$$

\begin{theorem}The map:
$$\begin{array}{cccc}
    F: & \mathcal{C}^{\infty}_2 & \rightarrow & \mathcal{C}^{\infty}{\rm \bf Rng} \\
     & (A,S) & \mapsto & A\{ S^{-1}\}\\
     & \xymatrix{(A,S) \ar[r]^{f_{ST}} & (B,T)} & \mapsto & \xymatrix{A\{ S^{-1}\} \ar[r]^{\widetilde{f_{ST}}} & B\{ T^{-1}\}} \\
  \end{array}$$
is a functor.
\end{theorem}

Now, given any $f: (A,S) \to (B,B^{\times})$ - that is, given any $f: A \to B$ such that $f[S]\subseteq B^{\times}$, there is a unique $\widetilde{f}: A\{ S^{-1}\} \to B$ such that $\widetilde{f}\circ \eta_S = f$.  Conversely, given any $\widetilde{g}: (A,S) \to \jmath(B)$, we can take $\widetilde{g}\circ \eta_S : A \to B$, which is such that $(\widetilde{g}\circ \eta_S)[S] \subseteq B^{\times}$. This process yields a natural bijection between ${\rm Hom}_{\mathcal{C}^{\infty}{\rm \bf Rng}}(F(A,S),B)$ and ${\rm Hom}_{\mathcal{C}^{\infty}_2}((A,S),\jmath(B))$, and thus $F \dashv \jmath$.\\

\begin{proposition}\label{par}Let $A$ and $B$ be two $\mathcal{C}^{\infty}-$rings and $S \subseteq A$ and $f: A \to B$ a $\mathcal{C}^{\infty}-$homomorphism. By the universal property of $\eta_S : A \to A\{ S^{-1}\}$ we have a unique $\mathcal{C}^{\infty}-$homomorphism $f_S : A\{ S^{-1}\} \to B\{ f[S]^{-1}\}$ such that the following square commutes:
$$\xymatrixcolsep{5pc}\xymatrix{
A \ar[r]^{\eta_S} \ar[d]_{f} & A\{ S^{-1}\} \ar[d]^{\exists ! f_S}\\
B \ar[r]_{\eta_{f[S]}} & B\{ f[S]^{-1}\}
}.$$

The diagram:
%$$\xymatrixcolsep{5pc}\xymatrix{
%B \ar[dr]^{\eta_{f[S]}}& \\
% & B \{ f[S]^{-1}\}\\
%A\{ S^{-1}\} \ar[ur]_{f_{S}}
%}$$
$$\xymatrixcolsep{5pc}\xymatrix{
B \ar[r]^{\eta_{f[S]}} & B\{ f[S]^{-1}\} & \ar[l]_{f_S}  A\{ S^{-1}\}
}$$

is a pushout of the diagram:
$$\xymatrixcolsep{5pc}\xymatrix{
 A\{ S^{-1}\} &\ar[l]_{\eta_S}  A \ar[r]^{f} & B
}$$
\end{proposition}

\begin{corollary}\label{papel}The following rectangle is a pushout:

$$\xymatrixcolsep{5pc}\xymatrix{
A  \ar@{-->}[dr]^{\eta_{a_i \cdot a_j}} \ar[r]^{\eta_{a_i}} \ar[d]_{\eta_{a_j}} & A\{ {a_i}^{-1}\} \ar[d]\\
A\{ {a_j}^{-1}\} \ar[r] & A\{(a_i \cdot a_j)^{-1}\}
}$$

\end{corollary}

\begin{proposition}Let $(A,S)$ and $(B,T)$ be any two pairs in $\mathcal{C}^{\infty}_2$ and let $\jmath_{ST} : (A,S) \to (B,T)$ be a  $\mathcal{C}^{\infty}-$monomorphism such that $\jmath_{ST}[S^{\infty-{\rm sat}}] = T^{\infty-{\rm sat}}$. Under these circumstances,  $ \widetilde{\jmath_{ST}}: A\{ S^{-1}\} \to B\{ T^{-1}\}$ is also a $\mathcal{C}^{\infty}-$monomorphism.
\end{proposition}

\begin{theorem}Let $(A_1,S_1)$ and $(A_2,S_2)$ be two pairs in ${\rm Obj}\,(\mathcal{C}^{\infty}_2)$, $\jmath_1: S_1 \hookrightarrow A_1$ and $\jmath_2 : S_2 \hookrightarrow A_2$ be the inclusion maps and let:
$$\xymatrixcolsep{5pc}\xymatrix{
A_1 \ar[r]^{k_1} & A_1\otimes_{\infty} A_2  & \ar[l]_{k_2} A_2 
}$$

be the coproduct of $A_1$ and $A_2$ in $\mathcal{C}^{\infty}{\rm \bf Rng}$. We have:
$$(A_1 \otimes_{\infty} A_2)\{ (k_1[S_1]\cup k_2[S_2])^{-1}\} \cong A_1\{ {S_1}^{-1}\}\otimes_{\infty} A_2 \{ {S_2}^{-1} \},$$
that is,
$$F(A_1,S_1)\otimes_{\infty} F(A_2,S_2) \cong F(A_1 \otimes_{\infty} A_2, k_1[S_1]\cup k_2[S_2]).$$
\end{theorem}

\begin{remark}With an analogous proof, we can obtain a more general result: given any set of pairs, $\{ (A_i,S_i) | i \in I\}$, we have:

$$\bigotimes_{\infty}A_i\{ {S_i}^{-1}\} \cong \left( \bigotimes_{\infty} A_i\right)\left\{ \left(\bigcup_{i \in I}k_i[S_i]\right)^{-1}\right\}$$

where $k_i: A_i \rightarrow \otimes_{\infty}A_i$, for each $i$, is the canonical $\mathcal{C}^{\infty}-$homomorphism into the coproduct.
\end{remark}

\begin{theorem}\label{Newt}Let $B$ be the directed colimit of a system $\{ A_{\ell} \stackrel{t_{\ell j}}{\rightarrow} A_j | \ell, j \in I\}$ of $\mathcal{C}^{\infty}-$rings, that is,

$$\xymatrixcolsep{5pc}\xymatrix{
 & B =\varinjlim_{\ell \in I} A_{\ell} & \\
A_i \ar[ur]_{t_i} \ar[rr]^{t_{ij}} & & A_j \ar[ul]^{t_j}
}$$

is a limit co-cone. Given any $u \in B$, there are $j \in I$ and $u_j \in A_j$ such that $t_j(u_j) = u$.\\

Under those circumstances, we have:

$$\varinjlim_{k \geq j} A_{k}\{ u_{k}^{-1} \} \cong B \{ u^{-1}\}$$
where $u_k = t_{jk}(u_j)$.
\end{theorem}

\begin{theorem}\label{2137}Let $(I, \leq)$ be a directed poset and let $(A,S)$ be the directed colimit of a system $\{ (A_{\ell}, S_{\ell}) \stackrel{t_{\ell j}}{\rightarrow} (A_j, S_j) | \ell, j \in I\}$ of $\mathcal{C}^{\infty}_{2}$, that is,

$$\xymatrixcolsep{5pc}\xymatrix{
 & (A,S) = \varinjlim_{\ell \in I} (A_{\ell}, S_{\ell}) & \\
(A_i, S_i) \ar[ur]_{t_i} \ar[rr]^{t_{ij}} & & (A_j, S_j) \ar[ul]^{t_j}
}$$

We have:

$$F\left( \varinjlim_{\ell \in I} (A_{\ell}, S_{\ell})\right) \cong \varinjlim_{\ell \in I} F(A_{\ell}, S_{\ell}) = \varinjlim_{\ell \in I} A_{\ell}\{ S_{\ell}^{-1}\}$$

$$\varinjlim_{i \in I} A_i\{ S_i^{-1}\} \cong \left( \varinjlim_{i \in I} A_i\right)\left\{ \left(\varinjlim_{i \in I}S_i\right)^{-1}\right\}$$
\end{theorem}

As a consequence of the above proposition we have:

\begin{corollary}In the same context of \textbf{Theorem \ref{2137}}, we have:
  $$\left( \varinjlim_{i \in I} S_i\right)^{\infty-{\rm sat}} = \varinjlim_{i \in I} S_i^{\infty-{\rm sat}} \subseteq \varinjlim_{i \in I}A_i$$
\end{corollary}

Next we present some theorems about quotients of $\mathcal{C}^{\infty}-$rings by their ideals and some  relationships between quotients and $\mathcal{C}^{\infty}-$rings of fractions. We are going to prove that, in the category of $\mathcal{C}^{\infty}-$rings, taking quotients and taking rings of fractions are constructions which commute in a sense that will become clear in the following considerations.\\

Let $A$ be a $\mathcal{C}^{\infty}-$ring, $I$ be any of its ideals, $S \subseteq A$ be any subset of $A$ and consider the canonical map of the ring of fractions of $A$ with respect to $S$:
$$\eta_S : A \to A\{ S^{-1}\},$$
and consider $\widetilde{I_S} = \langle \eta_S[I]\rangle$.\\

Since $(q_{\widetilde{I_S}}\circ \eta_S)[I] \subseteq \widetilde{I_S}$, then by the \textbf{Theorem of Homomorphism} there exists a unique $\eta_{SI} : \dfrac{A}{I} \to \dfrac{A\{ S^{-1}\}}{\widetilde{I_S}}$ such that:

$$\xymatrix{
A \ar[r]^{\eta_S} \ar[d]^{q_I} & A\{ S^{-1}\} \ar[d]^{q_{\widetilde{I_S}}}\\
\dfrac{A}{I} \ar@{.>}[r]^{\eta_{SI}} & \dfrac{A\{ S^{-1}\}}{\widetilde{I_S}}
}$$

commutes.

\begin{proposition}The morphism $\eta_{SI}: \dfrac{A}{I} \to \dfrac{A\{ S^{-1}\}}{\widetilde{I_S}}$ has the same universal property which identifies the ring of quotients of $q_I[S] = S+I$.
\end{proposition}

\begin{proposition}Let $A$ be a $\mathcal{C}^{\infty}-$ring, $I$ any of its ideals and $S \subseteq A$ any of its subsets. Since $(\eta_{S+I} \circ q_I)[A^{\times}] \subseteq \left( \dfrac{A}{I}\lbrace \left( S+I\right)^{-1}\rbrace\right)^{\times}$, there is a unique $\widehat{q}_{IS}$ such that the following diagram commutes:
$$\xymatrixcolsep{5pc}\xymatrix{
A \ar[r]^{\eta_S} \ar[d]^{q_I} & A\{ S^{-1}\} \ar@{.>}[d]^{\widehat{q}_{IS}}\\
\dfrac{A}{I} \ar[r]^{\eta_{S+I}} & \left( \dfrac{A}{I} \right) \lbrace \left( S+I\right)^{-1}\rbrace
}$$

The $\mathcal{C}^{\infty}-$homomorphism $\mu^{-1}: \dfrac{A\{S^{-1} \}}{\widetilde{I_S}} \to \dfrac{A}{I} \lbrace \left( S+I \right)^{-1}\rbrace$ is an isomorphism.
\end{proposition}

\begin{corollary}\label{Jeq}Let $A$ be a $\mathcal{C}^{\infty}-$ring, $I$ any of its ideals and $S \subseteq A$ any of its subsets. There are unique isomorphisms $\mu: \dfrac{A}{I} \lbrace \left( S + I\right)^{-1}\rbrace \to \dfrac{A\{ S^{-1}\}}{\langle \eta_S[I] \rangle}$ such that the following pentagons commute:
$$\begin{array}{lr}
\xymatrixcolsep{3pc}\xymatrix{
A \ar[d]^{q_I} \ar[r]^{\eta_S} & A\{ S^{-1}\} \ar[r]^{q_{\widetilde{I_S}}} & \dfrac{A\{ S^{-1}\}}{\langle \eta_S[I] \rangle} \\
\dfrac{A}{I} \ar[r]^{\eta_{S+I}} & \dfrac{A}{I} \lbrace \left( S + I\right)^{-1}\rbrace \ar[ur]^{\mu}
}&
\xymatrixcolsep{3pc}\xymatrix{
A \ar[d]^{q_I} \ar[r]^{\eta_S} & A\{ S^{-1}\} \ar[r]^{q_{\widetilde{I_S}}} & \dfrac{A\{ S^{-1}\}}{\langle \eta_S[I] \rangle} \ar[dl]^{\mu^{-1}} \\
\dfrac{A}{I} \ar[r]^{\eta_{S+I}} & \dfrac{A}{I} \lbrace \left( S+I\right)^{-1}\rbrace
}
\end{array}$$
\end{corollary}

\section{Distinguished Classes of $\mathcal{C}^{\infty}-$rings} \label{classes-sec}

\hspace{0.5cm}In this section we present some distinguished classes of $\mathcal{C}^{\infty}-$rings, such as $\mathcal{C}^{\infty}-$fields, $\mathcal{C}^{\infty}-$domains, $\mathcal{C}^{\infty}-$local rings and reduced $\mathcal{C}^{\infty}-$rings. The concept of von Neumann regular $\ci$-ring is explored in details in \cite{VNBM}).

In \cite{rings1}, we find definitions of $\mathcal{C}^{\infty}-$fields, $\mathcal{C}^{\infty}-$domains and $\mathcal{C}^{\infty}-$local rings. We describe these concepts in the following:

\begin{definition}\label{distintos}Let ${\rm \bf CRing}$ be the category of all commutative unital rings, and consider the forgetful functor $\widetilde{U}: \mathcal{C}^{\infty}{\rm \bf Rng} \to {\rm \bf CRing}$. We say that a $\mathcal{C}^{\infty}-$ring $A$ is:
\begin{itemize}
  \item{a \index{$\mathcal{C}^{\infty}-$field}\textbf{$\mathcal{C}^{\infty}-$field} whenever $\widetilde{U}(A) \in {\rm Obj}\,({\rm \bf CRing})$ is a field;}
  \item{a \index{$\mathcal{C}^{\infty}-$domain}\textbf{$\mathcal{C}^{\infty}-$domain} whenever $\widetilde{U}(A) \in {\rm Obj}\,({\rm \bf CRing})$ is a domain;}
  \item{a \index{local $\mathcal{C}^{\infty}-$ring}\textbf{local $\mathcal{C}^{\infty}-$ring} whenever $\widetilde{U}(A) \in {\rm Obj}\,({\rm \bf CRing})$ is a local ring.}
  %\item{a \textbf{von Neumann regular $\mathcal{C}^{\infty}-$ring} whenever $\widetilde{U}(A) \in {\rm Obj}\,({\rm \bf CRing})$ is a von Neumann regular ring.}
\end{itemize}
\end{definition}

Thus a $\mathcal{C}^{\infty}-$field, a $\mathcal{C}^{\infty}-$domain and a local $\mathcal{C}^{\infty}-$ring are $\mathcal{C}^{\infty}-$rings such that their underlying $\R-$algebras are  fields, domains and  local rings in the ordinary sense, respectively.\\

The next fact gives us a general method to obtain finitely generated local $\mathcal{C}^{\infty}-$rings using prime filters on the set of the closed parts of $\mathbb{R}^n$.\\

From now on, we are going to denote the set of all prime proper $\mathcal{C}^{\infty}-$radical ideals of a $\mathcal{C}^{\infty}-$ring $A$ by ${\rm Spec}^{\infty}(A)$.\\

\begin{proposition}\label{quetiro}Let $A$ be a $\mathcal{C}^{\infty}-$ring and $\p \in {\rm Spec}^{\infty}\,(A)$. Then $A_{\p} := A\{ {A \setminus \p}^{-1}\}$ is a local $\mathcal{C}^{\infty}-$ring whose unique maximal ideal is given by $\mathfrak{m}_{\p} = \left\{ \frac{\eta_{A\setminus \p}(x)}{\eta_{A\setminus \p}(y)} |  (x \in \p) \& (y \in A \setminus \p)\right\}$.
\end{proposition}
\begin{proof}
See p. 69 of \cite{BMTSCA}.
\end{proof}

\begin{proposition}\label{sizi}Let $A$ be a $\mathcal{C}^{\infty}-$ring. The following assertions are equivalent.
\begin{itemize}
\item[i)]{$A$ is a $\mathcal{C}^{\infty}-$field;}
\item[ii)]{For every subset $S \subseteq A\setminus \{0\}$, the canonical map:
$$\begin{array}{cccc}
{\rm Can}_S : & A & \rightarrow & A\{ S^{-1} \}
\end{array}$$
is a $\mathcal{C}^{\infty}-$ring isomorphism;}
\item[iii)]{For any $a \in A\setminus \{ 0\}$, we have that:
$$\begin{array}{cccc}
{\rm Can}_a : & A & \rightarrow & A\{ a^{-1}\}
\end{array}$$ is a $\mathcal{C}^{\infty}-$isomorphism.}
\end{itemize}
\end{proposition}
\begin{proof} See the proof of \textbf{Proposition 21}, page 70 of \cite{BMTSCA}.
\end{proof}

In ordinary Commutative Algebra, given an element $x$ of a ring $R$, we say that $x$ is a nilpotent infinitesimal if and only if there is some $n \in \mathbb{N}$ such that $x^n=0$. Let $A$ be a  $\mathcal{C}^{\infty}-$ring and $a \in A$. D. Borisov and K. Kremnizer in \cite{Borisov} call $a$ an $\infty-$infinitesimal if, and only if $A\{ a^{-1}\} \cong 0$. The next definition describes the notion of a $\mathcal{C}^{\infty}-$ring being free of $\infty-$infinitesimals - which is analogous to the notion of ``reducedness'', of a commutative ring.\\

\begin{definition}A $\mathcal{C}^{\infty}-$ring $A$ is \index{$\mathcal{C}^{\infty}-$reduced}\textbf{$\mathcal{C}^{\infty}-$reduced} if, and only if, $\sqrt[\infty]{(0)} = (0)$.
\end{definition}

\begin{proposition}Every $\mathcal{C}^{\infty}-$field, $A$, is a $\mathcal{C}^{\infty}-$reduced $\mathcal{C}^{\infty}-$ring.
\end{proposition}
\begin{proof}
Suppose $A$ is a $\mathcal{C}^{\infty}-$field. For every $a \in A$, we have $\dfrac{A}{(0)}\{ (a+(0))^{-1}\} \cong A\{ a^{-1}\}$, so $A\{ a^{-1}\} \cong A$ occurs if, and only if, $a \neq 0$. Hence $\sqrt[\infty]{(0)}=(0)$
\end{proof}

\begin{theorem}\label{Cindy}Let $A$ be a $\mathcal{C}^{\infty}-$ring, $S \subseteq A$, and $I \subset A$ any ideal. Then:
$$\langle \eta_S[I]\rangle = \left\{ \frac{\eta_S(b)}{\eta_S(d)} | b \in I \& d \in S^{\infty-{\rm sat}}\right\}$$
\end{theorem}
\begin{proof}
Given $h \in \langle \eta_S[I]\rangle$, there are $n \in \mathbb{N}$, $b_1, \cdots, b_n \in I$ and $\alpha_1, \cdots, \alpha_n \in A\{S^{-1}\}$ such that:

$$h = \sum_{i=1}^{n} \alpha_i \cdot \eta_{S}(b_i) $$

For each $i \in \{ 1, \cdots, n\}$ there are $c_i \in A$ and $d_i \in S^{\infty-{\rm sat}}$ such that:

$$\alpha_i \cdot \eta_S(d_i) = \eta_S(c_i)$$

so

$$h = \sum_{i=1}^{n} \alpha_i \cdot \eta_{S}(b_i) = \sum_{i=1}^{n} \dfrac{1}{\eta_S(d_i)}\cdot \eta_S(c_i)\cdot \eta_S(b_i),$$

and denoting $b_i':= c_i \cdot b_i \in I$, we get:

$$h = \sum_{i=1}^{n} \dfrac{\eta_S(b_i')}{\eta_S(d_i)}$$

For each $i=1, \cdots, n$, let
$$b_i'' = b_i' \prod_{j \neq i} d_j,$$
so
$$h = \dfrac{\eta_S \left( \prod_{j \neq 1}d_j\right) \eta_S(b_1')}{\eta_S(d_1 \cdots d_n)} + \dfrac{\eta_S \left( \prod_{j \neq 2}d_j \right) \eta_S(b_2')}{\eta_S(d_1 \cdots d_n)} + \cdots + \dfrac{\eta_S \left( \prod_{j \neq n}d_j\right) \eta_S(b_n')}{\eta_S(d_1 \cdots d_n)}$$
hence

$$h \cdot \eta_S(d_1 \cdots d_n)  =\eta_S \left( \prod_{j \neq 1}d_j\right) \eta_S(b_1') + \eta_S \left( \prod_{j \neq 2}d_j\right) \eta_S(b_2')+ \cdots + \eta_S \left( \prod_{j \neq n}d_j\right) \eta_S(b_n').$$

Let $b_i'' := \left( \prod_{j \neq i}d_i\right)\cdot \overbrace{b_i'}^{\in I} \in I$, so we have $h \cdot \eta_S(d_1 \cdots d_n) = \sum_{i=1}^{n} \eta_S(b_i'') = \eta_S\left( \sum_{i=1}^{n}b_i''\right)$

Since $\eta_S(d_1 \cdots d_n) \in A\{ S^{-1}\}^{\times}$, $d = d_1\cdots d_n \in S^{\infty-{\rm sat}}$, so taking $b = \sum_{i=1}^{n} b_i''$, we have $b \in I$, since it is a sum of elements of $I$, we can write:

$$h \cdot \eta_S(d) = \eta_S(b),$$

and $h = \dfrac{\eta_S(b)}{\eta_S(d)}$, with $b \in I$ and $d \in S^{\infty-{\rm sat}}$.

The other way round is immediate.
\end{proof}

As a consequence of \textbf{Theorem \ref{Cindy}}, we have:

\begin{corollary}If $A$ is a reduced $\mathcal{C}^{\infty}-$ring, then we have:
$$(\forall a \in A)((0 \in \{ a\}^{\infty-{\rm sat}}) \leftrightarrow (a=0))$$
\end{corollary}

Now, given  $I \subsetneq \mathcal{C}^{\infty}(\mathbb{R}^n)$ any ideal,
$$\widehat{I} = \{ A \subseteq \mathbb{R}^n | A \, \mbox{is closed and}\, (\exists f \in  I)(A = Z(f))\}$$ is a filter on the set of all the closed subsets of $\mathbb{R}^n$. Moreover, $I$ is a proper ideal if, and only if, $\widehat{I}$ is a proper filter (see p. 72 of \cite{BMTSCA}). Also, given any filter $\mathcal{F}$ on the set of all the closed subsets of $\mathbb{R}^n$, we define
$$\check{\mathcal{F}} = \{ f \in \mathcal{C}^{\infty}(\mathbb{R}^n) | Z(f) \in \mathcal{F} \},$$
which is an ideal of $\mathcal{C}^{\infty}(\mathbb{R}^n)$. Moreover, $\mathcal{F}$ is a proper filter if, and only if, $\check{\mathcal{F}}$ is a proper ideal. (see p. 73 of \cite{BMTSCA}).\\

Let $I \subseteq \mathcal{C}^{\infty}(\mathbb{R}^n)$ be any ideal. Then we have, by definition, $\widehat{I} = \{ A \subseteq \mathbb{R}^n | (\exists h \in I)(A = Z(h)) \}$, and therefore
$$\check{\widehat{I}} = \{ g \in \mathcal{C}^{\infty}(\mathbb{R}^n) | (\exists h \in I)(Z(g) = Z(h)) \}.$$

Let $\mathcal{F}$ be a proper filter on the set of all closed subsets of $\mathbb{R}^n$. Then, since every closed subset $B \subseteq \mathbb{R}^n$ is a zero set of some smooth function, we have:

\begin{multline*}\widehat{\check{\mathcal{F}}} = \{ A \subseteq \mathbb{R}^n | (\exists g \in \check{\mathcal{F}})(A=Z(g))\} = \\
= \{ A \subseteq \mathbb{R}^n | (\exists g \in \check{\mathcal{F}})((Z(g) \in \mathcal{F})\&(A=Z(g)))\} = \mathcal{F}
\end{multline*}

%\begin{proof}
%First we show that:
%$$\{ g \in \mathcal{C}^{\infty}(\mathbb{R}^n) | (\exists h \in I)(Z(g) \supseteq Z(h)) \} \subseteq \check{\widehat{I}}.$$

%Given $g$ such that $(\exists h \in I)(Z(g) \supseteq Z(h))$, consider the smooth function $h' = h \cdot \chi_{U_g}$, where $U_g = \mathbb{R}^n \setminus Z(g)$. We have:
%$$Z(h') = Z(h)\cup Z(\chi_{U_g}).$$
%Now, since $Z(h) \subset Z(g)$, we have $Z(h') = Z(h)\cup Z(\chi_{U_g}) = Z(h) \cup Z(g) = Z(g)$. Since $h \in I$, $h' = h \cdot \chi_{U_g} \in I$, and therefore $g \in \check{\widehat{I}}$

%\end{proof}

As a consequence of the discussion above, we have:

\begin{proposition}\label{49}Let $I \subseteq \mathcal{C}^{\infty}(\mathbb{R}^n)$ be an ideal. Then:
$$\check{\widehat{I}} = \{ g \in \mathcal{C}^{\infty}(\mathbb{R}^n) | (\exists h \in I)(Z(g) = Z(h))\} = \sqrt[\infty]{I}$$
In particular, the $\mathcal{C}^{\infty}-$radical of an ideal $I$ of the free $\mathcal{C}^{\infty}-$ring on finitely many generators, $\mathcal{C}^{\infty}(\mathbb{R}^n)$, is again an ideal.
\end{proposition}
\begin{proof}See p. 73 of \cite{BMTSCA}.
%Let $g \in \check{\widehat{I}}$, so there is some $h \in I$ such that $Z(g) = Z(h)$, so \textbf{a fortiori}, $Z(g) \supseteq Z(h)$, and $(h \displaystyle\upharpoonright_{U_g})$ is invertible in $\mathcal{C}^{\infty}(U_g)$, where $U_g = \mathbb{R}^n \setminus Z(g)$ (see \textbf{Proposition \ref{pani}}). Since there exists $h \in I$ such that $h \in (\mathcal{C}^{\infty}(\mathbb{R}^n)\{ g \}^{-1})^{\times}$ it follows that $g \in \sqrt[\infty]{I}$.\\

%Conversely, given $g \in \sqrt[\infty]{I}$, there must be some $h \in I$ such that $h \in (\mathcal{C}^{\infty}(\mathbb{R}^n)\{g^{-1}\})^{\times}$. But such an $h$ must not be zero when restricted to $U_g$ (otherwise $h$ would not be invertible), so $U_h \supseteq U_g$, and therefore $Z(g) \supseteq Z(h)$. Since $\widehat{I}$ is a filter, $Z(h) \in \widehat{I}$ and $Z(g) \supseteq Z(h)$, it follows that $Z(g) \in \widehat{I}$. But if $Z(g) \in \widehat{I}$ then there is $h' \in I$ such that $Z(g) = Z(h')$ so $g \in \check{\widehat{I}}$.
\end{proof}

The following proposition is a consequence of \textbf{Proposition \ref{49}}

\begin{corollary}\label{Guia}Let $A = \mathcal{C}^{\infty}(\mathbb{R}^n)$ be a finitely generated free $\mathcal{C}^{\infty}-$ring. $I \subseteq \mathcal{C}^{\infty}(\mathbb{R}^n)$ is a $\mathcal{C}^{\infty}-$radical ideal, that is, $\sqrt[\infty]{I}=I$, if, and only if:
$$(\forall g \in \mathcal{C}^{\infty}(\mathbb{R}^n))((g \in I) \leftrightarrow (\exists f \in I)(Z(f)=Z(g))).$$
\end{corollary}

The following proposition is a consequence of a comment made by Moerdijk and Reyes in the page 330 of \cite{rings2}:

\begin{corollary}\label{Gaia}Let $\mathcal{C}^{\infty}(\mathbb{R}^n)$ be the free $\mathcal{C}^{\infty}-$ring and let $I \subseteq \mathcal{C}^{\infty}(\mathbb{R}^n)$ be a finitely generated ideal, that is,
$$I = \langle g_1, \cdots, g_k\rangle,$$
for some $g_1, \cdots, g_k \in \mathcal{C}^{\infty}(\mathbb{R}^n)$.\\

$I$ is a $\mathcal{C}^{\infty}-$radical ideal if, and only if:

$$(\forall x \in Z(g_1, \cdots, g_k))(f(x)=0) \to (f \in I)$$
\end{corollary}
\begin{proof}
Suppose $I = \sqrt[\infty]{I}$ and let $f \in \mathcal{C}^{\infty}(\mathbb{R}^n)$ be such that:
$$(\forall x \in Z(g_1, \cdots, g_k))(f(x)=0).$$

We have $g=g_1^2 + \cdots + g_k^2 \in I$ such that $Z(g)\subseteq Z(f)$, so $Z(f)\in \widehat{I}$ (since $\widehat{I}$ is a filter) and by \textbf{Proposition \ref{49}}, $f \in \check{\widehat{I}} = \sqrt[\infty]{I}=I$.\\

Now, suppose

$$(\forall x \in Z(g_1, \cdots, g_k))(f(x)=0) \to (f \in I)$$

Given $h \in \sqrt[\infty]{I} = \check{\widehat{I}}$, we have $Z(h) \in \widehat{I}$, so there exists some $g \in I$ such that $Z(h)=Z(g)$, thus $(\forall x \in Z(g_1, \cdots, g_k))(h(x)=0)$. By hypothesis, this means that $h \in I$, so $\sqrt[\infty]{I} \subseteq I$. Since $I \subseteq \sqrt[\infty]{I}$ always holds, it follows that $I$ is a $\mathcal{C}^{\infty}-$radical ideal.
\end{proof}

By \textbf{Proposition 4.5} of \cite{MSIA}, it follows that any finitely generated $\mathcal{C}^{\infty}-$radical ideal $I$ of a free $\mathcal{C}^{\infty}-$ring on $n$ generators, $\mathcal{C}^{\infty}(\mathbb{R}^n)$
we have:

$$(\forall x \in Z(I))(f\upharpoonright_x \in I\upharpoonright_x \to f \in I),$$

where $Z(I):= \bigcap_{f \in I} Z(f).$

As a consequence, we have a particular version of the weak \textbf{Nullstellensatz}:

\begin{proposition}\label{Nullstellensatz}For any finitely generated $\mathcal{C}^{\infty}-$radical ideal $I$ of $\mathcal{C}^{\infty}(\mathbb{R}^n)$, we have:
$$1 \in I \iff Z(I)=\varnothing$$
\end{proposition}
\begin{proof}
Let $I = \langle g_1, \cdots, g_k\rangle$ be a $\mathcal{C}^{\infty}-$radical ideal of $\mathcal{C}^{\infty}(\mathbb{R}^n)$, so
$$(\forall x \in Z(g_1, \cdots, g_k))(f(x)=0) \to (f \in I)$$
Suppose $Z(g_1, \cdots,g_k)=\varnothing$, so we have, for each $f \in \mathcal{C}^{\infty}(\mathbb{R}^n)$:

$$\varnothing = Z(g_1, \cdots, g_k) \subseteq Z(f),$$

so

$$(\forall x \in Z(g_1\cdots g_k))(f(x)=0)$$

Thus $\mathcal{C}^{\infty}(\mathbb{R}^n) \subseteq I$ and, in particular, $1 \in I$.

On the other hand, if $1 \in I$, since $1^{\dashv}[\{ 0\}] = \varnothing \in \{ f^{\dashv}[\{ 0\}] \mid f \in I \}$, so $Z(I)=\bigcap \{ f^{\dashv}[\{ 0\}] \mid f \in I\} = \varnothing$.
\end{proof}

%\begin{example}As a consequence of the fact pointed out by Moerdijk and Reyes in \cite{rings2}, any countably generated prime ideal of $\mathcal{C}^{\infty}(\mathbb{R}^n)$, for $n \in \mathbb{N}$, is a $\mathcal{C}^{\infty}-$radical ideal.
%\end{example}

Let $\mathfrak{F}$ be the set of all the filters on the closed parts of $\mathbb{R}^n$ and $\mathfrak{I}$ be the set of all the ideals of $\mathcal{C}^{\infty}(\mathbb{R}^n)$. We have, so far, established that for every ideal $I \subseteq \mathcal{C}^{\infty}(\mathbb{R}^n)$ we have $\sqrt[\infty]{I} = \check{\widehat{I}}$ and for every filter $\mathcal{F} \in \mathfrak{F}$ we have $\widehat{\check{\mathcal{F}}}$.\\

Consider the following diagram:

$$\xymatrix{\mathfrak{F} \ar@/^/[rr]|{\vee} && \ar@/^/[ll]|{\wedge} \mathfrak{I}}$$

In what follows, we show that $\wedge \vdash \vee$ is a \index{Galois connection}Galois connection.

\begin{proposition}\label{pseu}The  adjunction $\wedge \vdash \vee$ is a covariant Galois connection between the posets $(\mathfrak{F}, \subseteq)$ and $(\mathfrak{I}, \subseteq)$, i.e.,
\begin{itemize}
  \item[a)]{Given $\mathcal{F}_1, \mathcal{F}_2 \in \mathfrak{F}$ such that $\mathcal{F}_1 \subseteq \mathcal{F}_2$ then $\check{\mathcal{F}_1} \subseteq \check{\mathcal{F}_2}$;}
  \item[b)]{Given $I_1, I_2 \in \mathfrak{I}$ such that $I_1 \subseteq I_2$ then $\widehat{I}_1 \subseteq \widehat{I}_2$;}
  \item[c)]{For every $\mathcal{F} \in \mathfrak{F}$ and every $I \in \mathfrak{I}$ we have:
  $$\widehat{I} \subseteq \mathcal{F} \iff I \subseteq \check{\mathcal{F}}.$$}
\end{itemize}
\end{proposition}
\begin{proof} See p. 75 of \cite{BMTSCA}.
\end{proof}

\begin{proposition}\label{Garland}Let $A = \mathcal{C}^{\infty}(\R^n)$ for some $n \in \mathbb{N}$. The Galois connection $\wedge \vdash \vee$ establishes a bijective correspondence between:
\begin{itemize}
  \item[(a)]{proper filters of $(\mathfrak{F}, \subseteq)$ and proper ideals of $(\mathfrak{I}, \subseteq)$;}
  \item[(b)]{maximal filters of $(\mathfrak{F}, \subseteq)$ and maximal ideals of $(\mathfrak{I}, \subseteq)$;}
  \item[(c)]{prime filters of $(\mathfrak{F}, \subseteq)$ and prime ideals of $(\mathfrak{I}, \subseteq)$.}
  \item[(d)]{filters on the closed parts of $\mathbb{R}^n$, $\mathfrak{F}$, and the set of all $\mathcal{C}^{\infty}-$radical ideals of $\mathcal{C}^{\infty}(\mathbb{R}^n)$, $\mathfrak{I}^{\infty} = \{ I \subseteq \mathcal{C}^{\infty}(\mathbb{R}^n) | \sqrt[\infty]{I} = I \}$.}
\end{itemize}
\end{proposition}
\begin{proof}See p. 75 of \cite{BMTSCA}.
\end{proof}

\begin{corollary} \label{Piqueri}Let $A = \mathcal{C}^{\infty}(\R^n)$ and $\p$ be a prime ideal of $A$. Under those circumstances $\sqrt[\infty]{\p}$ is a prime ideal.
\end{corollary}
\begin{proof}
By the \textbf{Corollary \ref{Garland}}, since $\p$ is a prime ideal, the filter associated with $\p$, $\widehat{\p}$ is a prime filter. Again by \textbf{Corollary \ref{Garland}}, it follows that $\check{\widehat{\p}} = \sqrt[\infty]{\p}$ is a prime ideal.
\end{proof}

\begin{corollary}\label{tumbleweed}The composite $\vee \circ \wedge : \mathfrak{I} \to \mathfrak{I}$ is a \index{closure operator}closure operator, thus we have the following three properties:
\begin{itemize}
  \item[(a)]{$I \subseteq \check{\widehat{I}}$;}
  \item[(b)]{$I_1 \subseteq I_2 \Rightarrow \check{\widehat{I_1}} \subseteq \check{\widehat{I_2}}$;}
  \item[(c)]{If $J = \check{\widehat{I}}$ then $\check{\widehat{J}} = J$.}
\end{itemize}
\end{corollary}

\begin{theorem}\label{yellow}Let $I, I_1, I_2 \subseteq \mathcal{C}^{\infty}(\mathbb{R}^n)$ be ideals. Then:
\begin{itemize}
  \item[(a)]{$\sqrt[\infty]{I}$ is an ideal of $\mathcal{C}^{\infty}(\mathbb{R}^n)$ and $I \subseteq \sqrt[\infty]{I}$;}
  \item[(b)]{$I_1 \subseteq I_2 \Rightarrow \sqrt[\infty]{I_1} \subseteq \sqrt[\infty]{I_2}$}
  \item[(c)]{$\sqrt[\infty]{\sqrt[\infty]{I}} = \sqrt[\infty]{I}$}
\end{itemize}
\end{theorem}
\begin{proof}
Ad (a):  From \textbf{Proposition \ref{49}} we have $\sqrt[\infty]{I} = \check{\widehat{I}}$, and from \textbf{Proposition \ref{47}} it follows that $\check{\widehat{I}}$ is an ideal. Moreover, since $\vee \circ \wedge$ is a closure operator, by item a) of \textbf{Corollary \ref{tumbleweed}} we have $I \subseteq \check{\widehat{I}} = \sqrt[\infty]{I}$.\\

Ad (b): From item (b) of \textbf{Corollary \ref{tumbleweed}}, $I_1 \subseteq I_2 \Rightarrow \check{\widehat{I_1}} \subseteq \check{\widehat{I_2}}$, so $\sqrt[\infty]{I_1} \subseteq \sqrt[\infty]{I_2}$.\\

Ad (c): From item (c) of \textbf{Corollary \ref{tumbleweed}} it follows immediately by the idempotence of $\vee \circ \wedge$ that $\sqrt[\infty]{\sqrt[\infty]{I}} = \sqrt[\infty]{I}$.
\end{proof}

The following definition will be helpful to prove that whenever $I$ is an ideal of \textbf{any} $\mathcal{C}^{\infty}-$ring $A$, then $\sqrt[\infty]{I}$ is also an ideal.\\

\begin{definition}Let $A$ be a $\mathcal{C}^{\infty}-$ring. We say that $A$ is \textbf{admissible} if for every ideal $I \subseteq A$, $\sqrt[\infty]{I}$ is an ideal in $A$.
\end{definition}

In what follows we sketch a proof that every $\mathcal{C}^{\infty}-$ring $A$ is admissible.\\

\begin{lemma}\label{l1}Let $A$ and $A'$ be  $\mathcal{C}^{\infty}-$rings such that $A \cong A'$. Then if $A$ is an admissible $\mathcal{C}^{\infty}-$ring, $A'$ is admissible too.
\end{lemma}

\begin{lemma}\label{l2}The free $\mathcal{C}^{\infty}-$ring in $n$ generators, $\mathcal{C}^{\infty}(\mathbb{R}^n)$, is an admissible $\mathcal{C}^{\infty}-$ring.
\end{lemma}
\begin{proof}It follows immediately from  \textbf{Proposition \ref{49}}.
\end{proof}

\begin{lemma}\label{l3}Let $A$ be an admissible $\mathcal{C}^{\infty}-$ring and $I,J \subset A$ be ideals such that $I \subseteq J$. Then $\dfrac{A}{J}$ is an admissible $\mathcal{C}^{\infty}-$ring and:
$$\sqrt[\infty]{I+J} = \sqrt[\infty]{I}+J$$

\end{lemma}

\begin{lemma}\label{l4}Let $\{ A_i \stackrel{h_{ij}}{\rightarrow} A_j \}$ be a filtered diagram of admissible $\mathcal{C}^{\infty}-$rings. Then $\varinjlim A_i$ is an admissible $\mathcal{C}^{\infty}-$ring.
\end{lemma}

\begin{theorem}\label{little}Every $\mathcal{C}^{\infty}-$ring is admissible.
\end{theorem}
\begin{proof}Let $A$ be any $\mathcal{C}^{\infty}-$ring. We know that every $\mathcal{C}^{\infty}-$ring is a filtered colimit of finitely presented $\mathcal{C}^{\infty}-$rings. Since every finitely presentable $\mathcal{C}^{\infty}-$ring is admissible, the result follows from \textbf{Lemma \ref{l4}}.\\

Let $B$ be any finitely presentable $\mathcal{C}^{\infty}-$ring. We know that there exist some $n \in \mathbb{N}$ and some ideal $J \subset \mathcal{C}^{\infty}(\mathbb{R}^n)$ such that $A \cong \dfrac{\mathcal{C}^{\infty}(\mathbb{R}^n)}{J}$. From \textbf{Lemma \ref{l1}}, if we prove that $\dfrac{\mathcal{C}^{\infty}(\mathbb{R}^n)}{J}$ is admissible, then it follows that $A$ is admissible.\\

From \textbf{Lemma \ref{l2}} we have that $\mathcal{C}^{\infty}(\mathbb{R}^n)$ is admissible, and from \textbf{Lemma \ref{l3}} it follows that $\dfrac{\mathcal{C}^{\infty}(\mathbb{R}^n)}{J}$ is admissible.
\end{proof}

Now we present some properties of taking the $\mathcal{C}^{\infty}-$radical of an ideal.\\

\begin{proposition}Let $A$ be a $\mathcal{C}^{\infty}-$ring, $I,J \subseteq A$ any of its ideals. Then:
\begin{itemize}
\item[(i)]{$I \subseteq J \Rightarrow \sqrt[\infty]{I} \subseteq \sqrt[\infty]{J}$}
\item[(ii)]{$I \subseteq \sqrt[\infty]{I}$}
\end{itemize}
\end{proposition}
\begin{proof}
Ad (i): Given $a \in \sqrt[\infty]{I}$, there is $b \in I$ such that $\eta_a(b) \in (A\{ a^{-1}\})^{\times}$. Since $I \subseteq J$, the same $b$ is a witness of the fact that $a \in \sqrt[\infty]{J}$, for $b \in J$ and $\eta_a(b) \in (A\{ a^{-1}\})$.\\

Ad (ii): We are going to show that $A \setminus \sqrt[\infty]{I} \subseteq A \setminus I$.\\

Given $a \in A \setminus \sqrt[\infty]{I}$, we have that $(\forall b \in I)(\eta_a(b) \notin (A\{ a^{-1}\})^{\times})$, so $\eta_a[I] \cap (A\{ a^{-1}\})^{\times} = \varnothing$. Since $\eta_a(a) \in (A\{ a^{-1}\})^{\times}$, it follows that $\eta_a(a) \notin \eta_a[I]$, so $a \notin I$.
\end{proof}

\begin{proposition}Let $B$ be a $\mathcal{C}^{\infty}-$ring and $J \subseteq B$ any of its ideals. We have the following equality:
$$\sqrt[\infty]{\{ 0_{\frac{B}{J}}\}} = \dfrac{\sqrt[\infty]{J}}{J}$$
\end{proposition}
\begin{proof}
$$a \in \sqrt[\infty]{J} \iff (\exists b \in J)(\eta(b) \in (B\{ a^{-1}\})^{\times}),$$
hence
\begin{multline*}a + J = \overline{a} \in \sqrt[\infty]{\{ 0_{\frac{B}{J}}\}} \iff (\exists \overline{b} \in \{ 0_{\frac{B}{J}}\})(\overline{\eta_a}(\overline{b}) \in \left( \frac{B}{J}\{(a + J)^{-1} \}\right)^{\times}) \iff\\
\iff \dfrac{B}{J}\{(a+J)^{-1} \} \cong \{ 0\} \iff \\
 \iff a \in \sqrt[\infty]{J}\iff  a+ J \in \dfrac{\sqrt[\infty]{J}}{J}
 \end{multline*}

Now, since $J \subseteq \sqrt[\infty]{J}$, if $a'+J = a+J$, then $a \in \sqrt[\infty]{J} \iff a' \in \sqrt[\infty]{J}$.
\end{proof}

\begin{corollary}Let $A$ be a $\mathcal{C}^{\infty}-$ring. We have:
\begin{itemize}
  \item[(a)]{An ideal $J \subseteq A$  is a $\mathcal{C}^{\infty}-$radical ideal if, and only if, $\dfrac{A}{J}$ is a $\mathcal{C}^{\infty}-$reduced $\mathcal{C}^{\infty}-$ring}
  \item[(b)]{A proper prime ideal $\mathfrak{p} \subseteq A$ is $\mathcal{C}^{\infty}-$radical if, and only if, $\dfrac{A}{\mathfrak{p}}$ is a $\mathcal{C}^{\infty}-$reduced $\mathcal{C}^{\infty}-$domain.}
\end{itemize}
\end{corollary}

\begin{proposition} \label{redi}
Let $A', B'$ be two $\mathcal{C}^{\infty}-$rings and $\jmath: A' \to B'$ be a monomorphism. If $B'$ is $\mathcal{C}^{\infty}-$reduced, then $A'$ is also $\mathcal{C}^{\infty}-$reduced.
\end{proposition}
\begin{proof}See p. 81 of \cite{BMTSCA}.
\end{proof}

%\begin{remark}\label{corpo} Up to $\mathcal{C}^{\infty}-$isomorphisms, the only finitely generated $\mathcal{C}^{\infty}-$field is $\mathbb{R}$. In fact, let $F = \dfrac{\mathcal{C}^{\infty}(\mathbb{R}^n)}{I}$ be a finitely generated $\mathcal{C}^{\infty}-$field, so $I$ must be a maximal ideal of $\mathcal{C}^{\infty}(\R^n)$ (otherwise the quotient would not even be an ordinary field). By \textbf{Proposition 3.6} of \cite{MSIA}, every maximal ideal $I$ of $\mathcal{C}^{\infty}(\R^n)$  has the form $I = \ker {\rm ev}_x $, where $x \in Z(I)$ and:
%$$\begin{array}{cccc}
%{\rm ev}_x: & \mathcal{C}^{\infty}(\R^n) & \rightarrow & \mathbb{R} \\
%            & f                                  & \mapsto     & f(x)
%\end{array}$$
%Since ${\rm ev}_x$ is surjective, by the Fundamental Theorem of the $\mathcal{C}^{\infty}-$Homomorphism we have:
%$$F = \dfrac{\mathcal{C}^{\infty}(\R^n)}{I} = \dfrac{\mathcal{C}^{\infty}(\R^n)}{\ker {\rm ev}_x} \cong \R$$
%\end{remark}

As a consequence of \textbf{Proposition \ref{redi}} and of \textbf{Proposition \ref{alba}}, we have the following:

\begin{corollary}\label{medusa}Every $\mathcal{C}^{\infty}-$subring of a $\mathcal{C}^{\infty}-$field is a $\mathcal{C}^{\infty}-$reduced $\mathcal{C}^{\infty}-$domain.
\end{corollary}

Next we show that the directed limit of reduced $\mathcal{C}^{\infty}-$rings is also reduced.\\

\begin{proposition}\label{Xango}Let $(I, \leq)$ be a directed set and suppose that $\{A_i\}_{i \in I}$ is a directed family of $\mathcal{C}^{\infty}-$reduced $\mathcal{C}^{\infty}-$rings. Then
$$B = \varinjlim_{i \in I} A_i$$
is a $\mathcal{C}^{\infty}-$reduced $\mathcal{C}^{\infty}-$ring.
\end{proposition}
\begin{proof}See  p. 82 of \cite{BMTSCA}.
\end{proof}

\begin{theorem}\label{preim}Let $A$ and $B$ be two $\mathcal{C}^{\infty}-$rings, $J \subseteq B$ a $\mathcal{C}^{\infty}-$radical ideal in $B$ and $f: A \to B$ any $\mathcal{C}^{\infty}-$homomorphism. Then $f^{\dashv}[J]$ is a $\mathcal{C}^{\infty}-$radical ideal in $A$.
\end{theorem}
\begin{proof}See the proof of \textbf{Theorem 20}, p. 82 of \cite{BMTSCA}.
\end{proof}

We register that, in general, holds:

\begin{proposition}\label{usa4africa}Let $A,B$ be $\mathcal{C}^{\infty}-$rings, $f: A \to B$ a $\mathcal{C}^{\infty}-$homomorphism and $J \subseteq B$ any ideal. Then:
$$\sqrt[\infty]{f^{\dashv}[J]} \subseteq f^{\dashv}[\sqrt[\infty]{J}].$$
\end{proposition}
\begin{proof}
See the proof of \textbf{Proposition 33}, p. 83 of \cite{BMTSCA}.
\end{proof}

\begin{remark}With exactly the same method used in the proof of \textbf{Theorem \ref{little}}, one proves that whenever $\mathfrak{p} \subseteq A$ is a prime ideal, $\sqrt[\infty]{\mathfrak{p}}$ is also a prime ideal.
\end{remark}

At this point it is natural to look for a $\mathcal{C}^{\infty}-$analog of the Zariski spectrum of a commutative unital ring. With this motivation, we give the following:

\begin{definition}For a $\mathcal{C}^{\infty}-$ring $A$, we define:
$${\rm Spec}^{\infty}\,(A) = \{ \mathfrak{p} \in {\rm Spec}\,(A) | \mathfrak{p} \, \mbox{is} \, \mathcal{C}^{\infty}-\mbox{radical} \}$$
equipped with the smooth Zariski topology defined by its basic open sets:
$$D^{\infty}(a) = \{ \mathfrak{p} \in {\rm Spec}^{\infty}\,(A) | a \notin \mathfrak{p} \} $$
\end{definition}

A more detailed study of the \index{smooth Zariski spectrum}smooth Zariski spectrum will be given in the next chapter.\\

Now we prove some properties of $\mathcal{C}^{\infty}-$radical ideals of a $\mathcal{C}^{\infty}-$ring $A$

\begin{proposition}\label{egito}The following results hold:
\begin{itemize}
\item[(a)]{Suppose that $(\forall \alpha \in \Lambda)(I_{\alpha} \in \mathfrak{I}^{\infty}_A)$. Then $\bigcap_{\alpha \in \Lambda} I_{\alpha} \in \mathfrak{I}^{\infty}_{A}$, that is, if $(\forall \alpha \in \Lambda)(I_{\alpha} \in \mathfrak{I}^{\infty}_A)$, then:
    $$\sqrt[\infty]{\bigcap_{\alpha \in \Lambda}I_{\alpha}} = \bigcap_{\alpha \in \Lambda}I_{\alpha} = \bigcap_{\alpha \in \Lambda} \sqrt[\infty]{I_{\alpha}}$$}
\item[(b)]{Let $\{ I_{\alpha} | \alpha \in \Sigma \}$ an upward directed family of elements of $\mathfrak{I}^{\infty}_{A}$. Then $\bigcup_{\alpha \in \Sigma} I_{\alpha} \in \mathfrak{I}^{\infty}_{A}$.}
\end{itemize}
\end{proposition}
\begin{proof}
\begin{itemize}
\item[(a)]{Let $u \in \sqrt[\infty]{\bigcap_{\alpha \in \Lambda} I_{\alpha}}$. There exists $b \in \bigcap_{\alpha \in \Lambda} I_{\alpha}$ such that ${\rm Can}_u(b) \in (A\{ u^{-1}\})^{\times}$.\\

    Since $b \in \bigcap_{\alpha \in \Lambda} I_{\alpha}$, then $(\forall \alpha \in \Lambda)(\exists b_{\alpha} \in I_{\alpha})(b_{\alpha = b})({\rm Can}_u(b) \in (A\{ u^{-1}\})^{\times})$, so $u \in \sqrt[\infty]{I_{\alpha}} =  I_{\alpha}$, and therefore $u \in \bigcap_{\alpha \in \Lambda} I_{\alpha}$. This proves that $\sqrt[\infty]{\bigcap_{\alpha \in \Lambda} I_{\alpha}} \subseteq \bigcap_{\alpha \in \Lambda} I_{\alpha}$, so:
    $$\sqrt[\infty]{\bigcap_{\alpha \in \Lambda} I_{\alpha}} = \bigcap_{\alpha \in \Lambda} I_{\alpha},$$
    and
    $$\bigcap_{\alpha \in \Lambda} I_{\alpha} \in \mathfrak{I}^{\infty}_{A}.$$}
\item[(b)]{Given any $u \in \sqrt[\infty]{\bigcup_{\alpha \in \Sigma} I_{\alpha}}$, then there is some $b \in \bigcup_{\alpha \in \Sigma} I_{\alpha}$ such that ${\rm Can}_u(b) \in (A\{ u^{-1}\})^{\times}$. Since $\Sigma$ is upward directed, there is $\alpha_0 \in \Sigma$ such that $b \in I_{\alpha_0}$, so ${\rm Can}_u(b) \in (A\{ u^{-1}\})^{\times}$, therefore $u \in \sqrt[\infty]{I_{\alpha_0}} = I_{\alpha_0} \subseteq \bigcup_{\alpha \in \Sigma} I_{\alpha}$.Thus,

    $$\sqrt[\infty]{\bigcup_{\alpha \in \Sigma} I_{\alpha}} = \bigcup_{\alpha \in \Sigma}I_{\alpha}$$
    }
\end{itemize}
\end{proof}

\begin{remark}Due to the previous proposition, we have that $\mathfrak{I}^{\infty}_{A}$ is a \index{complete Heyting algebra}complete Heyting algebra.
\end{remark}

\begin{lemma}\label{2328}Let $A$ be a $\mathcal{C}^{\infty}-$ring and $S \subset A$ any multiplicative subset. Consider the following posets:
$$({\rm Spec}^{\infty}(A\{ S^{-1}\}), \subseteq)$$
and
$$(\{ \mathfrak{p} \in {\rm Spec}^{\infty}(A) | \mathfrak{p}\cap S = \varnothing \}, \subseteq)$$
The following poset maps:
$$\begin{array}{cccc}
{\rm Can}_S^{*}: & ({\rm Spec}^{\infty}(A\{ S^{-1}\}), \subseteq) & \rightarrow & (\{ \mathfrak{p} \in {\rm Spec}^{\infty}(A) | \mathfrak{p}\cap S = \varnothing \}, \subseteq)\\
                 & Q & \mapsto & {\rm Can}_S^{\dashv}[Q]
\end{array}$$

and

$$\begin{array}{cccc}
{{\rm Can}_S}_{*}: & (\{ \mathfrak{p} \in {\rm Spec}^{\infty}(A) | \mathfrak{p}\cap S = \varnothing \}, \subseteq) & \rightarrow & ({\rm Spec}^{\infty}(A\{ S^{-1}\}), \subseteq)\\
                 & P & \mapsto &\langle {\rm Can}_S[P] \rangle
\end{array}$$
are poset isomorphisms, each one inverse of the other.
\end{lemma}
\begin{proof}See p. 286 of \cite{rings1}.
\end{proof}

\begin{proposition}\label{coisinha}Let $A$ be a $\mathcal{C}^{\infty}-$ring. For every $\mathfrak{p} \in {\rm Spec}^{\infty}(A)$ let $\hat{\mathfrak{p}}$ denote the maximal ideal of $A_{\{ \mathfrak{p}\}} = A\{ {A \setminus \p}^{-1}\}$ and consider:
$${\rm Can}_{\mathfrak{p}}: A \to A_{\{ \mathfrak{p}\}}.$$
We have the following equalities:
$${\rm Can}_{\mathfrak{p}}^{\dashv}[\hat{\mathfrak{p}}] = \mathfrak{p}$$
and
$$({\rm Can}_{\mathfrak{p}}[\mathfrak{p}]) = \hat{\mathfrak{p}}$$
\end{proposition}
\begin{proof}Let us take $S = A \setminus \mathfrak{p}$. Since $\hat{\mathfrak{p}}$ is a maximal ideal, it is the largest element of ${\rm Spec}^{\infty}(A\{ (A \setminus \mathfrak{p})^{-1}\})$. Hence ${\rm Can}_{\mathfrak{p}}^{\dashv}[\hat{\mathfrak{p}}]$ is the largest element of $(\{ \mathfrak{p}' \in {\rm Spec}^{\infty}(A) | \mathfrak{p}'\cap (A\setminus \mathfrak{p}) = \varnothing\})$. Thus, by \textbf{Lemma \ref{2328}}, ${\rm Can}_{\mathfrak{p}}^{\dashv}[\hat{\mathfrak{p}}] = \mathfrak{p}$.
\end{proof}

\begin{proposition}\label{Margarida} If $D$ is a reduced $\mathcal{C}^{\infty}-$domain, then $D\{ a^{-1}\} \cong \{ 0\}$ implies $a=0$.
\end{proposition}
\begin{proof} If $D\{ a^{-1}\} \cong \{ 0\}$ then $\left( \dfrac{D}{(0)}\right)\{ {a + (0)}^{-1} \} \cong \{ 0\}$, so $a \in \sqrt[\infty]{(0)} = (0)$ (for $D$ is reduced), and $a = 0$.
\end{proof}

\begin{proposition}Any free $\mathcal{C}^{\infty}-$ring is a reduced $\mathcal{C}^{\infty}-$ring.
\end{proposition}
\begin{proof}See p. 87 of \cite{BMTSCA}.
\end{proof}

\begin{proposition}\label{Quico}Let $\{ A_i\}_{i \in I}$ be a directed family of $\mathcal{C}^{\infty}-$rings, so we have the diagram:
$$\xymatrixcolsep{5pc}\xymatrix{
 & \varinjlim_{i \in I} A_i & \\
A_i \ar[ur]^{\alpha_i} \ar[rr]_{\alpha_{ij}}& & A_j \ar[ul]_{\alpha_j}
}$$
and let $(\mathfrak{p}_i)_{i \in I}$ be a compatible family of prime $\mathcal{C}^{\infty}-$radical ideals, that is:
$$(\mathfrak{p}_i)_{i \in I} \in \varprojlim_{i \in I} {\rm Spec}^{\infty}(A_i).$$
Under those circumstances,
$$\varinjlim_{i \in I} \mathfrak{p}_i = \bigcup_{i \in I} \alpha_i[\mathfrak{p}_i]$$
is a $\mathcal{C}^{\infty}-$radical prime ideal of $\varinjlim_{i \in I} A_i$.
\end{proposition}
\begin{proof}
  It suffices to prove that
  $$\dfrac{\varinjlim_{i \in I} A_i}{\varinjlim_{i \in I}\mathfrak{p}_i}$$
  is a $\mathcal{C}^{\infty}-$reduced domain.\\

  It is a fact that colimits commute with quotients, so:

  $$\dfrac{\varinjlim_{i \in I} A_i}{\varinjlim_{i \in I}\mathfrak{p}_i} \cong \varinjlim_{i \in I} \dfrac{A_i}{\mathfrak{p}_i}$$

  Now, since every $\mathfrak{p}_i$ is a $\mathcal{C}^{\infty}-$radical prime ideal of $A_i$, we have, for every $i \in I$, that $\dfrac{A_i}{\mathfrak{p}_i}$ is a $\mathcal{C}^{\infty}-$reduced $\mathcal{C}^{\infty}-$domain.\\

  The colimit of $\mathcal{C}^{\infty}-$reduced $\mathcal{C}^{\infty}-$domains is again a $\mathcal{C}^{\infty}-$reduced $\mathcal{C}^{\infty}-$domain, so $\dfrac{\varinjlim_{i \in I} A_i}{\varinjlim_{i \in I}\mathfrak{p}_i}$ is a domain and:
  $$\varinjlim_{i \in I} \mathfrak{p}_i$$
  is a $\mathcal{C}^{\infty}-$radical prime ideal of $\varinjlim_{i \in I} A_i$.
\end{proof}

\section{Separation Theorems for Smooth Commutative Algebra} \label{septheo-sec}

From the notions and results previously established, we are ready to state and prove the The main result of this work:\\

\begin{theorem}[\index{Separation Theorems}\textbf{Separation Theorems}]\label{TS} Let $A$ be a $\mathcal{C}^{\infty}-$ring, $S \subseteq A$ be a subset of $A$ and $I$ be an ideal of $A$. Denote by $\langle S \rangle$ the multiplicative submonoid of $A$ generated by $S$. We have:
\begin{itemize}
  \item[(a)]{If $I$ is a  $\mathcal{C}^{\infty}-$radical ideal, then:

$$I \cap \langle S \rangle = \varnothing \iff {I} \cap S^{\infty-{\rm sat}} = \varnothing$$}
  \item[(b)]{If  $S \subseteq A$ is a $\mathcal{C}^{\infty}$-saturated subset, then:
$$I \cap S = \varnothing \iff \sqrt[\infty]{I} \cap S = \varnothing$$
}
 \item[(c)]{If $\mathfrak{p} \in {\rm Spec}^{\infty}\,(A)$ then $A\setminus \mathfrak{p} = (A \setminus \mathfrak{p})^{\infty-{\rm sat}}$}
 \item[(d)]{If  $S \subseteq A$ is a $\mathcal{C}^{\infty}$-saturated subset, then:
$$I \cap S = \varnothing \iff (\exists \mathfrak{p} \in {\rm Spec}^{\infty}\,(A))((I \subseteq \p)\& (\p \cap S = \varnothing)).$$}
 \item[(e)]{$\sqrt[\infty]{I} = \bigcap \{ \mathfrak{p} \in {\rm Spec}^{\infty}\,(A) | I \subseteq \mathfrak{p} \}$}
\end{itemize}
\end{theorem}
\begin{proof}
Ad (a): Since $\langle S \rangle \subseteq S^{\infty-{\rm sat}}$, it is clear that $(ii) \to (i)$.\\

We are going to show that $(i) \to (ii)$ via contraposition.\\

Suppose there exists some $b \in I \cap S^{\infty-{\rm sat}}$, so $\eta_S(b) \in (B\{ S^{-1}\})^{\times}.$\\

We have:
$$B\{ S^{-1}\} \cong_{\psi} \varinjlim_{S' \stackrel{\subseteq}{{\rm fin}} S} B\{{S'}^{-1} \} \cong_{\varphi} \varinjlim_{S' \stackrel{\subseteq}{{\rm fin}} S} B\left\{{\prod S'}^{-1} \right\},$$
so $(\varphi \circ \psi)(\eta_S(b)) \in \varinjlim_{S' \stackrel{\subseteq}{{\rm fin}} S} B \{ {\prod S'}^{-1}\}$ implies that there is some finite $S'' \subseteq S$ such that $\eta_{S''}(b) \in B\{ {S''}^{-1} \} \cong B\{{\prod S''}^{-1} \}$. Let $a = \prod S''$, so $a \in \langle S \rangle$. We have that $\eta_{a}(b) \in B \{ {S''}^{-1} \}$ implies $\eta_a(b) \in (B\{ a^{-1}\})^{\times}$ and by hypothesis, $b \in I$ so $a \in \sqrt[\infty]{I} = I$. Hence $a \in I \cap \langle S \rangle \neq \varnothing$, and the result is proved.\\

Ad (b): Given $b \in \sqrt[\infty]{I}\cap S$, there must be some $x \in I$ such that $\eta_b(x) \in A\{ b^{-1}\}^{\times}$, so $ x \in \{ b\}^{\infty-{\rm sat}} \subseteq S^{\infty-{\rm sat}} = S$. Thus $x \in I \cap S$ and $I \cap S = \varnothing$. The other way round is immediate since $I \subseteq \sqrt[\infty]{I}$.\\

Ad (c): Since $A\setminus \mathfrak{p} = \langle A \setminus \mathfrak{p}\rangle$, by item (a) we have:

$$\mathfrak{p} \cap (A\setminus \mathfrak{p}) \iff \mathfrak{p}\cap (A \setminus \mathfrak{p})^{\infty-{\rm sat}} = \varnothing$$

so $(A \setminus \mathfrak{p})^{\infty-{\rm sat}} \subseteq A \setminus \mathfrak{p}$. The other inclusion always holds, so $A \setminus \mathfrak{p} = (A \setminus \mathfrak{p})^{\infty-{\rm sat}}$.\\

Ad (d): Consider the set:

$$\Gamma_S := \{ J \in \mathfrak{I}(A) | (I \subseteq J)\& (S \cap J = \varnothing)\}$$

ordered by inclusion. It is straightforward to check that $(\Gamma_S, \subseteq)$ satisfies the hypotheses of \textbf{Zorn's Lemma}. Let $M$ be a maximal member of $\Gamma_S$. We are going to show that $M \in {\rm Spec}^{\infty}\,(A)$.\\

\textbf{Claim:} $M$ is a proper prime ideal of $A$.\\

In fact, $M$ is proper since $1 \in S$ and $S \cap M = \varnothing$.\\

The proof that $M$ is prime is made by contradiction.\\

If $a,a' \notin M$, then by maximality there are $\alpha, \alpha' \in A$ and $m,m' \in M$ such that $m+\alpha \cdot a \in S$ and $m'+\alpha'\cdot a' \in S$. Since $\mathcal{C}^{\infty}-$saturated sets are submonoids, it follows that $(m + \alpha \cdot a)\cdot (m' + \alpha' \cdot a') \in S$. Thus,

$$\underbrace{(m\cdot m' + m \cdot \alpha' \cdot a' + m' \cdot \alpha \cdot a)}_{\in M} + \underbrace{(\alpha \cdot \alpha')\cdot (a \cdot a')}_{\in S} \in S$$

If $a,a' \in M$, we get $M \cap S \neq \varnothing$, a contradiction. Thus, $a,a' \notin M$, and $M$ is a prime ideal.\\

\textbf{Claim:} $M = \sqrt[\infty]{M}$.\\

Since $M \cap S = \varnothing$ and $S$ is $\mathcal{C}^{\infty}-$saturated, by item (b) it follows that $\sqrt[\infty]{M}\cap S = \varnothing$, so $\sqrt[\infty]{M} \in \Gamma_S$. Since $M \subseteq \sqrt[\infty]{M}$, $\sqrt[\infty]{M} \in \Gamma_S$ and $M$ is a maximal element of $\Gamma_S$, it follows that $M = \sqrt[\infty]{M}$.\\

Hence, $M \in {\rm Spec}^{\infty}\,(A)$.\\

Ad (e): Clearly,

$$\sqrt[\infty]{I} \subseteq  \bigcap \{ \mathfrak{p} \in {\rm Spec}^{\infty}\,(A) | I \subseteq \mathfrak{p} \}$$

so we need only to prove the reverse inclusion.\\

Let $a \notin \sqrt[\infty]{I}$, so $S_a = \{ 1, a, a^2, \cdots\} \cap \sqrt[\infty]{I} = \varnothing$. Since $\dfrac{A}{I}\{ (a+I)^{-1}\} \cong \dfrac{A}{I}\{ (a^k+I)^{-1}\}$ for any $k \in \mathbb{N}$ such that $k \geq 1$ and since $\sqrt[\infty]{I}$ is a $\mathcal{C}^{\infty}-$radical ideal, by item (a), we have $(S_a)^{\infty-{\rm sat}}\cap \sqrt[\infty]{I}=\varnothing$, and by item (d), there is some $\mathfrak{p} \in {\rm Spec}^{\infty}\,(A)$ such that $\mathfrak{p} \supseteq \sqrt[\infty]{I} \supseteq I$ such that $a \notin \mathfrak{p}$.
\end{proof}

\begin{proposition} \label{equacaocobertura} Let $A$ be a
$C^\infty$-ring and let
$\{a_i : i \in I\} \subseteq A$. Denote ${\cal I} := \langle \{a_i : i
\in I\}\rangle$. Then the following are equivalent:

(a) ${\rm Spec}^\infty(A) = \bigcup_{i \in I} D^\infty(a_i)$

(b) $A = {\cal I}$

\end{proposition}

\begin{proof}

(b) $\Ra$ (a): Since there exists $\{i_1, \cdots, i_n\} \subseteq I$
such that $ 1_A = \sum_{j = 1}^n \lambda_j. a_{i_j}$
for some $\lambda_1, \cdots, \lambda_n \in A$, then
there is no $\mathfrak{p} \in {\rm Spec}^{\infty}(A)$ such that $\{a_{i_1},
\cdots, a_{i_n} \} \subseteq \mathfrak{p}$, i.e.
${\rm Spec}^\infty(A) \subseteq \bigcup_{i \in I} D^\infty(a_i)$.

(a) $\Ra$ (b): Suppose that $A \neq \mathcal{I}$. Then
$$ A^\times \cap \mathcal{I} = \emptyset.$$

Since $A^\times \subseteq A$ is a $\infty$-saturated subset of $A$
($A^\times = \eta_{1}^{-1} [ (A \{1^{-1}\})^\times]$),
then by the {\bf Separation Theorem
\ref{TS}.(b)},

$$A^\times \cap \sqrt[\infty]{\cal I}= \emptyset.$$

By the {\bf Separation Theorem \ref{TS}.(d)}, there is $\mathfrak{p}
\in {\rm Spec}^\infty(A)$ such that

$$ \sqrt[\infty]{\cal I} \subseteq \mathfrak{p}$$

thus $\mathfrak{p} \in {\rm Spec}^\infty(A) \setminus \bigcup_{i \in I}
D^\infty(a_i)$.

\end{proof}

\section{Order Theory of $\mathcal{C}^{\infty}-$rings}
\label{order-sec}

The class of $\mathcal{C}^{\infty}-$rings carries good notions of order theory for rings. As pointed out by Moerdijk and Reyes in \cite{rings2}, every $\mathcal{C}^{\infty}-$ring $A$ has a canonical pre-order - this and other aspects of the order theory of $\mathcal{C}^{\infty}-$rings  are developed in \cite{berni2020order}. The key point of this section is to introduce the notion of the smooth real spectrum of a $\mathcal{C}^{\infty}-$rings (made in \cite{tese}, denoted by ${\rm Sper}^{\infty}$) and to describe, as a consequence of the separation theorems presented in the previous section, a spectral bijection from the smooth Zariski spectrum to the real smooth spectrum of a $\mathcal{C}^{\infty}-$ring.  We begin by describing the notion of ``order'':\\

\begin{definition}\label{villa}Given a $\mathcal{C}^{\infty}-$ring $(A,\Phi)$, we write:

$$(\forall a \in A)(\forall b \in A)(a \prec b \iff (\exists c \in A^{\times})(b-a = c^2))$$
Note that if $A \neq \{ 0\}$, $\prec$ is an irreflexive relation, i.e., $(\forall a \in A)(\neg (a \prec a))$.
\end{definition}

According to \cite{rings1}, we have the following facts about the order $\prec$ defined above.\\

\begin{fact}Let $A=\dfrac{\mathcal{C}^{\infty}(\mathbb{R}^{E})}{I}$ for some set $E$. Then, given any $f+I,g+I \in A$, with respect to the relation $\prec$, given in \textbf{Definition \ref{villa}}, we have:
$$f \prec g \iff (\exists \varphi \in I)((\forall x \in Z(\varphi))(f(x)<g(x)))$$
so $\prec$ is compatible with the ring structure which underlies $A$, \textit{i.e.}:
\begin{itemize}
  \item[(i)]{$0 \prec f+I,g+I \Rightarrow 0 \prec (f+I) \cdot (g+I)  $;}
  \item[(ii)]{$0 \prec f+I,g+I \Rightarrow 0 \prec f+g + I $}
\end{itemize}
\end{fact}

\begin{fact}Let $A=\dfrac{\mathcal{C}^{\infty}(\mathbb{R}^{E})}{I}$  for some set $E$ be a $\mathcal{C}^{\infty}-$field, so $I = \sqrt[\infty]{I}$. The relation $\prec$, given in \textbf{Definition \ref{villa}}, is such that:
 $$(\forall f+I \in A)(f+I \neq 0+I \rightarrow (f + I \prec 0)\vee (0 \prec f+I))$$
\end{fact}

We have the following:

\begin{proposition}\label{doria}For any  $\mathcal{C}^{\infty}-$ring $A$, we have:
$$1+ \sum A^2 \subseteq A^{\times},$$

where $\sum A^2 = \{ \sum_{i=1}^{n}a_i^2 | n \in \mathbb{N}, a_i \in A\}$. In particular, every $\mathcal{C}^{\infty}$ ring $A$ is such that its underlying commutative unital ring, $\widetilde{U}(A)$ is a semi-real ring.
\end{proposition}
\begin{proof}
First suppose $A = \mathcal{C}^{\infty}(\mathbb{R}^n)$ for some $n \in \mathbb{N}$. Given any $k \in \mathbb{N}$ and any $k-$tuple, $f_1, \cdots, f_k \in \mathcal{C}^{\infty}(\mathbb{R}^n)$, we have $Z\left( 1 + \sum_{i=1}^{k}f_i^2\right) = \varnothing$, so $1 + \sum_{i=1}^{k}f_i^2 \in \mathcal{C}^{\infty}(\mathbb{R}^n)^{\times}$.\\

Suppose, now, that $A$ is a finitely generated $\mathcal{C}^{\infty}-$ring, that is, that there are $n \in \mathbb{N}$ and an ideal $I \subseteq \mathcal{C}^{\infty}(\mathbb{R}^n)$ such that:
$$A = \dfrac{\mathcal{C}^{\infty}(\mathbb{R}^n)}{I}$$

We have, for every $k \in \mathbb{N}$ and for every $k-$tuple $f_1+I, \cdots, f_k +I \in \dfrac{\mathcal{C}^{\infty}(\mathbb{R}^n)}{I}$:

$$(1+I) + \sum_{i=1}^{k}(f_i+I)^2 = \left( 1 + \sum_{i=1}^{k}f_i^2\right) + I$$

Since $1 + \sum_{i=1}^{k}f_i^2 \in \mathcal{C}^{\infty}(\mathbb{R}^n)^{\times}$, it follows that $\left( 1 + \sum_{i=1}^{k}f_i^2\right) + I \in \left( \dfrac{\mathcal{C}^{\infty}(\mathbb{R}^n)}{I}\right)^{\times}$.\\

Finally, for an arbitrary $\mathcal{C}^{\infty}-$ring $A$, we can always write:

$$A = \varinjlim_{i \in I} A_i$$

for the directed family $\{ A_i | i \in I\}$ consisting of its finitely generated $\mathcal{C}^{\infty}-$subrings. We are going to denote the canonical colimit arrow by $\alpha_i: A_i \to \varinjlim_{i \in I}A_i$, for every $i \in I$.\\

For every $k \in \mathbb{N}$ and for every $k-$tuple, $[a_{i_1}], \cdots, [a_{i_k}] \in \varinjlim_{i \in I}A_i$, there is some $\ell \geq i_1, \cdots, i_k$ such that for every $j \in \{ 1, \cdots, k\}$ $[a_{i_j}] = [\alpha_{\ell}(\alpha_{i_j\ell}(a_{i_j}))]$. Since $A_{\ell}$ is finitely generated, it follows that $1 + \sum_{j=1}^{k} \alpha_{i_j\ell}(a_{i_j})^2 \in A_{\ell}^{\times}$, so
$$1 + \sum_{j=1}^{k}[a_{i_j}]^2 = \alpha_{\ell}\left( 1 + \sum_{j=1}^{k} \alpha_{i_j\ell}(a_{i_j})^2\right) \in \left( \varinjlim_{i \in I}A_i\right)^{\times}.$$
\end{proof}

\begin{fact}\label{famosa}Let $n \in \mathbb{N}$. We have:

$$h + I \in \left( \dfrac{\mathcal{C}^{\infty}(\mathbb{R}^n)}{I}\right)^{\times} \iff (\exists \varphi \in I)(\forall x \in Z(\varphi))(h(x) \neq 0)$$
\end{fact}

Recall that a totally ordered field $(F, \leq)$ is \index{real closed field}\textbf{real closed} if it satisfies:
\begin{itemize}
  \item[(a)]{$(\forall x \in F)(0 < x \rightarrow (\exists y \in F)(x=y^2))$;}
  \item[(b)]{every polynomial of odd degree has, at least, one root;}
\end{itemize}

\begin{remark}A \index{$\mathcal{C}^{\infty}-$polynomial}\textbf{$\mathcal{C}^{\infty}-$polynomial in one variable} is an element of $\mathcal{C}^{\infty}(\mathbb{R}^n)\{ t\}$. More generally, a \textbf{$\mathcal{C}^{\infty}-$polynomial in set $S$ of variables} is an element of $A\{ S\}$.
\end{remark}

As pointed out in \textbf{Theorem 2.10} of \cite{rings1}, we have the following:

\begin{fact}Every $\mathcal{C}^{\infty}-$field, $F$, together with its canonical preorder $\prec$ given in \textbf{Definition \ref{villa}}, is such that $\widetilde{U}(F)$ is a real closed field.
\end{fact}

As a consequence of the above fact, given any polynomial $p \in F[x]$ and any $f,g \in F$ such that $f \prec g$, if $p(f)<p(g)$ then there is some $h \in ]f,g[$ such that $p(h)=0$.\\

We have the $\mathcal{C}^{\infty}-$analog of the notion of ``real closedness'':\\

\begin{definition}Let $(F,\prec)$ be a $\mathcal{C}^{\infty}-$field. We say that $(F,\prec)$ is \index{$\mathcal{C}^{\infty}-$real closed}\textbf{$\mathcal{C}^{\infty}-$real closed} if, and only if:

$$(\forall f \in F\{ x\})((f(0)\cdot f(1)<0)\&(1 \in \langle \{ f, f'\}\rangle \subseteq F\{ x\} ) \rightarrow (\exists \alpha \in ]0,1[ \subseteq F)(f(\alpha)=0))$$
\end{definition}

\begin{fact}As proved in \textbf{Theorem 2.10'} of \cite{rings1}, every $\mathcal{C}^{\infty}-$field is $\mathcal{C}^{\infty}-$real closed.
\end{fact}

From the preceding proposition, we conclude that $\overline{f}\prec \overline{g}$ occurs if, and only if, there is a ``witness'' $\varphi \in I$ such that $(\forall x \in Z(\varphi))(f(x)<g(x))$.

Given any $\mathcal{C}^{\infty}-$ring $A$ and any $\mathfrak{p} \in {\rm Spec}^{\infty}\,(A)$, let:

$$k_{\p}:= \left( \dfrac{A}{\p} \right)\left\{ {\dfrac{A}{\p} \setminus \{ 0 + \p\}}^{-1}\right\},$$

that is, $k_{\p}(A)$ is the $\mathcal{C}^{\infty}-$field obtained by taking the quotient $\dfrac{A}{\p}$:

$$q_{\p}: A \to \dfrac{A}{\p}$$

and then taking its $\mathcal{C}^{\infty}-$ring of fractions with respect to $\dfrac{A}{\p}\setminus \{ 0 + \p\}$,

$$\eta_{\dfrac{A}{\p} \setminus \{ 0 + \p\}}: \dfrac{A}{\p} \rightarrow \left( \dfrac{A}{\p}\right)\left\{ {\dfrac{A}{\p} \setminus \{ 0 + \p\}}^{-1}\right\}.$$

The family of $\mathcal{C}^{\infty}-$fields  $\{ k_{\p}(A)| \p \in {\rm Spec}^{\infty}\,(A)\}$ has the following multi-universal property:

``Given any $\mathcal{C}^{\infty}-$homomorphism $f: A \to \mathbb{K}$, where $\mathbb{K}$ is a $\mathcal{C}^{\infty}-$field, there is a unique $\mathcal{C}^{\infty}-$radical prime ideal $\mathfrak{p}$ and a unique  $\mathcal{C}^{\infty}-$homomorphism $\widetilde{f}: k_{\p}(A) \to \mathbb{K}$ such that the following diagram commutes:
$$\xymatrixcolsep{5pc}\xymatrix{
A \ar[r]^{\alpha_{\p}} \ar[dr]^{f} & k_{\p}(A) \ar@{.>}[d]^{\widetilde{f}}\\
  & K},$$

where $\alpha_{\p} = \eta_{\dfrac{A}{\p} \setminus \{ 0 + \p\}} \circ q_{\p}: A \rightarrow k_{\p}(A)$.''

Thus, given $f: A \to K$, take $\p = \ker(f)$, so $\overline{f}$ is injective, $\dfrac{A}{\p}$ is a $\mathcal{C}^{\infty}-$reduced $\mathcal{C}^{\infty}-$ring. By the universal property of the smooth fraction field $k_{\p}(A)$, there is a unique arrow $\widetilde{f}: k_{\p}(A) \to {K}$ such that the following diagram commutes:

$$\xymatrixcolsep{5pc}\xymatrix{
\dfrac{A}{\p} \ar[r]^{\alpha_{\p}} \ar[dr]_{f} & k_{\p}(A) \ar@{-->}[d]^{\widetilde{f}}\\
 & {K}
}$$

\begin{definition} \label{rel}
Let $\mathcal{F}$ be the (proper) class of all the $\mathcal{C}^{\infty}-$homomorphisms of  $A$ to some $\mathcal{C}^{\infty}-$field. We define the following relation $\mathcal{R}$: given $h_1: A \to F_1$ and $h_2: A \to F_2$, we say that $h_1$ is related with $h_2$ if, and only if, there is some $\mathcal{C}^{\infty}-$field $\widetilde{F}$ and some $\mathcal{C}^{\infty}-$fields homomorphisms $\mathcal{C}^{\infty}$ $f_1: F_1 \to \widetilde{F}$ and $f_2: F_2 \to \widetilde{F}$ such that the following diagram commutes:

$$\xymatrixcolsep{5pc}\xymatrix{
  & F_1 \ar[dr]^{f_1} & \\
A \ar[ur]^{h_1} \ar[dr]_{h_2} & & \widetilde{F} \\
  & F_2 \ar[ur]_{h_2}
}$$

The relation $\mathcal{R}$ defined above is symmetric and reflexive.
\end{definition}

Let $h_1: A \to F_1$ and $h_2: A \to F_2$ be two $\mathcal{C}^{\infty}-$homomorphisms of  $A$ to the $\mathcal{C}^{\infty}-$fields $F_1,F_2$ such that $(h_1,h_2) \in \mathcal{R}$, and let $f_1: F_1 \to \widetilde{F}$ and $f_2: F_2 \to \widetilde{F}$ be two $\mathcal{C}^{\infty}-$homomorphisms to the $\mathcal{C}^{\infty}-$field $\widetilde{F}$ such that $f_1 \circ h_1= f_2\circ h_2$, so:
$$(f_1 \circ h_1)^{\dashv}[\{ 0\}] = (f_2 \circ h_2)^{\dashv}[\{ 0\}]$$
Then
$$\ker(h_1) = h_1^{\dashv}[\{ 0\}] = h_1^{\dashv}[f_1^{\dashv}[\{ 0\}]] =  h_2^{\dashv}[f_2^{\dashv}[\{ 0\}]] = h_2^{\dashv}[\{ 0\}] = \ker(h_2).$$

The above considerations prove the following:

\begin{proposition}If $h_1: A \to F_1$ and $h_2: A \to F_2$ be two $\mathcal{C}^{\infty}-$homomorphisms from the $\mathcal{C}^{\infty}-$ring  $A$ to the $\mathcal{C}^{\infty}-$fields $F_1,F_2$ such that $(h_1,h_2) \in \mathcal{R}$, then $\ker(h_1)=\ker(h_2)$.
\end{proposition}

The above proposition  has the following immediate consequence:

\begin{corollary}Keeping the same notations of the above result, let $\mathcal{R}^t$ be the transitive closure of $\mathcal{R}$. Then $(h_1,h_2) \in \mathcal{R}^t$ implies $\ker(h_1) = \ker(h_2)$. Thus, $\mathcal{R}^t$ is an equivalence relation on $\mathcal{F}$. We are going to denote the quotient set $\dfrac{\mathcal{F}}{\mathcal{R}^t}$ by $\widetilde{\mathcal{F}}$.
\end{corollary}

Let $F_1, F_2, \widetilde{F}$ be $\mathcal{C}^{\infty}-$fields  and $f_1: F_1 \to \widetilde{F}$ and $f_2: F_2 \to \widetilde{F}$ be two following $\mathcal{C}^{\infty}-$fields homomorphisms. Since $\mathcal{C}^{\infty}-$fields homomorphisms must be injective maps, we have the following:\\

\begin{proposition}The following relation
$$\beta = \{ ([h: A \to F], \ker (h)) | F \, \mbox{is a}\, \mathcal{C}^{\infty}-\mbox{field} \} \subseteq \widetilde{\mathcal{F}} \times {\rm Spec}^{\infty}(A)$$

is a functional relation whose domain is $\widetilde{\mathcal{F}}$.
\end{proposition}
\begin{proof}
Suppose $[h: A \to F_1] = [g: A \to F_2]$, so there are maps $f_1: F_1 \to \widetilde{F}$ and $f_2: F_2 \to \widetilde{F}$ for some $\mathcal{C}^{\infty}-$field $\widetilde{F}$, such that the following diagram commutes:

$$\xymatrix{
    & F \ar[rd]^{f_1} &    \\
A \ar[ur]^{h} \ar[dr]^{g} &  & \widetilde{F}\\
        & F_2 \ar[ur]^{f_2} & \\
}.$$

Now, if $([h: A \to F_1], \ker (h)), ([g: A \to F_2], \ker (g)) \in \beta$ are such that $[h: A \to F_1] = [g: A \to F_2]$, then:

$$f_1 \circ h = f_2 \circ g$$
so
$$\ker(h) =  h^{\dashv}[\{ 0 \}] = h^{\dashv}[f_1^{\dashv}[\{0\}]] = \ker(f_1 \circ h) = \ker (f_2 \circ g) = g^{\dashv}[f_2^{\dashv}[\{ 0\}]] = g^{\dashv}[\{ 0 \}] = \ker(g)
$$
\end{proof}

\begin{definition}Let $A$ be a $\mathcal{C}^{\infty}-$ring. A \index{$\mathcal{C}^{\infty}-$ordering}$\mathcal{C}^{\infty}-$\textbf{ordering} in $A$ is a subset $P \subseteq A$ such that:
\begin{itemize}
  \item[(O1)]{$P + P \subseteq P$;}
  \item[(O2)]{$P \cdot P \subseteq P$;}
  \item[(O3)]{$P \cup (-P) = A$}
  \item[(O4)]{$P \cap (-P) = \mathfrak{p} \in {\rm Spec}^{\infty}\,(A)$}
\end{itemize}
\end{definition}

\begin{definition}Let $A$ be a $\mathcal{C}^{\infty}-$ring. Given a $\mathcal{C}^{\infty}-$ordering $P$ in $A$, the \index{$\mathcal{C}^{\infty}-$support}$\mathcal{C}^{\infty}-$\textbf{support} of $A$ is given by:

$${\rm supp}^{\infty}(P):= P \cap (-P)$$
\end{definition}

\begin{definition}\label{Santo}Let $A$ be a $\mathcal{C}^{\infty}-$ring. The \index{$\mathcal{C}^{\infty}-$real spectrum of  $A$}\textbf{$\mathcal{C}^{\infty}-$real spectrum of  $A$} is given by:

$${\rm Sper}^{\infty}\,(A)=\{ P \subset A | P \, \mbox{is an ordering of the elements of}\,\, A\}$$
together with the (spectral) topology generated by the sets:

$$H^{\infty}(a) = \{ P \in {\rm Sper}^{\infty}\,(A) | a \in P \setminus {\rm supp}^{\infty}\,(P)\}$$

for every $a \in A$. The topology generated by these sets will be called ``\index{smooth Harrison topology}smooth Harrison topology'', and will be denoted by ${\rm Har}^{\infty}$.
\end{definition}

\begin{remark}\label{avio}Given a $\mathcal{C}^{\infty}-$ring $A$, we have a  function given by:

$$\begin{array}{cccc}
    {\rm supp}^{\infty}: & ({\rm Sper}^{\infty}(A), {\rm Har}^{\infty}) & \rightarrow & ({\rm Spec}^{\infty}\,(A), {\rm Zar}^{\infty}) \\
     & P & \mapsto & P \cap (-P)
  \end{array}$$

which is spectral, and thus continuous, since given any $a \in A$, ${{\rm supp}^{\infty}}^{\dashv}[D^{\infty}(a)] = H^{\infty}(a)\cup H^{\infty}(-a)$.
\end{remark}

Contrary to what happens to a general commutative ring $R$, for which the mapping:

$$\begin{array}{cccc}
    {\rm supp} : & {\rm Sper}\,(R) & \rightarrow & {\rm Spec}\,(R) \\
     & P & \mapsto & P \cap (-P)
  \end{array}$$

is seldom surjective or injective, within the category of $\mathcal{C}^{\infty}-$rings ${\rm supp}^{\infty}$ is, as matter of fact, a bijection.\\

In order to prove this fact, we are going to need some preliminary results, given below.\\

\begin{lemma}\label{90}Let $A$ be a $\mathcal{C}^{\infty}-$ring and $\mathfrak{p}$ any  $\mathcal{C}^{\infty}-$radical prime ideal, and let $\hat{\mathfrak{p}}$ be the maximal ideal of $A_{\{ \mathfrak{p}\}}$. Then:
$${\rm Can}_{\mathfrak{p}}^{\dashv}[\hat{\mathfrak{p}}] = \mathfrak{p}.$$
\end{lemma}
\begin{proof}
Let $a \in \mathfrak{p}$, then ${\rm Can}_{\mathfrak{p}}(a) \in \hat{\mathfrak{p}}$, and $\mathfrak{p} \subseteq {\rm Can}_{\mathfrak{p}}^{\dashv}[\hat{\mathfrak{p}}]$. Now, if $a \in A \setminus \mathfrak{p}$ then ${\rm Can}_{\mathfrak{p}}(a) \in \mathcal{U}(A\{ (A\setminus \mathfrak{p})^{-1} \})$. Since $A_{\{ \mathfrak{p} \}}$ is a local ring, $A_{\{\mathfrak{p} \}} = \hat{\mathfrak{p}}\stackrel{\cdot}{\cup} \mathcal{U}(A\{ (A \setminus \mathfrak{p})^{-1} \})$, so ${\rm Can}_{\mathfrak{p}}(a) \in A\{ (A\setminus \mathfrak{p})^{-1} \}\setminus \hat{\mathfrak{p}}$, and therefore ${\rm Can}_{\mathfrak{p}}^{\dashv}[\hat{\mathfrak{p}}] \subseteq \mathfrak{p}$.
\end{proof}

\begin{theorem}\label{1}Let $A$ be a $\mathcal{C}^{\infty}-$ring and $\mathcal{R}^t$ be the relation defined above. The function:
$$\begin{array}{cccc}
\alpha': & {\rm Sper}^{\infty}\,(A) & \rightarrow & \frac{\mathcal{F}}{\mathcal{R}^t}\\
        & P & \mapsto & [\eta_{P \cap (-P)}]
\end{array}$$
is the inverse function of:
$$\begin{array}{cccc}
\beta': & \frac{\mathcal{F}}{\mathcal{R}^t} & \to & {\rm Sper}^{\infty}\,(A)\\
       & [h: A \to K] & \mapsto & h^{\dashv}[K^2]
\end{array}$$
\end{theorem}
\begin{proof}
Note that:
$$(\alpha' \circ \beta')([h: A \to F]) = \alpha'(h^{\dashv}[F^2]) = [\eta_{{\rm supp}\,(h^{\dashv}[F^2])}],$$
where ${\rm supp}^{\infty}(h^{\dashv}[F^2]) = h^{\dashv}[F^2] \cap(- h^{\dashv}[F]^2) = h^{\dashv}[\{0\}] = \ker (h)$.\\

Thus we have:
$$ (\alpha' \circ \beta')([h: A \to F])= [\eta_{\ker(h)}: A \to k_{\ker(h)}(A)].$$

We claim that $\beta'$ is the left inverse function for $\alpha'$, that is:

$$(\forall P \in \, {\rm Sper}^{\infty}(A))((\beta' \circ \alpha')(P) = P).$$

Thus, it will follow that $\alpha'$ is injective and $\beta'$ is surjective.\\

We have $\beta'(\alpha'(P)) = \eta_{{\rm supp}^{\infty}(P)}^{\dashv}[k_{\p}(A)^{2}]$, so we need to show that:

$$\eta_{{\rm supp}(P)}^{\dashv}[k_{\p}(A)^2] = P$$

Let $\mathfrak{p} = {\rm supp}^{\infty}(P)$.\\

\textit{Ab absurdo}, suppose

$$\eta_{\mathfrak{p}}^{\dashv}[k_{\p}(A)^2] \nsubseteq P$$

There must exist some $x \in A$ such that $x \in \eta_{\mathfrak{p}}^{\dashv}[k_{\p}(A)^2]$ and $x \notin P$. We have:

$$\eta_{\mathfrak{p}}(x) \in \left( \frac{A\{ A \setminus {\mathfrak{p}}^{-1}\}}{\widehat{\mathfrak{p}}}\right)^2 \,\,\, {\rm and} \,\,\, x \notin P.$$

Now, since by \textbf{Theorem 24}, p. 97 of \cite{BMTSCA} (denoting $\mathfrak{m}_{\p}$ by $\widehat{\p}$ instead)  the following diagram commutes:

$$\xymatrixcolsep{5pc}\xymatrix{
 & \dfrac{A\{ {A\setminus \p}^{-1}\}}{\widehat{\p}} \ar@<1ex>[dd]^{\varphi_{\p}}\\
A \ar[ur]^{\eta_{\p}'} \ar[dr]_{\eta_{\p}} & \\
   & \left( \dfrac{A}{\p}\right)\left\{ {\dfrac{A}{\p} \setminus \{ 0 + \p\}}^{-1}\right\} \ar@<1ex>[uu]^{\psi_{\p}}}$$

where $\eta_{\p} = \eta_{\dfrac{A}{\p} \setminus \{ 0 + \p\}} \circ q_{\p}$ and $\eta_{\p}' = q_{\widehat{\p}} \circ \eta_{A\setminus \p}$ and $\varphi_{\p}$ and $\psi_{\p}$ are the isomorphisms described in that theorem. Thus, we have:

$$\eta_{\mathfrak{p}}(x) \in \left( k_{\p}(A)\right)^2 \Rightarrow \psi_{\p}(\eta_{\p}(x)) = \eta_{\p}'(x) \in \left(\dfrac{A\{ {A\setminus \p}^{-1}\}}{\widehat{\p}}\right)^2$$

and

$$\eta_{\mathfrak{p}}(x) \in \left( k_{\p}(A)\right)^2 \Rightarrow (\exists (g+ \widehat{\mathfrak{p}}) \in \frac{A\{ {A \setminus \mathfrak{p}}^{-1}\}}{\widehat{\mathfrak{p}}})(\eta_{\mathfrak{p}}'(x) = g^2 + \widehat{\mathfrak{p}})$$

Since $q_{\widehat{\mathfrak{p}}}$ is surjective, given this $g + \widehat{\mathfrak{p}} \in  \left( \frac{A\{ A \setminus {\mathfrak{p}}^{-1}\}}{\widehat{\mathfrak{p}}}\right)$, there is some $\theta \in A\{ {A \setminus \mathfrak{p}}^{-1}\}$ such that $q_{\widehat{\mathfrak{p}}}(\theta) = g + \widehat{\mathfrak{p}}$.\\

By \textbf{Theorem 1.4}, item (i) of \cite{rings1}, given this $\theta \in A\{ {A \setminus \mathfrak{p}}^{-1}\}$, there are $a \in A$ and $b \in {A\setminus \mathfrak{p}}^{\infty-{\rm sat}}$, that is,
$${\rm Can}_{\mathfrak{p}}(b) \in (A\{ {A \setminus \mathfrak{p}}^{-1}\})^{\times}$$
such that:
$$\theta = \dfrac{{\rm Can}_{\mathfrak{p}}(a)}{{\rm Can}_{\mathfrak{p}}(b)}$$
or equivalently, since $(A \setminus \mathfrak{p})^{\infty-{\rm sat}} = A \setminus \mathfrak{p}$:
\begin{equation}\label{absu} b \notin \mathfrak{p}
\end{equation}

Hence,

$$\eta_{\mathfrak{p}}'(x) = g^2 + \widehat{\mathfrak{p}} = q_{\widehat{\mathfrak{p}}}\left( \frac{{\rm Can}_{\mathfrak{p}}(a)}{{\rm Can}_{\mathfrak{p}}(b)} \right)^2 = \left( \frac{{\rm Can}_{\mathfrak{p}}(a)}{{\rm Can}_{\mathfrak{p}}(b)} \right)^2 + \widehat{\mathfrak{p}}$$

$$\eta_{\mathfrak{p}}'(x)\cdot ({\rm Can}_{\mathfrak{p}}^2(b) + \widehat{\mathfrak{p}}) = {\rm Can}_{\mathfrak{p}}^2(a) + \widehat{\mathfrak{p}}$$

$${\rm Can}_{\mathfrak{p}}(x\cdot b^2 - a^2) \in \widehat{\mathfrak{p}}.$$

$$(x \cdot b^2 - a^2) \in {\rm Can}_{\mathfrak{p}}^{\dashv}[\widehat{\mathfrak{p}}].$$

By \textbf{Lemma \ref{90}}, $\mathfrak{p} = {\rm Can}_{\mathfrak{p}}^{\dashv}[\widehat{\mathfrak{p}}]$, so
$$x \cdot b^2 - a^2 \in \mathfrak{p} \subseteq P$$

Let $y = x \cdot b^2 + (-a^2) \in \mathfrak{p} \subseteq P$. Note that since $x \notin P$, $x \in (-P)\setminus \mathfrak{p}$ and
\begin{equation}\label{londonbeat}
x \cdot b^2 \in (-P)
\end{equation}

Since $y \in P$,
$$x \cdot b^2 = \underbrace{y}_{\in P} + \overbrace{a^2}^{\in P} \in P,$$

\begin{equation}\label{madonna}
  x \cdot b^2 \in P
\end{equation}

By \eqref{londonbeat} and \eqref{madonna}, it follows that $x \cdot b^2 \in P \cap (-P) = \mathfrak{p}$. Since $\mathfrak{p}$ is prime, either $x \in \mathfrak{p}$ or $b^2 \in \mathfrak{p}$. However, since $x \notin P$, \textit{a fortiori}, $x \notin \mathfrak{p}$, so we must have $b^2 \in \mathfrak{p}$. Once again, since $\mathfrak{p}$ is prime, it follows that $b \in \mathfrak{p}$ which contradicts \eqref{absu}. Hence,

$${\eta_{{\rm supp}^{\infty}(P)}'}^{\dashv}\left[\left( \frac{A\{ A \setminus {\rm supp}^{\infty}(P)^{-1}\}}{\widehat{\mathfrak{p}}}\right)^2\right] \subseteq P$$.\\

Now we claim that:\\

$$P \subseteq {\eta_{{\rm supp}(P)}'}^{\dashv}\left[\left( \frac{A\{ A \setminus {\rm supp}^{\infty}(P)^{-1}\}}{\widehat{\mathfrak{p}}}\right)^2\right] $$

Conversely, suppose, \textit{ab absurdo} that

\begin{equation}\label{Adma}
P \nsubseteq {\eta_{{\rm supp}(P)}'}^{\dashv}\left[\left( \frac{A\{ A \setminus {\rm supp}^{\infty}(P)^{-1}\}}{\widehat{\mathfrak{p}}}\right)^2\right]
\end{equation}

so there must exist some $x \in P$ such that

\begin{equation}\label{sarue}
 (\forall (g + \widehat{\mathfrak{p}}) \in \dfrac{A\{{A \setminus \mathfrak{p}}^{-1}\}}{\widehat{\mathfrak{p}}})(\eta_{\mathfrak{p}}'(x) \neq g^2 + \widehat{\mathfrak{p}})
\end{equation}

Equivalently, there must exist some $(h+\widehat{\mathfrak{p}}) \in \frac{A\{ {A \setminus \mathfrak{p}}^{-1}\}}{\widetilde{\mathfrak{p}}}$ such that $\eta_{\mathfrak{p}}(x) = - h^2 + \widehat{\mathfrak{p}}$.\\

Thus, since $q_{\widehat{\mathfrak{p}}}$ is surjective, given such an $h + \widehat{\mathfrak{p}} \in \dfrac{A\{{A \setminus \mathfrak{p}}^{-1}\}}{\widehat{\mathfrak{p}}}$ there must be some $\zeta \in A\{ {A \setminus \mathfrak{p}}^{-1}\}$ such that $q_{\widetilde{\mathfrak{p}}}(\zeta) = h + \widehat{\mathfrak{p}}$. By item (i) of \textbf{Theorem \ref{340}}, there are $c \in A$ and $d \in A$ with:
$${\rm Can}_{\mathfrak{p}}(d) \in A\{ {A \setminus \mathfrak{p}}^{-1}\}^{\times},$$
equivalently
$$d \in (A \setminus \mathfrak{p})^{\infty-{\rm sat}},$$
and since $(A \setminus \mathfrak{p})^{\infty-{\rm sat}} = A \setminus \mathfrak{p}$,
\begin{equation}\label{criat}
d \notin \mathfrak{p}
\end{equation}

such that:
$$\zeta = \dfrac{{\rm Can}_{\mathfrak{p}}(c)}{{\rm Can}_{\mathfrak{p}}(d)}.$$

Hence,

$$\eta_{\mathfrak{p}}(x) = - \dfrac{{\rm Can}_{\mathfrak{p}}^2(c)}{{\rm Can}_{\mathfrak{p}}^2(d)} + \widehat{\mathfrak{p}}$$

$${\rm Can}_{\mathfrak{p}}^2(c) + {\rm Can}_{\mathfrak{p}}(x) \cdot {\rm Can}_{\mathfrak{p}}^2(d) \in \widehat{\mathfrak{p}}$$

$${\rm Can}_{\mathfrak{p}}(c^2 + x \cdot d^2) \in \widehat{\mathfrak{p}}$$

so

$$\underbrace{c^2}_{\in P} + \overbrace{x \cdot d^2}^{\in P} = z = {\rm Can}_{\mathfrak{p}}^{\dashv}[\widehat{\mathfrak{p}}] = \mathfrak{p} \subseteq (-P)$$

$$x\cdot d^2 = z - c^2 \in (-P)$$

Since $x \in P$, we also have $x \cdot d^2 \in P$, hence $x \cdot d^2 \in \mathfrak{p}$ .  Since $\mathfrak{p}$ is prime, either $x \in \mathfrak{p}$ or $d^2 \in \mathfrak{p}$. Now, if $x \in \mathfrak{p}$ then ${\rm Can}_{\mathfrak{p}}(x) = 0^2+\mathfrak{p}$, which contradicts our hypothesis \eqref{sarue}. On the other hand, if $d^2 \in \mathfrak{p}$, then $d \in \mathfrak{p}$, and this contradicts \eqref{criat}. Hence $x \notin \mathfrak{p}$ and $d^2 \notin \mathfrak{p}$, so $x \cdot d^2 \notin \mathfrak{p}$. Thus we achieved an absurdity: $(x \cdot d^2 \in \mathfrak{p}) \& (x \cdot d^2 \notin \mathfrak{p})$. It follows that our premise \eqref{Adma} must be false, so:

$$P \subseteq {\eta_{{\rm supp}(P)}'}^{\dashv}\left[\left( \frac{A\{ A \setminus {\rm supp}^{\infty}(P)^{-1}\}}{\widehat{\mathfrak{p}}}\right)^2\right].$$

Hence $P = {\eta_{{\rm supp}(P)}'}^{\dashv}\left[\left( \frac{A\{ A \setminus {\rm supp}^{\infty}(P)^{-1}\}}{\widehat{\mathfrak{p}}}\right)^2\right]$.\\

%It is not the case that $x \in \p$, for if it was, $x \in \p = {\rm Can}_{\p}^{\dashv}[\widehat{\p}]$ so  ${\rm Can}_{\p}(x) \in \widehat{\p}$ and $\eta_{\p}(x) = q_{\widehat{\p}}({\rm Can}_{\p}(x)) =\widehat{\p} = 0^2 + \widehat{\p}$, which contradicts \eqref{sarue}. Thus, $d^2 \in \p$, and since $\p$ is a prime ideal, it follows that $d \in \p$. Under these circumstances, we would have:

%$${\rm Can}_{\p}(d) \in \widehat{\p}$$

%and since $A\{ {A \setminus \p}^{-1}\} = \widehat{\p} \stackrel{\cdot}{\cup} (A \{ {A \setminus \p}^{-1}\})^{\times}$,

%\begin{equation}\label{pitanga}
%{\rm Can}_{\p}(d) \notin (A\{ {A \setminus \p}^{-1}\})^{\times}
%\end{equation},

%which contradicts \eqref{criat}. Hence:

%$$P \subseteq \eta_{{\rm supp}(P)}^{\dashv}[\left( \frac{A\{ A \setminus \supp^{\infty}(P)^{-1}\}}{\widehat{\p}}\right)^2]$$

%It follows that:

%$$ \eta_{{\rm supp}(P)}^{\dashv}[\left( \frac{A\{ A \setminus \supp^{\infty}(P)^{-1}\}}{\widehat{\p}}\right)^2] = P$$

%$$(\forall P \in \, {\rm Sper}^{\infty}(A))((\beta' \circ \alpha')(P) = \eta_{{\rm supp}(P)}^{\dashv}[\left( \frac{A\{ A \setminus \supp^{\infty}(P)^{-1}\}}{\widehat{\p}}\right)^2] = P).$$

Now we need only to show that $\alpha' \circ \beta' = {\rm id}_{\widetilde{\mathcal{F}}}$. \\

Let $[h: A \to F] \in \widetilde{\mathcal{F}}$. We have:

$$(\alpha' \circ \beta')([h: A \to F]) = \alpha'(h^{\dashv}[F^2]) = [\eta_{{\rm supp}^{\infty}(h^{\dashv}[F^2])}: A \to k_{{\rm supp}^{\infty}(h^{\dashv}[F^2])}(A)].$$

It suffices to show that $[h]=[\eta_{{\rm supp}^{\infty}(h^{\dashv}[F^2])}]$. Note that ${\rm supp}^{\infty}(h^{\dashv}[F^2]) = h^{\dashv}[F^2 \cap (-F^2)] = \ker(h)$.\\

By the universal property of the $\mathcal{C}^{\infty}-$field of fractions of $\left( \frac{A}{\ker(h)}\right)$, $k_{\ker(h)}(A)$, since $h \left[ A^{\times}\right] \subseteq F^{\times}$ (for $\mathcal{C}^{\infty}-$homomorphisms preserve invertible elements), there is a unique $\mathcal{C}^{\infty}-$homomorphism $\widetilde{h}: k_{\ker(h)}(A) \to F$ such that the following diagram commutes:

$$\xymatrixcolsep{5pc}\xymatrix{
A \ar[r]^{\eta_{\p}} \ar[rd]^{h} & k_{\ker(h)}(A) \ar[d]^{\widetilde{h}}\\
  & F
}$$

We have, then, the following commutative diagram:

$$\xymatrixcolsep{5pc}\xymatrix{
  & F \ar[dr]^{{\rm id}_F} &  \\
  A \ar[ur]^{h} \ar[dr]_{\eta_{\mathfrak{p}}} & & F\\
    & k_{\ker(h)} \ar@{.>}[ur]_{\widetilde{h}}
}$$

so $[h: A \to F] = [\eta_{\ker(h)}: A \to k_{\ker(h)}]$ and

$$(\alpha' \circ \beta')([h: A \to F]) = [h: A \to F].$$

Hence it follows that $\alpha'$ and $\beta'$ are inverse bijections of each other.
\end{proof}

\begin{theorem}\label{2}The map:
$$\begin{array}{cccc}
\alpha :& {\rm Spec}^{\infty}\,(A) & \rightarrow & \widetilde{\mathcal{F}}\\
        & \mathfrak{p}  & \mapsto & [\eta_{\mathfrak{p}}]= [q \circ {\rm Can}_{\mathfrak{p}}]
\end{array}$$
is a bijection whose inverse is given by:
$$\begin{array}{cccc}
\beta: & \widetilde{\mathcal{F}} & \rightarrow & {\rm Spec}^{\infty}\, (A)\\
       & [h: A  \to F] & \mapsto & \ker(h)
\end{array}$$
\end{theorem}
\begin{proof}
First we are going to show that $\alpha \circ \beta = {\rm id}_{\widetilde{\mathcal{F}}}$, so $\beta$ is an injective map and $\alpha$ is its left inverse, hence it is surjective.\\

Let $[h: A \to F] \in \widetilde{F}$. We have:

$$(\alpha \circ \beta)([h: A \to F]) = \alpha(\ker(h)) = \left[ \eta_{\ker(h)}': A \to \left( \frac{A\{ {A \setminus \ker(h)}^{-1} \}}{\widehat{\ker(h)}}\right)\right]$$

It suffices to show that $[h: A \to F] = \left[ \eta_{\ker(h)}': A \to \left( \frac{A\{ {A \setminus \ker(h)}^{-1} \}}{\widehat{\ker(h)}}\right)\right]$.\\

Since $\left( \frac{A\{ {A \setminus \ker(h)}^{-1} \}}{\widehat{\ker(h)}}\right)$ is (up to $\mathcal{C}^{\infty}-$isomorphism) the $\mathcal{C}^{\infty}-$field of fractions of $\left( \frac{A}{\ker(h)}\right)$, $k_{\ker(h)}(A)$, there is a unique $\mathcal{C}^{\infty}-$homomorphism $\widetilde{h}: \left( \frac{A\{ {A \setminus \ker(h)}^{-1} \}}{\widehat{\ker(h)}}\right) \to F$ such that the following diagram commutes:

$$\xymatrixcolsep{5pc}\xymatrix{
  & \left( \frac{A\{ {A \setminus \ker(h)}^{-1} \}}{\widehat{\ker(h)}}\right) \ar[dr]^{\widetilde{h}} &   \\
A \ar[ur]^{\eta_{{\rm supp}^{\infty}(h^{\dashv}[F^2])'}} \ar[dr]_{h} & & F\\
  & F \ar@{.>}[ur]_{{\rm id}_F} &}$$
and the equality holds, \textit{i.e.},

$$[h: A \to F] = \left[ \eta_{\ker(h)}: A \to \left( \frac{A\{ {A \setminus \ker(h)}^{-1} \}}{\widehat{\ker(h)}}\right)\right]$$

It follows that $\alpha \circ \beta = {\rm id}_{\widetilde{F}}$, $\alpha$ is a surjective map and $\beta$ is an injective map.\\

On the other hand, given $\mathfrak{p} \in {\rm Spec}^{\infty}(A)$, we have:

$$(\beta \circ \alpha)(P) = \beta([\eta_{\mathfrak{p}}: A \to k_{\mathfrak{p}}]) = \ker(\eta_{\mathfrak{p}}) = {\rm Can}_{\mathfrak{p}}^{\dashv}[\widehat{\mathfrak{p}}] = \mathfrak{p}$$

so:

$$(\beta \circ \alpha) = {\rm id}_{{\rm Spec}^{\infty}\,(A)}$$
\end{proof}

as a Corollary of the theorems \textbf{\ref{1}} and \textbf{\ref{2}}, we have:

\begin{lemma}Let $A$ be a $\mathcal{C}^{\infty}-$ring, and define:

$$\begin{array}{cccc}
{\rm supp}^{\infty} : & {\rm Sper}^{\infty}\,(A) & \rightarrow & {\rm Spec}^{\infty}\,(A)\\
      & P & \mapsto & P \cap(-P)
\end{array}$$

The following diagram commutes:

$$\xymatrixcolsep{5pc}\xymatrix{
   & \,\, {\rm Sper}^{\infty}\,(A) \ar@<1ex>[dl]^{\alpha'} \ar[dd]_{{\rm supp}^{\infty}}\\
   \widetilde{\mathcal{F}} \ar@<1ex>[ur]^{\beta'} \ar@<1ex>[dr]^{\beta} &  \\
      & {\rm Spec}^{\infty}\,(A) \ar@<1ex>[ul]^{\alpha}
}$$

that is to say that:
$$\alpha \circ {\rm supp}^{\infty} = \alpha'$$
and
$${\rm supp}^{\infty}\circ \beta' = \beta$$
\end{lemma}
\begin{proof}
Note that if we prove that $\alpha \circ {\rm supp}^{\infty} = \alpha'$, then composing both sides with $\beta'$ yields:
$$(\alpha \circ {\rm supp}^{\infty})\circ \beta' = \alpha' \circ \beta' = {\rm id}_{\widetilde{\mathcal{F}}}$$
so
$$\alpha \circ ({\rm supp}^{\infty} \circ \beta') = {\rm id}_{\widetilde{\mathcal{F}}}$$
and by the uniqueness of the inverse of $\alpha$, it follows that:
$${\rm supp}^{\infty} \circ \beta' = \beta.$$

Now we are going to prove that $\alpha \circ {\rm supp}^{\infty} = \alpha'$.\\

Given $P \in {\rm Sper}^{\infty}\,(A)$ we have:

$$(\alpha \circ {\rm supp}^{\infty})(P) = \alpha({\rm supp}^{\infty}(P)) = [\eta_{{\rm supp}^{\infty}\,(P)}: A \to k_{{\rm supp}^{\infty}\,(P)}(A)] =: \alpha'(P),$$

so the result holds.
\end{proof}

As an important result of the theory of $\mathcal{C}^{\infty}-$rings which distinguishes it from the theory of the rings, we have the following:

\begin{theorem}\label{exv}Let $A$ be a $\mathcal{C}^{\infty}-$ring. The following map:
$$\begin{array}{cccc}
{\rm supp}^{\infty} : & {\rm Sper}^{\infty}\,(A) & \rightarrow & {\rm Spec}^{\infty}\,(A)\\
      & P & \mapsto & P \cap(-P)
\end{array}$$
is a spectral \underline{bijection}.
\end{theorem}
\begin{proof}
In \textbf{Remark \ref{avio}}, we have already seen that ${\rm supp}^{\infty}$ is a spectral function, so we need only to show that is is a bijection.\\

Just note  that ${\rm supp}^{\infty} = \beta \circ \alpha' = \alpha \circ \beta'$, and since ${\rm supp}^{\infty}$ is a composition of bijections, it is a bijection.
\end{proof}

\end{document}